\documentclass[ijoc,nonblindrev]{informs3wp} % current default for manuscript submission

\OneAndAHalfSpacedXI % current default line spacing
\usepackage{endnotes}
\let\footnote=\endnote

\usepackage{natbib}
 \bibpunct[, ]{(}{)}{,}{a}{}{,}%
 \def\bibfont{\small}%
\TheoremsNumberedThrough % Preferred (Theorem 1, Lemma 1, Theorem 2)
\ECRepeatTheorems

\EquationsNumberedThrough% Default: (1), (2), ...
\usepackage[pdfpagemode={UseOutlines},bookmarks=true,bookmarksopen=true,
bookmarksopenlevel=0,bookmarksnumbered=true,hypertexnames=false,
colorlinks,linkcolor={blue},citecolor={black},urlcolor={red},
pdfstartview={FitV},unicode,breaklinks=true]{hyperref}
\hypersetup{urlcolor=blue, colorlinks=true}

\usepackage{multirow}% http://ctan.org/pkg/multirow
\usepackage{hhline}% http://ctan.org/pkg/hhline
\usepackage{footnote}

\newcommand{\transpose}{{\mbox{\tiny T}}}

\newcommand{\cA}{{\mathcal{A}}}
\newcommand{\cU}{{\mathcal{U}}}
\newcommand{\cG}{{\mathcal{G}}}
\newcommand{\cV}{{\mathcal{V}}}

\newcommand{\cP}{{\mathcal{P}}}
\newcommand{\cL}{{\mathcal{L}}}

\newcommand{\cK}{{\mathcal{K}}}
\newcommand{\cF}{{\mathcal{F}}}
\newcommand{\cN}{{\mathcal{N}}}
\newcommand{\cQ}{{\mathcal{Q}}}

\newcommand{\cT}{{\mathcal{T}}}
\newcommand{\cI}{{\mathcal{I}}}

\newcommand{\cH}{{\mathcal{H}}}
\newcommand{\cW}{{\mathcal{W}}}
\newcommand{\cO}{{\mathcal{O}}}

\newcommand{\bbe}{{\textbf{e}}}

\newcommand{\bX}{\textbf{X}}
\newcommand{\bY}{\textbf{Y}}
\newcommand{\bx}{\textbf{x}}
\newcommand{\by}{\textbf{y}}
\newcommand{\bs}{\textbf{s}}
\newcommand{\bd}{\textbf{d}}

\newcommand{\bu}{\textbf{u}}
\newcommand{\bc}{\textbf{c}}
\newcommand{\bU}{\textbf{U}}

\newcommand{\bp}{\textbf{p}}
\newcommand{\ba}{\textbf{a}}
\newcommand{\bz}{\textbf{z}}
\newcommand{\bb}{\textbf{b}}

\newcommand{\bw}{\textbf{w}}
\newcommand{\bap}{\pmb{\alpha}} 

\newcommand{\bld}{{\pmb{\lambda}}}
\newcommand{\dg}{\text{diag}}

\newcommand{\bbR}{\mathbb{R}}

\newcommand{\bbI}{\mathbb{I}}

\newtheorem{note}{Note}

\newcommand{\T}{\text{\tiny T}}
\usepackage{mathtools}

\newif\ifnotes\notestrue
\def\mtien#1{{{#1}}}

\def\ctien#1{{\color{red}{\textsc{[#1]}}}}
\def\htien#1{}
\def\here{{\color{red}{\textsc{[--------I AM HERE------------]}}}}

\newcommand{\RO}{\textsc{RO}}

\newcommand{\DET}{\textsc{DET}}

\def\me#1{\langle#1\rangle}

\begin{document}
%%%%%%%%%%%%%%%%

% Outcomment only when entries are known. Otherwise leave as is and
%default values will be used.
%\setcounter{page}{1}
%\VOLUME{00}%
%\NO{0}%
%\MONTH{Xxxxx}% (month or a similar seasonal id)
%\YEAR{0000}% e.g., 2005
%\FIRSTPAGE{000}%
%\LASTPAGE{000}%
%\SHORTYEAR{00}% shortened year (two-digit)
%\ISSUE{0000} %
%\LONGFIRSTPAGE{0001} %
%\DOI{10.1287/xxxx.0000.0000}%

% Author's names for the running heads
% Sample depending on the number of authors;
% \RUNAUTHOR{Jones}
% \RUNAUTHOR{Jones and Wilson}
% \RUNAUTHOR{Jones, Miller, and Wilson}
% \RUNAUTHOR{Jones et al.} % for four or more authors
% Enter authors following the given pattern:
\RUNAUTHOR{Mai T. and Jaillet P.}

% Title or shortened title suitable for running heads. Sample:
% \RUNTITLE{Bundling Information Goods of Decreasing Value}
% Enter the (shortened) title:
\RUNTITLE{Robust Product-line Pricing under Generalized Extreme Value Models}

% Full title. Sample:
% \TITLE{Bundling Information Goods of Decreasing Value}
% Enter the full title:
\TITLE{Robust Product-line Pricing under Generalized Extreme Value Models}

% Block of authors and their affiliations starts here:
% NOTE: Authors with same affiliation, if the order of authors allows,
%should be entered in ONE field, separated by a comma.
%\EMAIL field can be repeated if more than one author
\ARTICLEAUTHORS{%
\AUTHOR{Tien Mai}
\AFF{School of Information Systems, Singapore Management University, \EMAIL{atmai@smu.edu.sg}} %, \URL{}}
\AUTHOR{Patrick Jaillet}
\AFF{EECS, Massachusetts Institute of Technologies, \EMAIL{jaillet@mit.edu}}
% Enter all authors
} % end of the block

\ABSTRACT{%
We study robust versions of pricing problems where customers choose products according to a generalized extreme value (GEV) choice model, and the choice parameters are not known exactly but lie in an uncertainty set.
We show that, when the robust problem is unconstrained and the price sensitivity parameters are homogeneous, the robust optimal prices have a constant markup over products and we provide  formulas that allow to compute this constant markup  by bisection.  
We further show that, in the case that the price sensitivity parameters are only  homogeneous in each partition of the products, under the assumption that the choice probability generating function and the uncertainty set are partition-wise separable, a robust solution will have a constant markup in each subset, and this constant-markup vector can be found efficiently by convex optimization.
We provide numerical results  to illustrate the advantages of our robust approach in protecting  from bad scenarios. {Our results generally hold for convex and bounded uncertainty sets,} and for any arbitrary GEV model, including the multinomial logit, nested or cross-nested logit.
}%

% Sample
%\KEYWORDS{deterministic inventory theory; infinite linear programming duality;
%  existence of optimal policies; semi-Markov decision process; cyclic schedule}

% Fill in data. If unknown, outcomment the field
%\HISTORY{This paper wasfirst submitted on April 12, 1922 and has been with the authors for 83 years for 65 revisions.}

\maketitle
%%%%%%%%%%%%%%%%%%%%%%%%%%%%%%%%%%%%%%%%%%%%%%%%%%%%%%%%%%%%%%%%%%%%%%

% Samples of sectioning (and labeling) in OPRE
% NOTE: (1) \section and \subsection do NOT end with a period
% (2) \subsubsection and lower need end punctuation
% (3) capitalization is as shown (title style).
%
%\section{Introduction.}\label{intro} %%1.
%\subsection{Duality and the Classical EOQ Problem.}\label{class-EOQ} %% 1.1.
%\subsection{Outline.}\label{outline1} %% 1.2.
%\subsubsection{Cyclic Schedules for the General Deterministic SMDP.}
%  \label{cyclic-schedules} %% 1.2.1
%\section{Problem Description.}\label{problemdescription} %% 2.

% Text of your paper here

%%%%%%%%%%

\section{Introduction}
In revenue management, pricing  is an important problem that refers to the selection of prices for a set of products in order to maximize an expected revenue. This is motivated by the fact that prices are key features that may significantly  affect demand for products. The literature of multi-product pricing has seen a large number of papers focusing on how to set prices when customers purchase products according to a discrete choice model \citep[e.g.][]{Talluri2004RM,Gallego2014multiproduct,zhang2018multiproduct}. 
To the best of our knowledge, prior work all assumes that the parameters of the choice models are known in advance or can be estimated exactly from data. Thus, the corresponding pricing optimization models are built based on pre-determined parameters and ignore any uncertainty in case the parameters are estimated. Nevertheless, in practice, the parameter estimates may vary significantly for different customer types or in different purchasing periods of the year. Thus, ignoring such uncertainties may  lead to bad pricing decisions.
To deal with the uncertainty issue,  one may consider a stochastic approach, i.e., a model aiming at maximizing an average expected revenue over a finite number of scenarios of the choice parameters.
 This would require of course a trusted assumption and/or a solid optimization of these parameters in each of these scenarios.  Moreover, such a stochastic optimization model would be computationally difficult to handle, as the objective function does not have nice properties to derive tractable solutions as in the deterministic case; e.g., a stochastic objective function would be non-unimodal and non-concave when defined in terms of purchase probabilities  \citep{Li2018mixed_pricing}.

{In this paper, we formulate  and solve pricing optimization problems under uncertainty in a robust manner. 
That is, we assume customers' behavior is driven by any choice model in the Generalized Extreme Value (GEV) family such as the Multinomial Logit (MNL) or nested logit model, and the parameters of the choice model are not known exactly but belong to an uncertainty set.
The goal here is to maximize the worst-case expected revenue when the choice parameters vary in their support set. We consider 
problems where the price sensitivity parameters (PSP) are homogeneous or partition-wise homogeneous, i.e., the set of products can be separated into disjoint subsets and the PSP are the same in each subset but can be different over subsets.
 For the latter, we assume that the choice probability generating function \citep{Fosgerau2013GPGF} has a separable structure and the uncertainty set is partition-wise separable.   
We also look at expected-sale requirements in pricing decisions and  argue that the model with expected-sale constraints is not appropriate in our robust setting. Therefore, we propose an alternative formulation by adding a penalty term to the objective function for violated expected-sale constraints. 
We are able to show that the models can then be solved in  a tractable way. Our results 
generally hold for any convex and bounded uncertainty set, and  for any choice model in the GEV family. }

From now on, when saying ``a GEV model'', we refer to any choice model in the GEV family. Each GEV model can be represented by a choice probability generating function (CPGF) $G(\cdot)$ (see our detailed definition in the next section). To relax the homogeneity of the PSP, we need to assume that the CPGF  has a separable structure, which means that 
$G(\cdot)$ can  be written as a sum of sub-CPGFs, each corresponding to a subset of products. 

\noindent

%\ctien{Revise the contributions .... }

%\ctien{Now we do not restrict ourselves to rect. uncertainty sets}

%\ctien{Better motivate the robust model: the stochastic one is not tractable..., and the distribution of the choice parameters is typically hard to estimate ...}

\textbf{Our contributions:} We consider robust versions of the standard pricing optimization problem under GEV models. The setting here is to assume that the parameters of the choice model are not known with certainty and the aim is to find optimal prices associated with  products, which maximize the worst-case expected revenue when the choice parameters vary in an uncertainty set. 
For the unconstrained problem with homogeneous PSP, we show that if the uncertainty set is convex and compact, the robust optimal prices have a constant markup with respect to the products costs, i.e., the robust optimal price of a product is equal to its unit cost plus a constant that is the same over all products. We also provide formulas that allow efficient computation of that constant markup by binary search. This finding  generalizes the results for the deterministic unconstrained problem with homogeneous PSP considered  in  \cite{zhang2018multiproduct}.
We also provide comparative insights showing how the robust optimal revenue and the robust optimal constant markups change as functions of the uncertainty level (i.e., the size of the uncertainty set).

\mtien{For the pricing problem with \textit{non-homogeneous} PSP, %similarly to previous studies \citep{zhang2018multiproduct}, 
we assume that the  CPGF is partition-wise separable and in each partition, the PSP are homogeneous.
Moreover, the uncertainty set is also assumed to be partition-wise separable. 
We show that the robust problem can be converted equivalently into a \textit{reduced optimization problem}, which  can be conveniently solved by convex optimization. 
As a result, the robust optimal prices have partition-wise constant markups, i.e., in each partition, the robust optimal prices have a constant markup with respect to their costs,  and these constant markups can be obtained by  convex optimization.
We also provide comparative insights for the robust optimal prices and solutions when the size of the uncertainty set varies.}

\mtien{For both cases (i.e., homogeneous PSP and partition-wise homogeneous PSP), we further show that the robust optimal solutions form saddle points of the robust problems, leading to an equality between the objective functions of the \textit{max-min} problem and its \textit{min-max} counterpart. }

Previous studies \citep{zhang2018multiproduct,Song2007demand,Zhang2013assessing}  have been looking at constraints on the expected sales,  as motivated by applications with inventory considerations \citep{Gallego1997multiproductPricing}. In this context, the aim  is to select prices that maximize the expected revenue while requiring that the expected  sales of products lie in a convex set. The advantage of such constraints is that the pricing problem can be reformulated equivalently as a convex program where the decision variables are the purchase probabilities. {However, the final decision is a vector of prices and there may be no fixed prices under which the resulting purchase probabilities always satisfy the expected sale constraints when the choice parameters vary. For this reason, the use of the constrained formulation is not appropriate}
in our robust setting. 
Thus, we propose an alternative formulation in which, instead of requiring that the expected sale constraints be satisfied, we add a penalty cost to the objective function for violated constraints. Our formulation, called pricing with over-expected-sale penalties, is more general than the constrained formulation, in the sense that 
%we show that 
if the penalty parameters increase to infinity,  then the corresponding optimal solutions will converge to those from the constrained problem, and with \textit{zero} penalty parameters, the pricing problem becomes the unconstrained one. 
%{Since the robust version under over-expected-sale penalties does not seek a price solution that satisfies the expected sale constraints under choice parameter uncertainty, it is more appropriate to use than the constrained version.} 
We show that if the  CPGF and the uncertainty set are partition-wise separable, then the robust problem can be converted into a reduced optimization problem, which can be  conveniently solved by convex optimization. % We also investigate the convergence of the optimal values, in both robust and deterministic cases, when the penalty parameters go to infinity.   

\mtien{In summary, we show that the robust versions of the pricing problem with homogeneous PSP and partition-wise PSP, with and without over-expected-sale penalties, can be solved in tractable ways by bisection and convex optimization. Our results generally holds for convex and compact uncertainty sets, and for any choice model in the GEV family. In Table \ref{tab:summart-DE-RO} below we give a summary and comparison of the solution methods used to solve the robust pricing problems and their deterministic counterparts, under different settings. The solution methods proposed in this paper are highlighted in bold.} 

\begin{table}[htb]
\centering
\begin{tabular}{l|l|l}
Settings & Deterministic pricing & \textbf{Robust pricing}\\ 
\hline
Unconstrained and homogeneous PSP& Closed-form solutions & \begin{tabular}[c]{@{}l@{}}\textbf{Bisection and }\\\textbf{convex optimization}\end{tabular}  \\\hline
\begin{tabular}[c]{@{}l@{}}Unconstrained and \\partition-wise homogeneous PSP\end{tabular} & Bisection & \textbf{Convex optimization}\\\hline
Expected-sale constraints& Convex optimization& {Not appropriate~}\\\hline
Over-expected-sale penalties& Convex optimization& \textbf{Convex optimization} 
\end{tabular}
\caption{\mtien{Solution  methods for deterministic pricing and robust pricing problems under different settings.}}
\label{tab:summart-DE-RO}
\end{table}

\noindent
\textbf{Literature review:} The GEV family includes most of the parametric discrete choice models  in the demand modeling and operations research literatures. The simplest and most popular member is the MNL \citep{Mcfadden1978modeling,Mcfadden1980econometric} and it is well-known that the MNL model retains the independence from irrelevant alternatives (IIA) property, which does not hold in  many contexts. There are a number of GEV models that relax this property and  provide flexibility in modeling the correlation between alternatives, for example, the nested logit model \citep{Ben1985discrete,Ben1973structure}, the cross-nested logit \citep{Vovsha1998link}, the generalized nested logit \citep{Wen2001generalized}, the paired combinatorial logit \citep{koppelman2000paired}, the ordered generalized extreme value \citep{Small1987discrete}, the specialized compound generalized extreme value models \citep{Bhat1998accommodating,whelan2002flexible} and network-based GEV \citep{Daly2006general,Mai2017dynamic} models.  
\cite{Fosgerau2013GPGF} show that the cross-nested logit model and its generalized version (i.e. network-based GEV) are fully flexible in the sense that they can approximate arbitrarily close any random utility maximization model.  
Beside the GEV family, it is worth noting that the mixed logit model \citep{Mcfadden2000mixed} is also popular due to its flexibility in capturing utility correlation. There is a fundamental trade-off between  the flexibility and the generality of the choice models and the complexity of their estimation andapplication in operational problems. 
For the case of GEV models, even being flexible in modeling choice behavior, the resulting  operational problems (e.g., product assortment or pricing) are often nonlinear and non-convex, leading to difficulties solving them in practice.

There is a large amount of research on unconstrained pricing under different discrete choice models. For example, \cite{Hopp2005}  and \cite{Dong2009dynamic} consider  the  pricing problem under the MNL model, \cite{Li2011pricing} consider the  nested logit model, \cite{li2015d_nested} consider the pricing problem under the paired combinatorial logit model, and \cite{zhang2018multiproduct} consider the pricing problem under any choice model in the GEV family. Under the assumption that the PSP are the same over product, these authors show that the prices have a constant markup with respect to the product costs and provide formulas to explicitly computed this constant markup.  

There are some papers trying to get over the assumption that the PSP are homogeneous over products. \cite{Li2011pricing} study the pricing problem under the nested logit model and assume that the PSP are homogeneous only in each nest and can be different over nests. They then show that the PSP in each nest have a constant markup. \cite{zhang2018multiproduct} generalize these results by considering the pricing problem under GEV models, in which the CPGF is partition-wise separable and the PSP are assumed to be homogeneous in each partition. The authors also show that, in this case, the optimal prices have a constant markup in each partition. 

There are also publications considering the pricing problem with arbitrary PSP. \cite{Gallego2014dynamicPricing} show that the pricing optimization problem under the nested logit model can have multiple local optimal solutions if the PSP are arbitrarily heterogeneous and provide sufficient conditions to ensure unimodality of the expected revenue function.  \cite{li2015d_nested} and \cite{Huh2015pricing} consider the pricing problem under the  $d$-nested and paired combinatorial logit models and also provide sufficient conditions on the PSP to ensure unimodality of the expected revenue function. 

The constrained pricing problem where the prices are required to lie in a feasible set is difficult to solve as the expected revenue function is nonlinear and non-concave in the prices. Motivated by applications with inventory considerations \citep{Gallego1997multiproductPricing} and the observation that the expected revenue function is concave in the purchase probabilities, researchers have consider the pricing problem with constraints on the expected sales.
For example, \cite{Song2007demand,Zhang2013assessing} consider the pricing problem under the MNL model and show that the expected revenue is concave in the purchase probabilities if the PSP are homogeneous. \cite{Keller2013} consider the pricing problem under the MNL and nested logit models and show that the expected revenue function is concave in the purchase probabilities under the MNL and arbitrary PSP, and establish sufficient conditions on the PSP  to ensure that the expected revenue under the nested logit model is concave. \cite{zhang2018multiproduct} also generalizes all these results by showing that, under any GEV model, if the PSP are homogeneous or partition-wise homogeneous, then the expected revenue is concave in purchasing probabilities, making the pricing problem with expected sale constraints tractable. 

All above publications assume that the parameters of the choice model is given in advance and ignore any uncertainty associated with such parameters in the pricing problem. However, the choice parameters typically need to be inferred from data and uncertainties may occur, for instance, due to the heterogeneity of the market. 
In this work, we explicitly take into consider this issue by considering robust versions of the unconstrained and constrained pricing problems, with homogeneous and partition-wise homogeneous PSP. Our results directly generalize the results for deterministic pricing  from \cite{zhang2018multiproduct}, which already covers most of the pricing optimization studies in the literature. 

Our work is concerned with robust solutions for the pricing problem under uncertainty, so it is directly related to the concept of robust optimization, an important research area in operations research which has received a growing attention over the past two decades. Robust optimization is motivated by the fact that many real-world decision problems arising in engineering and
management science have uncertain parameters due
to limited data or noisy measurements. The literature on robust optimization includes a larger number of excellent studies \citep[see][for instance]{Ben1998robustconvex,Ben2000robust,Ben2006extending}. Most of the studies in the literature of robust optimization  focus on linear, piece-wise linear or convex objective functions.  In our context, the expected revenue is nonlinear and non-convex/non-concave in the prices, implying that existing robust optimization results do not apply (except the part where we consider the constrained pricing problem under uncertain expected-sale constraints in Section \ref{sec:robust-con}), and making our robust problem challenging to solve in a tractable way. It is worth noting that our work is relevant to \cite{Rusmevichientong2012robust}
where the authors consider  robust versions of the assortment planing problem. The main difference is that the decision variables in \cite{Rusmevichientong2012robust} are discrete  (i.e., a set of products).

\noindent
\textbf{Paper outline:} We organize the paper as follows.
In Section \ref{sec:der-pricing}, we present the deterministic pricing problem under GEV models and recall some results from  previous work.  In Section \ref{sec:robust-homo} and \ref{sec:uncon-hete}, we present our results for the robust pricing problem under homogeneous PSP
and partition-wise homogeneous PSP.  In  Section \ref{sec:expertiments} we provide some experimental results and in Section \ref{sec:conclude} we conclude. In the appendix,  Section \ref{apd:sec-proofs} provides   detailed proofs for our main claims and  Section \ref{sec:robust-penalties}  investigates the robust pricing problem with over-expected-sale penalties.

\noindent
\textbf{Notation:}
Boldface characters represent matrices (or vectors), and $a_i$ denotes the $i$-th element of vector $\ba$. We use $[m]$, for any $m\in \mathbb{N}$, to denote the set $\{1,\ldots,m\}$. For any vector $\bb$ with all equal elements, we use $\me{\bb}$ to denote the value of one element of the vector. Given two vectors of the same size $\ba,\bb \in \bbR^m$, $\ba \succeq \bb$ is equivalent to $\ba-\bb \in\bbR_+^m$, and  $\ba \preceq \bb$ is equivalent to $\bb \succeq \ba$. 

\section{Background: Deterministic Pricing under Generalized Extreme Value Models}
\label{sec:der-pricing}
We denote by $\cV =  \{1,\ldots,m\}$ the set of $m$ available products. There is a \textit{non-purchase item}  indexed by 0, so the set of all possible \textit{products} is $\cV \cup \{0\}$. We also denote by $x_i$ and $c_i$ the price and the cost of product $i$, respectively. The random utility maximization (RUM) framework \citep{Mcfadden1978modeling} is the most popular approach to model discrete choice behavior. Under this framework, each product $i\in\cV$ is assigned with a random utility $U_i$ and the additive RUM framework \citep{Fosgerau2013GPGF,Mcfadden1978modeling} assumes that each random utility can be expressed as a sum of two part $U_i =  u_i+\varepsilon_i$, where the term $u_i$ is deterministic and can include values representing characteristics of the product, and the term $\varepsilon_i$ is unknown to the analyst. The RUM principle then assume that  the selections are made by maximizing these utilities and the probability that a product $i$ (including the non-purchase item) is selected can be computed  as $P(U_i \geq U_j,\ \forall j\in\cV\cup\{0\} )$.

In our context, we are interested in the effect of the prices on the expected revenue. So we assume that the deterministic terms  
$u_i$, $\forall i\in\cV$,  can be expressed as $u_i = a_i-b_i x_i$, where $b_i$ is the PSP associated with product $i$ and $a_i$ can include other information that may affect customer's demand such as the brand, size or color of the items. 
These values can be obtained by fitting the choice model with observation data. 
%Note that here we assume that the utilities $u_i$ depend linearly on the prices, which is a popular assumption in most of the existing pricing studies. %Nonlinear formulation would make the pricing problem much more difficult to deal with, but would be interesting to look at in future research.

A GEV model can be represented by a choice probability generating function (CPGF) $G(\bY)$, where $\bY$ is a vector of size $m$ with entries $Y_i = e^{u_i}$, for all $i\in\cV$.
Given $i_1,\ldots,i_k \in [m]$, let $\partial G_{i_1,\ldots,i_k}(\bY)$ be the mixed partial derivatives of $G$ with respect to $Y_{i_1},\ldots,Y_{i_k}$.  
It is well-known that 
%To be consistent with the RUM principle, 
the CPGF  $G(\cdot)$ and the mixed partial derivatives have the
the following properties \citep{Mcfadden1978modeling,Ben1985discrete}.
\begin{remark}[\textbf{Properties of GEV-CPGF}]\label{propert:GEV-CPGF}
{\it A GEV-CPGF $G(\bY)$ has the following properties.
\begin{itemize}
 \item[(i)] $G(\bY) \geq 0,\ \forall \bY\in \bbR^m$,
 \item[(ii)] $G$ is homogeneous of degree one, i.e., $G(\lambda \bY) = \lambda G(\bY)$
 \item[(iii)] $G(\bY)\rightarrow \infty$ if $Y_i\rightarrow \infty$
 \item[(iv)]  Given $i_1,\ldots,i_k \in [m]$ distinct from each other,
  $\partial G_{i_1,\ldots,i_k}(\bY)>0$ 
 if $k$ is odd, and $\leq$ if $k$ is even
 \item[(v)] $G(\bY) = \sum_{i\in\cV} Y_i\partial G_i(\bY)$
 \item[(vi)] $\sum_{j\in\cV} Y_j\partial G_{ij} (\bY) = 0$, $\forall i\in \cV$.
\end{itemize}
%where $\partial G_i(\bY) = \partial G(\bY)/\partial Y_i$
}
\end{remark}
\mtien{Here we note that  (i)-(iv) are basic properties of the CPGF  to ensure that the choice model is consistent with the RUM principle \citep{Mcfadden1980econometric}. Properties (v) and (vi) are direct results from the  homogeneity property \citep{zhang2018multiproduct}. }   

Under a GEV model specified by a CPGF $G$, given any vector $\bY\in\bbR^m$,  the choice probability of product $i\in\cV$ is given by
\[
 P_i(\bY|G) = \frac{Y_i\partial G_i(\bY)}{1+G(\bY)}.
\]
Note that the above formulation also implies that the choice probability of the \textit{non-purchase} item  is $P_0(\bY|G)  = 1/(1+G(\bY))$. The GEV becomes the MNL model if $G(\bY) = \sum_{i=1}^m Y_i$, and it becomes the nested logit model if 
$G(\bY) = \sum_{n\in \cN} \left(\sum_{i\in C_n} (\sigma_{in}Y_i)^{\mu_n} \right)^{\mu/\mu_n}$, where $\cN$ is the set of nests,  $C_n$ is the set of items in  nest $n$ and $\sigma_{in},\mu>0,\mu_n>0$ are the parameters of the nested logit model. In the generalized version of the nested logit model proposed by \cite{Daly2006general},  called the network GEV, the corresponding CPGF can  be computed recursively based on a rooted and cycle-free graph representing the correlation structure of the items.

Under a GEV  model specified by a CPGF $G(\cdot)$, the deterministic version of the pricing problem is stated as  
\begin{equation}\label{prob:pricing-deterministic}\tag{P1}
 \max_{\bx \in \bbR^m} \left\{R(\bx) = \sum_{i=1}^m (x_i-c_i) P_i(\bY(\bx,\ba,\bb)|G)\right\},
\end{equation}
where $\bY(\bx,\ba,\bb) \in \bbR^m$ with entries $Y_i(\bx,\ba,\bb) = \exp(a_i-b_ix_i)$. The expected revenue $R(\bx)$ becomes more difficult to handle  as the GEV model becomes more complicated. By leveraging the properties of GEV models stated in Remark \ref{propert:GEV-CPGF}, \cite{zhang2018multiproduct} manage to show that if the PSP are homogeneous, i.e., $b_i = b_j$ for all $i,j\in\cV$ and if  $\bx^*$ is an optimal  solution to \eqref{prob:pricing-deterministic}, then
\begin{equation}\label{eq:optimal-solution-deterprob} 
x^*_i-c_i = \frac{1}{\me{\bb}}+ R(\bx^*),\forall i\in\cV \ \text{and }R(\bx^*) = \frac{W(\gamma e^{-1})}{\me{\bb}} 
\end{equation}
where $\gamma = G(Y_1(c_1),\ldots,Y_m(c_m))$ and $W(\cdot)$ is the Lambert-W function.
The results in \eqref{eq:optimal-solution-deterprob} indeed imply that a constant markup solution is optimal to  \eqref{prob:pricing-deterministic} and this constant markup can be computed explicitly. 
Moreover, if the PSP are partition-wise homogeneous and $G$ is separable, then \cite{zhang2018multiproduct} show that the optimal prices have a constant markup in each partition. These results also provide an explicit way to compute optimal prices for the pricing problem under the MNL with arbitrary PSP.
\cite{zhang2018multiproduct} also show that the expected revenue function is concave in the purchasing probabilities under any GEV model, making the pricing problem with expected sale constraints tractable. 

\section{Robust Pricing under Homogeneous Price Sensitivity Parameters}
\label{sec:robust-homo}
In this section, we study a robust version of the unconstrained pricing problem, under the setting that the choice parameters $(\ba,\bb)$ are not known exactly but belong to an uncertainty set. We focus here on  the case of  homogeneous PSP.
In our robust model, we aim at maximizing the worst-case expected revenue over all parameters in the uncertainty set. The  robust unconstrained pricing problem can be formulated as 
\begin{equation}\label{prob:RO}\tag{RO}
 \max_{\bx \in \bbR^m} \left\{g(x) =  \min_{(\ba,\bb)\in \cA}  \sum_{i=1}^m (x_i-c_i) P_i(\bY(\bx,\ba,\bb)|G),  \right\},
\end{equation}
where $\cA$ is the uncertainty set of the parameters $(\ba,\bb)$. 
We denote $\Phi(\bx,\ba,\bb) =  \sum_{i=1}^m (x_i-c_i) P_i(\bY(\bx,\ba,\bb)|G)$ for notational simplicity. 
We  assume that $\cA$ is \textit{convex  and bounded}. The convexity and boundedness assumptions are useful later in the section, as we need to show that, under a constant-markup style vector of prices, 
the objective function of the adversary's problem is convex on $\cA$, which in turn helps identify a saddle point of the robust problem. The boundedness assumption is realistic in the context, as  the choice parameters are often inferred from data and it is expected that they are finite.  We also assume that the PSP are positive, i.e., ${\bb}>0$ for any $(\ba,\bb) \in\cA$, which is conventional from a behavior point of view.

When the PSP are the same over all the products, we will show that the robust optimal prices
have a constant markup and this constant markup can be computed efficiently by binary search.
The idea is motivated by the observation that if we consider the \textit{min-max} counterpart of the robust problem
$
\min_{(a,b)\in\cA}\; \max_{\bx\in\bbR^{m} } \Big\{\Phi(\bx,\ba,\bb)\Big\},
$
then  we know that the adversary problem always yields a constant-markup optimal solution for any fixed choice parameters $(\ba,\bb)$  \citep{zhang2018multiproduct}. So, the \textit{min-max} counterpart is equivalent to
\begin{equation}
\label{eq:uncon-minmax-prob}
\min_{(a,b)\in\cA}\; \max_{\bx\in \bX}\Big\{\Phi(\bx,\ba,\bb)\Big\},
\end{equation}
where $\bX$ is the set of constant-markup solutions, i.e., $\bX = \{\bx\in\bbR^m|\ x_i-c_i = x_j-c_j, \forall i,j\in[m]\}$. This suggests that if there is a saddle point of the \textit{max-min} problem \eqref{prob:RO}, then it should have a constant-markup form.
%In the rest of the section, we will show that such a saddle point exists and can be computed efficiently by bisection. 

To prove the result, we will consider the robust unconstrained pricing problem with constant-markup prices, i.e., we only look at prices $\bx \in \bX$. Then we show that there exist constant-markup prices $\bx^*$ such that if $(\ba^*,\bb^*)$ is an optimal solution to the adversary's problem under prices $\bx^*$, then $\bx^*$ is also optimal to the deterministic unconstrained problem with choice parameters $(\ba^*,\bb^*)$. In other words, $(\bx^*,\ba^*,\bb^*)$ is a saddle point of \eqref{eq:uncon-minmax-prob} and $\bx^*$ is also an optimal solution to the robust problem. 

Given constant-markup prices $\bx\in\bX$ and choice parameters $(\ba, \bb)\in\cA$, the expected revenue becomes
\[
\begin{aligned}
\sum_{i=1}^m (x_i-c_i) P_i(\bY(\bx,\ba,\bb)|G) &=  \frac{z \sum_{i\in\cV} Y_i(\bx,\ba,\bb)\partial G_i(\bY(\bx,\ba,\bb))}{1+G(\bY(\bx,\ba,\bb))} \\
&=z\left(1 - \frac{1}{1+G(\bY(\bx,\ba,\bb))}\right),
\end{aligned}
\]
where $z = x_i-c_i$, $\forall i\in \cV$ and $\bY$ is a vector with entries $Y_i = \exp(a_i -b_i (z+c_i))$ for all $i\in \cV$. 
For the sake of simplicity, from now on we will write $\bY$ instead of $\bY(\bx,\ba,\bb)$.
The expected revenue is a function of $z$ and $(\ba,\bb)$, and if $(\ba^*(z), \bb^*(z))$ is an optimal solution to the adversary's problem, then we  also have 
\begin{equation}\label{prob:m1}
 	(\ba^*(z),\bb^*(z))  = \underset{\ba,\bb \in \cA}{\text{argmin}}\qquad G(\bY|z,\ba,\bb), 
\end{equation}
where $G(\bY|z,\ba,\bb) = G(Y_1,\ldots,Y_m)$ with $Y_i = e^{a_i-b_i(z+c_i)}$. 
In Proposition \ref{prop:G-strictly-convex} below, we first  show that $G(\bY|z,\ba,\bb)$ is strictly convex in $(\ba,\bb)$. As a result,  $(\ba^*(z),\bb^*(z))$ is always uniquely determined. This result is important to identify a saddle point of the robust problem. 
\begin{proposition}\label{prop:G-strictly-convex}
Given any $z\in \bbR_+$, $G(\bY|z,\ba,\bb)$ is strictly convex on $\cA$, Problem \ref{prob:m1} always has a unique solution, and $(\ba^*(z),\bb^*(z))$ determined in \eqref{prob:m1} is continuous in $z\in \bbR_+$. 
\end{proposition}
The proof is given in Appendix \ref{proof:G-strictly-convex}. The proposition plays an important role in our main claim, as in the theorem below we will show that a solution to the robust problem can be found by solving a 1-dimensional fixed-point problem. The continuity of $(\ba^*(z),\bb^*(z))$ guarantees that this fixed-point problem always has a solution that can be found efficiently by bisection. 
%Theorem \ref{theor:RO-GEV-homo} below is the main result of the section. It shows that the robust optimal prices have a constant markup and this constant  can be found efficiently by binary search. 
% The proof of the theorem is quite obvious, as we already show that there isa saddle point of the robust problem that has the constant-markup style (Lemma \ref{lemma:z*-exists}).  
\begin{theorem}[Constant markup is optimal to the robust problem)]\label{theor:RO-GEV-homo}
There always exists a unique solution $z^{*}\in\bbR$ to the fixed point problem
\begin{equation}\label{eq:z-equation}
z = \frac{1+ W(G(\bY|0,\ba^*(z),\bb^*(z))e^{-1})}{\me{\bb^*(z)}},
\end{equation}
where $W(\cdot)$ is the Lambert-W function and the constant-markup prices $\bx^{*}$ defined  as  $x^{*}_i = z^{*} + c_i ,\ \forall i\in[m]$, is the unique  robust solution of the robust problem \eqref{prob:RO}. Moreover, $(\bx^*,\ba^*(z^*), \bb^*(z^*))$ is a saddle point of \eqref{prob:RO} and the minimax equality holds, i.e., $$\max_{\bx\in\bbR^m}\min_{(\ba,\bb)\in\cA}\Phi(\bx,\ba,\bb) = \min_{(\ba,\bb)\in\cA}\max_{\bx\in\bbR^m}\Phi(\bx,\ba,\bb).$$
\end{theorem}
We highlight two important claims from Theorem \ref{theor:RO-GEV-homo}. First, the robustness preserves the constant-markup property of the solutions to  the deterministic pricing problem, and second, the \textit{minimax equality} holds. We note that the \textit{minimax equality} is not straightforward to see at first sight, as the objective function $\Phi(\bx,\ba,\bb)$ is not (quasi) concave in $\bx$ nor convex in $(\ba,\bb)$. 

We will make use of Lemmas \ref{lemma:G-bounded} -\ref{lemma:z*-exists} below to prove the theorem. In Lemma \ref{lemma:G-bounded} we show that function $G(\bY|z,\ba, \bb)$ is always bounded. Together with the results established in Proposition \ref{prop:G-strictly-convex}, it then becomes clear that there always exists a fixed point solution to \eqref{eq:z-equation}. As a result, if $z^*$ is a solution \eqref{eq:z-equation}, then it will form  optimal constant-markup prices for the robust problem. Now, let us go into details of the lemmas and proofs. 

%%%%%%%%%%%%%%%%%%%%%%%%%%
%%%%%%%%%%%%%%%%%%%%%%%%%%
%In the next lemma, we show that, given any $z\in \bbR_+$, if  the uncertainty set $\cA$ is bounded, then the function $G(\bY|z,{\ba},{\bb})$ is also bounded for all parameters $(\ba,\bb)\in \cA$. The lemma allows us to determine an finite interval where we can search the robust optimal constant markup. 
\begin{lemma}\label{lemma:G-bounded} 
If there are  $(\underline{\ba}, \underline{\bb}), (\overline{\ba}, \overline{\bb}) \in \bbR^{2m}$ such that $\underline{\ba} \leq \ba \leq \overline{\ba}$ and  $\underline{\bb} \leq \bb\leq  \overline{\bb}$ for all $(\ba,\bb) \in\cA$, then 
$
G(\bY|z,\underline{\ba},\overline{\bb}) \leq  G(\bY|z,\ba,\bb) \leq G(\bY|z,\overline{\ba},\underline{\bb}), \ \forall z\in \bbR_+,\ (\ba,\bb) \in\cA.
$
\end{lemma}
The proof can be done quite easily using the properties of function $G(\cdot)$ and we refer the reader to Appendix \ref{apdx:proof-Lemma-G-bounded} for details. We are now ready to show that there is a solution to the fixed point problem \eqref{eq:z-equation}.
In Lemma \ref{lemma:z*-exists} below we show this by making use of the continuity of $\ba^*(z),\bb^*(z)$ (showed above)  and the boundedness assumption on $\cA$ to identify an interval where we can find $z^*$. Without this assumption, one can simply choose 0 as a lower bound, as $f(0)$ is always less than 0. However, to identify an upper bound, one needs some limits from the uncertainty set. This is because even in the deterministic case, if the choice parameters  $\bb$ approach zero, or $\ba$ increase to infinity,  then the optimal constant markup will go to infinity (see Equation \ref{eq:optimal-solution-deterprob}).

%%%%%%%%%%%%%%%%%%%%%%%%%
%%%%%%%%%%%%%%%%%%%%%%%%%

\begin{lemma}
\label{lemma:z*-exists}
For any $i\in \cV$, there exists $z^*\in \bbR_+$ such that 
\[
z^* = \frac{1+ W(\tau(z^*))}{\me{\bb^*(z^*)}} \in \left[\underline{Z}^0,\overline{Z}^0 \right] 
\]
where
\[
\begin{aligned}
 \underline{Z}^0 &= \frac{1+W(G(\bY|0,\underline{\ba},\overline{\bb})e^{-1})}{\me{\overline{\bb}}} \\
 \overline{Z}^0 &=  \frac{1+W(G(\bY|0,\overline{\ba},\underline{\bb})e^{-1})}{\me{\underline{\bb}}} \\
 \tau(z^*) &= G(\bY|0,\ba^*(z^*),\bb^*(z^*))e^{-1}
\end{aligned}
\]
and  $\underline{\ba},\overline{\ba},\underline{\bb}, \overline{\bb} \in\bbR^m$ such that  $\underline{\ba} \leq \ba \leq \overline{\ba}$ and  $\underline{\bb} \leq \bb\leq  \overline{\bb}$ for all $(\ba,\bb) \in\cA$,  
 and $W(\cdot)$ is the is the Lambert-W function.
\end{lemma}
\proof{Proof:}
Let
$f(z) = z - ({1+W(\tau(z))}))/{\me{\bb^*(z)}}.
$ 
From Lemma \ref{lemma:G-bounded}, we have 
\[
\underline{Z}^0 \leq \frac{1+W(\tau(z))}{\me{\bb^*(z)}} \leq \overline{Z}^0,\ \forall z\in\bbR_+.
\]
Which means
\[
f(\underline{Z}^0) \leq 0; \ f(\overline{Z}^0) \geq 0  
\]
Since $f(z)$ is continuous in $z$ (Proposition  \ref{prop:G-strictly-convex}), equation $f(z) = 0$ always has a solution in the interval $\left[\underline{Z}^0,\overline{Z}^0 \right]$. 
\endproof
We are now ready for the proof of Theorem \ref{theor:RO-GEV-homo}. Basically, we will show that a $z^*$ determined in Lemma \ref{lemma:z*-exists} and $(\ba^*(z^*),\bb^*(z^*)$ will form a saddle point to the robust problem. 

\proof{Proof of Theorem \ref{theor:RO-GEV-homo}:}
We know that there always exists $z^{*}$ being a fixed point solution to \eqref{eq:z-equation} (Lemma \ref{lemma:z*-exists}).
Given $\bx^{*}$ and $z^{*}$, we
first remark that $(\ba^*(z^{*}), \bb^*(z^{*}))$ is also the unique solution of the adversary's problem
\[
\underset{(\ba,\bb) \in \cA}{\text{argmin}}\qquad \Phi(\bx^{*},\ba,\bb)
\]
Moreover,  according to the way $\bx^{*}$ is computed and Theorem 3.1 of \cite{zhang2018multiproduct}, $\bx^{*}$ is optimal to the following  problem 
\[
\max_{\bx\in\bbR^m}\quad \left\{ \Phi(\bx,\ba^*(z^{*}), \bb^*(z^{*}) ) \right\}.
\]
This leads to the fact that $(\bx^{*},\ba^*(z^{*}), \bb^*(z^{*}) )$ is a saddle point to the robust \textit{max-min} problem \eqref{prob:RO}. In other words, $\bx^{*}$  is an optimal solution to the robust problem.

Note that the deterministic version of the unconstrained pricing problem always has a unique solution, which is a constant markup one. So,
for any $\bx \neq \bx^{*}$ we have
\begin{align}
g(\bx^{*}) &= \Phi(\bx^{*},\ba^*(z^{*}), \bb^*(z^{*}) )\nonumber \\
&> \Phi(\bx,\ba^*(z^{*}), \bb^*(z^{*}) )\nonumber \\
&\geq g(\bx).
\end{align}
Thus, there is only one solution to the robust pricing problem \eqref{prob:RO} and there is only one solution to the equation \eqref{eq:z-equation}, as required.
Since there is a saddle point to the \textit{max-min} problem \eqref{prob:RO}, the minimax equality holds, i.e., $\min_{(\ba,\bb)\in\cA} \max_{\bx\in\bbR^{m} } \Big\{\Phi(\bx,\ba,\bb)\Big\} =  \max_{\bx\in\bbR^{m} } \min_{(\ba,\bb)\in\cA} \Big\{\Phi(\bx,\ba,\bb)\Big\}$. \mtien{Note that the existence of a saddle point  directly  implies the minimax equality (a.k.a minimax equality),
but the opposite  does not always hold.} 
%In fact, there is another way to prove the minimax equality using \textit{Sion's minimax theorem} \citep{Sion1958general}, as shown in Appendix \ref{apdx:another-proof-for-the-minimax-equality}.}  
\endproof

Theorem \ref{theor:RO-GEV-homo} implies that a solution to the  robust problem can be found by solving the equation
\begin{equation}\label{prob:RO-EQ}
f(z) = z- \frac{1+ W(G(\bY|0,\ba^*(z),\bb^*(z))e^{-1})}{\me{\bb^*(z)}}= 0,
\end{equation}
in the interval $\left[\underline{Z}^0,\overline{Z}^0 \right]$, in which $\underline{Z}^0,\overline{Z}^0$ are defined in Lemma \ref{lemma:z*-exists}. 
This is a one-dimensional problem which could be solved efficiently via bisection and convex optimization. That is,
 we use convex optimization to compute $f(x)$ for any given $z\in \left[\underline{Z}^0,\overline{Z}^0 \right]$ and use bisection to find $z^*$ such that $f(z^*) = 0$.
In comparison with its deterministic counterpart, the robust problem  requires an extra computing cost of $\delta \cO(\ln(1/\epsilon))$  to obtain a constant markup that is in the $\epsilon$-neighbourhood of the optimal solution, where $\delta$ is the computation cost to solve the adversary problem.

\section{{Robust Pricing under Partially Heterogeneous Price Sensitivity Parameters}}
\label{sec:uncon-hete}
We relax the assumption that the PSP are homogeneous. Completely relaxing this assumption makes the pricing  problem challenging, even for its deterministic version \citep{Gallego2014multiproduct}.   
Thus, we  assume that the products can be separated into partitions, and the PSP can be different over partitions. 
More specifically, we assume  that the products can be  partitioned into disjoint subsets and the products in each partition share the same PSP, and the CPGF is also partition-wise separable. This assumption has been  used in previous work to derive tractable solutions to the deterministic pricing problems \citep{zhang2018multiproduct}.
More precisely, we partition the set of all products $\cV$ into $N$  non-empty  subsets $\cV_1,\ldots,\cV_N$ such that $\cV = \bigcup_{n=1}^N \cV_n$ and $\cV_i \cap \cV_j = \emptyset, \ \forall i\neq j, i,j\in[N]$. Moreover, we separate the vector $\bY$ into sub-vectors $\bY^1,\ldots,\bY^N$ such that $\bY^n = \{Y_i|\ i\in \cV_n\}$ for all $n \in [N]$.
We  assume that the GEV-CPGF $G(\bY)$ can be separated into $N$ GEV-CPGFs as 
\[
G(\bY) =  \sum_{n=1}^N G^n(\bY^n).
\]
Note  that the nested logit model \citep{Ben1973structure}, one of the most widely-used GEV models in the literature, also has this separating structure. For notational convenience, we also separate  $(\ba,\bb) \in \bbR_+^{2m}$ into sub-vectors $(\ba^1,\bb^1),\ldots, (\ba^N,\bb^N)$ such that $(\ba^n,\bb^n) = \{(a_i,b_i)|\ i\in \cV_n\}$, $\forall n\in [N]$. 
 
 To deal with the robust problem, we further assume that the uncertainty set $\cA$ is also partition-wise separable, i.e., $\cA = \otimes_{n\in[N]} \cA^n$, where $\otimes$ is the Cartesian operation,  $\cA^n\subset \bbR^{2|\cV_n|}$ is the  uncertainty set for the sub-vector $(\ba^n,\bb^n)$, and $\cA^n$ are \textit{convex and bounded} for all $n\in[N]$. In other words, we assume that the vector of choice parameters can vary independently across partitions. 
 This assumption is a bit restrictive, but important to maintain the tractability of the robust problem. The reason is that if we use a general uncertainty set that allows for dependency between the choice parameters from different partitions, the adversary problem itself is generally not convex or quasi-convex in $(\ba,\bb)$, even under constant-markup prices, thus not tractable to solve. 
 On the other hand, the assumption will allow us to handle each function $G^n(\bY^n)$ independently, thus making it possible to convert the robust optimization problem into a convex one. In fact, one can construct a partition-wise separable uncertainty set 
 by collecting some samples of choice parameter estimates from each partition. We will discuss this in more detail in Section \ref{sec:expertiments}.

%Under the above assumptions, our main result in this section is to show that the robust problem can be converted equivalently into a convex optimization.

%In this context, we need to further assume that the uncertainty set is rectangular, because of the following reason. 
In this context, the difficulty lies in the fact that the optimal prices to the deterministic pricing problem do not have a single constant markup over all products.
%\citep{zhang2018multiproduct}.
As a consequence, the  robust optimal prices to  \eqref{prob:RO} would  generally not have a single constant markup over all the products and  the corresponding adversary's objective function  
would not be quasi-convex  and solutions to the adversary's problem may not be unique.
For this reason, we can not apply the techniques used in the previous section to identify a saddle point of the  robust problem. %In fact,
%one can show that a saddle point of the robust problem with partially homogeneous PSP can be found by solving a nonlinear system (Appendix \ref{appdx:no-constant-markup}), but this system is generally not tractable to handle.
 
In the rest of the section, we will show that the robust problem can be converted equivalently into a convex optimization.

To start our exposition, we note that, in analogy to the analysis in the case of homogeneous PSP, we also see that the \textit{min-max} counterpart always yields a partition-wise constant-markup solution \citep{zhang2018multiproduct}. Thus, if there is a saddle point in the robust problem, then it should have a partition-wise constant-markup form. This motivates us to find such a saddle point of the \textit{max-min} problem.

First, let us look at the robust problem where we only seek prices that have a constant markup in each partition, i.e., $\bx\in \bX^N$, where $\bX^N = \{\bx\in \bbR^m|\;x_i-c_i = x_j-c_j,\forall i,j\in\cV_n, n\in[N]\}$.
Let $z_n = x_i-c_i$ for all $i\in\cV_n$ and $n\in[N]$. The robust problem becomes
\begin{equation*}
 \max_{\bz \in \bbR^N} \left\{\min_{(\ba,\bb)\in \cA}  \sum_{n\in[N]} \sum_{i\in \cV_n} z_n P_i(\bY^n(z_n+\bc,\ba,\bb)|G^n) \right\},
\end{equation*}
or equivalently 
\begin{equation}\label{prob:RO-GEV-hete}
 \max_{\bz \in \bbR^N} \left\{\min_{(\ba,\bb)\in \cA}  \frac{\sum_{n\in[N]} z_n G^n(\bY^n| z_n,\ba^n,\bb^n)}{1 + \sum_{n\in[N]} G^n(\bY^n|z_n,\ba^n,\bb^n)} \right\},
\end{equation}
where $G^n(\bY^n|z,\ba^n,\bb^n) = G^n(Y_i,\ i\in\cV_n)$ with $Y_i = e^{a_i-b_i(z_n+c_i)}$, for all $i\in \cV_n$. For notational brevity, let $$\rho(\bz,\ba,\bb) = \frac{\sum_{n\in[N]} z_n G^n(\bY^n|z_n,\ba^n,\bb^n)}{1 + \sum_{n\in[N]}  G^n(\bY^n|z_n,\ba^n,\bb^n)}.$$
%%%%%%%%%%%%%%%%%%%%%%%%%%%%%%%%%%%%%%%%%%
Let us  also denote 
\[
%\overline{\cG}^n(z_n) = \max_{(\ba^n,\bb^n) \in\cA^n}\; \Big\{ G^n(\bY^n| z_n,\ba^n,\bb^n)\Big\};\; 
\underline{\cG}^n(z_n) = \min_{(\ba^n,\bb^n) \in\cA^n}\; \Big\{ G^n(\bY^n| z_n,\ba^n,\bb^n)\Big\}.
\]
Since $G^n(\bY^n| z_n,\ba^n,\bb^n)$ is strictly convex in $(\ba^n,\bb^n)$ (Proposition \ref{prop:G-strictly-convex}), we see that  $\underline{\cG}^n(z_n)$ is continuous and differentiable in $z_n$. Let $(\ba^{n*}(z_n), \bb^{n*}(z_n) =  \text{argmin}_{(\ba^n,\bb^n) \in\cA^n}\; \Big\{ G^n(\bY^n| z_n,\ba^n,\bb^n)\Big\}$, which are always uniquely determined given any $z_n\in\bbR$ (Proposition \ref{prop:G-strictly-convex}). 
To handle the robust problem \eqref{prob:RO-GEV-hete}, let us consider the following\textit{ reduced optimization problem}, which is obtained by forcing each component $G^n(\cdot)$ to its minimum value over $\cA^n$.  
\begin{equation}
\label{prob:RO-hete-simlified}
 \max_{\bz\in\bbR^N}  \left\{ \cW(\bz) = \frac{\sum_{n\in[N]} z_n \underline{\cG}^n(z_n)}{1 + \sum_{n\in[N]} \underline{\cG}^n(z_n)} \right\}.
\end{equation}
In the rest of the section, we will focus on solving the robust problem \eqref{prob:RO} by making use of Problems \eqref{prob:RO-GEV-hete} and \eqref{prob:RO-hete-simlified}. More specifically, we will prove the following chain of results.
\begin{itemize}
\item[\textbf{(i)}] The \textit{reduced problem} \eqref{prob:RO-hete-simlified} always yields a unique solution, and this solution can be found by convex optimization (Theorem \ref{theor;RO-hete-convexness-simplified-model}).
\item[\textbf{(ii)}] Any optimal solution to \eqref{prob:RO-hete-simlified} is also a robust solution to \eqref{prob:RO-GEV-hete}  and vice-versa (Theorem \ref{theor:convert-simplified-model}).
\item[\textbf{(iii)}] A solution to \eqref{prob:RO-GEV-hete} forms an optimal solution to robust problem \eqref{prob:RO} (Theorem \ref{theor:main-hete}).
\end{itemize}
To make the technical results easier to follow, we separate the rest of the section into two subsections, where Section \ref{sec:concave-reduced-model} will focus on the \textit{reduced problem}, and Section \ref{sec:hetePSP-solving-robust-problem} shows 
how to  convert the original robust problem  \eqref{prob:RO-GEV-hete} into the reduced one, which eventually leads to the result that \eqref{prob:RO-GEV-hete} and  \eqref{prob:RO} can be solved by convex optimization. 

%%%%%%%%%%%%%%%%%%%%%%%%%%%%%%%%
%%%%%%%%%%%%%%%%%%%%%%%%%%%%%%%%
%%%%%%%%%%%%%%%%%%%%%%%%%%%%%%%%

\subsection{Convexity of the Reduced Problem}
\label{sec:concave-reduced-model}
The \textit{reduced problem} is indeed not convex if it is defined in terms of the prices $\bx$. Nevertheless, we can show that it becomes convex if we view it under purchase probabilities. More precisely, we will do some change of variables. Let use denote a vector $\bp^{G} \in \bbR^{N}$ with entries
\begin{equation}
\label{eq:pG-z}
p^{G}_n =  \frac{\underline{\cG}^n(z_n)}{1 + \sum_{n\in[N]} \underline{\cG}^n(z_n)},\;\forall
n\in[N],
\end{equation}
then the objective function in \eqref{prob:RO-hete-simlified} can be written as  $\cW(\bz) = \sum_{n\in[N]} z_np^G_n$. 
This vector $\bp^G$ can be interpreted as an aggregated purchase probabilities for the partitions, i.e., $p^G_n = \sum_{i\in\cV_n} p_i$, where $p_i$ is the purchase probability of item $i\in\cV$.  
In Theorem \ref{theor;RO-hete-convexness-simplified-model} below, we show that, given any $\bp^G \in \cP^G= \{\bp^G\in \bbR^N_+|\ \sum_{n\in[N]} p^G_n <1 |\}$, there is a unique $\bz(\bp^G) \in \bbR^{N}$ satisfying \eqref{eq:pG-z}. Moreover, Problem \eqref{prob:RO-hete-simlified} can be formulated as a convex optimization program of variables $\bp^G$. Note that a similar result has been shown  previously \citep{zhang2018multiproduct} for the case that the choice parameters $(\ba,\bb)$ are fixed. In our setting, $(\ba,\bb)$ are  a solution to convex optimization problems parameterized by $\bz(\bp^G)$, thus requiring a new and more complicated  proof.
\begin{theorem}[Convexity of the reduced problem]
\label{theor;RO-hete-convexness-simplified-model}
Given any $\bp^G\in \cP^G$, there is a unique vector $\bz(\bp^G) \in\bbR^N$ satisfying \eqref{eq:pG-z}, and this vector can be found by bisection. Moreover, $\cW(\bz(\bp^G))$ is strictly concave in $\bp^G$. 
\end{theorem} 
To prove the result, we first show that each function $\underline{\cG}^n(z_n)$ is invertible. That is, for any $\alpha>0$ there is a unique $z_n\in \bbR_+$ such that $\underline{\cG}^n(z_n) = \alpha$. This allows us to define the inverse function $(\underline{\cG}^n)^{-1}$ such that $(\underline{\cG}^n)^{-1}(\underline{\cG}^n(z_n)) = z_n$. This inverse function can be computed by bisection.
The existence of the inverse function is necessary for the claim that there is always a unique vector $\bz$ that yields a given purchase probability vector $\bp^G \in\cP^G$ and this vector can be computed as
\[
\bz(\bp^G)_n = (\underline{\cG}^n)^{-1}\left(\frac{p^G_n}{1-\sum_{l\in[N]}  p^G_l}\right).
\]
To show the convexity of $\cW(\bz(\bp^G))$ in $\bp^G$, we first validate convexity of its deterministic counterpart, i.e., the version in which all the choice parameters are given 
$
\widetilde{\cW}(\widetilde{\bz}(\bp^G|\ba,\bb)) = \sum_{n\in[N]} \widetilde{\bz}(\bp^G|\ba,\bb)_n p^G_n,
$
where $\widetilde{\bz}(\bp^G|\ba,\bb) \in \bbR^N$ are a vector of constant-markups that archive vector $\bp^G$  as 
\begin{equation}
\label{eq:pG-z-fixed-ab}
p^{G}_n =  \frac{G^n(\bY^n| \widetilde{z}_n,\ba^n,\bb^n)}{1 + \sum_{l\in[N]} G^l(\bY^l| \widetilde{z}_l,\ba^l,\bb^l)},\;\forall n\in [N].
\end{equation}
Once the convexity of $\widetilde{\bz}(\bp^G|\ba,\bb)$ is validated,  we can further take the derivatives of $\bz(\bp^G)$ with respect to $\bp^G$ and show that they are equal to the derivative values of a deterministic function. This is the key result to show that the second-order derivative of $\cW(\bz(\bp^G))$ is positive-definite,  leading to the convexity of $\cW(\bz(\bp^G))$. We provide the detailed proof in Appendix \ref{proof:RO-hete-convexness-simplified-model}.   

We further characterize a solution to \eqref{prob:RO-hete-simlified}. In Proposition \ref{prop:RO-hete-fixed-point} below we show that  \eqref{prob:RO-hete-simlified} always has a unique local optimal $\bz^*$ (i.e., $\cW(\bz)$ is unimodal), and this solution will satisfy a fixed point system that is an extended version of the one shown in Theorem \ref{theor:RO-GEV-homo}. Note  that the uniqueness of a \textit{local optimal solution} of \eqref{prob:RO-hete-simlified} defined in terms of $\bp^G$ is straightforward due to the concavity of $\cW(\bz(\bp^G))$. It is however not trivial when the objective function is defined in terms of $\bz$.  
\begin{proposition}
\label{prop:RO-hete-fixed-point}
Problem  \ref{prob:RO-hete-simlified} always yields a unique local optimal solution $\bz^*$ (i.e., $\cW(\bz)$ is unimodal) and this solution satisfies the following fixed point system
\begin{equation}
\label{eq:temp3452}
z_n =\frac{1}{\me{{\bb}^{n*}(z_n)}}  + \sum_{l\in[N]}\frac{\underline{\cG}^l(z_l)}{\me{{\bb}^{l*}(z_l)}},\; \forall n\in[N].
\end{equation}
\end{proposition}
%\proof{Proof of Proposition \ref{prop:RO-hete-fixed-point}}
Proposition \ref{prop:RO-hete-fixed-point} implies that solving the fixed-point problem \eqref{eq:temp3452} will yield a solution to \eqref{prop:RO-hete-fixed-point}. However, directly solving  
\eqref{eq:temp3452} would be not tractable. Instead,  Theorem \ref{theor;RO-hete-convexness-simplified-model} show that it can be solved conveniently by convex optimization. The fixed-point system in Proposition \ref{prop:RO-hete-fixed-point} is however important to establish the saddle point result in the next section (Proposition \ref{prop:duality}). 
%%%%%%%%%%%%%%%%%%%%%%%%%%%%%%%%%%%%%%
%%%%%%%%%%%%%%%%%%%%%%%%%%%%%%%%%%%%%%

\subsection{Solving the Robust Problem}
\label{sec:hetePSP-solving-robust-problem}
We know from the previous section that the reduced problem is tractable to solve. 
We now move to the second part showing that the original robust optimization problem can be converted into the \textit{reduced problem}, for which a solution can be found  by convex optimization. We first state the following result connecting the \textit{reduced problem} and \eqref{prob:RO-GEV-hete}.
\begin{theorem}[Equivalence between \eqref{prob:RO-GEV-hete} and the reduced problem]
\label{theor:convert-simplified-model}
Any optimal solution to \eqref{prob:RO-GEV-hete} is also optimal to \eqref{prob:RO-hete-simlified} and vice-versa.  
\end{theorem}
The general idea to prove the theorem is to show that, under the optimal price solution, the adversary will force each component $G^n(\bY^n|z_n,\ba^n,\bb^n)$ of the objective function to its minimum values. We refer the reader to Appendix \ref{apdx:proof-concert-to-reduced-model} for a detailed proof. 

We now come back to the original robust problem \eqref{prob:RO} with partition-wise homogeneous PSP. We will gather all the results established above to show how we can get an optimal solution of \eqref{prob:RO} by convex optimization. 
%We will also provide some comparative analyses when the size of the uncertainty set varies.
Before stating the main theorem, let us introduce the following result saying that a solution obtained by solving the reduced problem \eqref{prob:RO-hete-simlified} forms a saddle point to \eqref{prob:RO-GEV-hete}, thus the minimax equality holds. 
\begin{proposition}[Saddle point of \eqref{prob:RO-GEV-hete}]
\label{prop:duality}
If $\bz^*$ is a solution to \eqref{prob:RO-GEV-hete}, then $(\bz^*,\ba^*(\bz^*),\bb^*(\bz^*))$ is a saddle point of the \textit{max-min} problem \eqref{prob:RO-GEV-hete}. As a result, the minimax equality holds, i.e., 
\[
\max_{\bz\in\bbR^N}\min_{(\ba,\bb) \in\cA} \rho(\bz,\ba,\bb) = \min_{(\ba,\bb) \in\cA} \max_{\bz\in\bbR^N} \rho(\bz,\ba,\bb).
\]
\end{proposition}
\proof{Proof:}
It is clear from Theorem \ref{theor:convert-simplified-model} that if $\bz^*$ to a solution to \eqref{prob:RO-GEV-hete}, then  $(\ba^*(\bz^*),\bb^*(\bz^*))$ is a solution to the corresponding adversary's problem. Moreover, from Proposition \ref{prop:RO-hete-fixed-point} and Theorem C1 of \cite{zhang2018multiproduct}, we also see that $\bz^*$ is a solution to the pricing problem under fixed parameters $(\ba^*(\bz^*),\bb^*(\bz^*))$, i.e., 
$\max_{\bz \in\bbR^N} \rho(\bz,\ba^*(\bz^*),\bb^*(\bz^*))
$. 
Thus, $(\bz^*,\ba^*(\bz^*),\bb^*(\bz^*))$ is clearly a saddle point of \eqref{prob:RO-GEV-hete}. The minimax equality follows directly from the existence of a saddle point. 
\endproof

We now gather all the previous results to establish our main theorem.
Theorem \ref{theor:main-hete} below states that a solution to the robust problem \eqref{prob:RO} will have a constant-markup style and this constant-markup vector can be found be convex optimization. The proof can be done easily given all the claims we have in Section \ref{sec:concave-reduced-model} and Theorem \ref{theor:convert-simplified-model} above.
\begin{theorem}
\label{theor:main-hete}
{\bf(A partition-wise constant-markup solution is optimal to the robust problem).}
\label{theor:RO-GEV-hete}
Under partition-wise homogeneous PSP and partition-wise decomposable uncertainty sets, the robust problem \eqref{prob:RO} yields a unique partition-wise constant-markup solution $\bx^*$ such that $x^*_i = z^*_n + c_i$, $\forall n\in[N], i\in\cV_i$, where $\bz^*$ is a unique solution to Problem \eqref{prob:RO-GEV-hete}, which can be solved by convex optimization. Moreover, $(\bx^*,\ba^*(\bz^*),\bb^*(\bz^*))$ is a saddle point of \eqref{prob:RO} and the minimax equality holds, i.e.,
\[
\max_{\bx\in\bbR^m}\min_{(\ba,\bb) \in\cA} \Phi(\bx,\ba,\bb) = \min_{(\ba,\bb) \in\cA} \max_{\bx\in\bbR^m} \Phi(\bx,\ba,\bb).
\]
\end{theorem}
\proof{Proof:}
We first prove the minimax equality property by the chain
\begin{align}
  \max_{\bx\in\bbR^m}\min_{(\ba,\bb) \in\cA} \Phi(\bx,\ba,\bb) &\stackrel{(a)}{\leq}\min_{(\ba,\bb) \in\cA}\max_{\bx\in\bbR^m} \Phi(\bx,\ba,\bb) \nonumber \\
  &\stackrel{(b)}{ = }\min_{(\ba,\bb) \in\cA}\max_{\bz\in\bbR^N} \rho(\bz,\ba,\bb) \nonumber\nonumber \\
  &\stackrel{(c)}{= }\max_{\bz\in\bbR^N} \min_{(\ba,\bb) \in\cA} \rho(\bz,\ba,\bb) \nonumber\nonumber \\
&{\leq}\max_{\bx\in\bbR^m} \min_{(\ba,\bb) \in\cA} \Phi(\bx,\ba,\bb), \nonumber
\end{align}
where $(a)$ is from the well-known  max-min inequality, $(b)$ is from the property that any deterministic pricing problem (with fixed $(\ba, \bb)$) always yields a partition-wise constant-markup solution, $(c)$ is due to the minimax equality of \eqref{prob:RO-GEV-hete} shown in Proposition \ref{prop:duality} above.
This chain of (in)equalities leads to the minimax equality property of \eqref{prob:RO} and the result that $\bx^*$ defined by a constant-markup solution $\bz^*$ of \eqref{prob:RO-GEV-hete} is a robust solution to \eqref{prob:RO} and $(\bx^*,\ba^*(\bz^*),\bb^*(\bz^*))$ forms a saddle point of the \textit{max-min} problem
\eqref{prob:RO}. 
\endproof

%\ctien{Say something}
%Theorem \ref{theor:RO-GEV-hete} indicates that one needs to solve a convex optimization problem to find a robust solution to \eqref{prob:RO}. In the case of deterministic pricing, an optimal solution can be computed much easier by bisection \citep{zhang2018multiproduct}. 

We now discuss in detail how to solve the reduced problem \eqref{prob:RO-hete-simlified}. Since the problem is convex when  the objective function is defined in terms of the purchase probability $\bp^G$, we show how to compute  $\cW(\bz(\bp^G))$ and its gradients, which are crucial for the optimization process. 
Given a purchase probability $\bp^G\in\cP^G$, from Lemma \ref{lemma:zp-unique-deter}, we can compute $\cW(\bz(\bp^G))$ as
\[
\cW(\bz(\bp^G)) = \sum_{n \in[N]} z(\bp^G)_n p^G_n  = \sum_{n\in[N]} (\underline{\cG}^n)^{-1}\left(\frac{p^G_n}{1-\sum_{l\in[N]}  p^G_l}\right) p^G_n
\]
where the inverse function $\underline{\cG}^n)^{-1}(\cdot)$ can be computed efficiently by \textit{bisection}. The gradients of 
$\cW(\bz(\bp^G))$ are more difficult to get and we show how to do it in Proposition \ref{prop:gradient-cW} below.
\begin{proposition}[Gradients of $\cW(\bz(\bp^G))$]
\label{prop:gradient-cW}
For any $\bp^G\in\cP^G$, we have
\begin{equation}
\label{eq:gradient-cW}
\frac{\partial \cW(\bz(\bp^G))}{\partial p^G_n}  =  z(\bp^G)_n - \frac{1}{\me{\bb^{k*}(z_k)}} -\frac{1}{ (1-\bbe^\T \bp^G)} \sum_{k\in [N]} \frac{p^G_k}{\me{\bb^{k*}(z_k)}},\; \forall n\in [N],
\end{equation}
where $ (\ba^{k*}(z_k),\bb^{k*}(z_k)) =  \text{argmin}_{(\ba^k,\bb^k) \in\cA^k}\; \Big\{ G^k(\bY^k| z_k,\ba^k,\bb^k)\Big\}$. 
\end{proposition}
The proof (details in Appendix \ref{apdx:proof-gradient-cW}) can be done by directly taking the derivatives of $\cW(\bz(\bp^G))$ with respect to $\bp^G$ and using \eqref{eq:grad-cG}. The computation of $\cW(\bz(\bp^G))$ for a given $\bp^G\in\cP^G$ can be done by performing the following steps: (i) compute $\bz(\bp^G)$ using Lemma \ref{lemma:zp-unique-deter}, (ii) compute $(\ba^*(\bz),\ba^*(\bz))$  as (unique) optimal solutions of the problems $ \min_{(\ba^n,\bb^n) \in\cA^n}\; \Big\{ G^n(\bY^n| z_n,\ba^n,\bb^n)\Big\}$, $n\in [N]$, (iii) compute $\cW(\bz(\bp^G)) = \bz(\bp^G)^\T \bp^G$  and its gradients by \eqref{prop:gradient-cW}. Since the objective function is strictly concave, we know that the optimization problem can be solved efficiently by a convex optimization solver.
When the uncertainty set is rectangular, the reduced optimization problem can be further simplified, as a solution to $ \min_{(\ba^n,\bb^n) \in\cA^n}\; \Big\{ G^n(\bY^n| z_n,\ba^n,\bb^n)\Big\}$  can be identified, thus the reduced problem can be transformed equivalently to a deterministic pricing problem with fixed choice parameters  and it is known that such a deterministic pricing problem yields closed form solutions \cite{zhang2018multiproduct}.
We state this result in the following corollary.
\begin{corollary}[Rectangular uncertainty sets]
If the uncertainty is rectangular, i.e., $$\cA = \left\{(\ba,\bb)|\ \ba\in [\underline{\ba},\overline{\ba}],\ \bb \in [\underline{\bb},\overline{\bb}], \ b_i = b_j,\ \forall i,j\in \cV_n,\ \ \forall n \right\},$$
then the robust problem \eqref{prob:RO} is equivalent to the deterministic pricing problem  
$
 \max_{\bx_\in\bbR^m} \Phi(\bx,\underline{\ba},\overline{\bb}). 
$
\end{corollary}
The result is easy to validate, as from Lemma \eqref{lemma:G-bounded}  we see that $(\underline{\ba}^n,\overline{\bb}^n) = \text{argmin}_{(\ba,\bb)\in\cA} G^n(\bY^n|z_n,
\ba^n,\bb^n)$ for any $z_n\in\bbR$.

\section{Numerical experiments}
\label{sec:expertiments}
We provide experimental results to show how the robust model considered above (i.e., robust unconstrained pricing with homogeneous and partition-wise homogeneous PSP)  protect us from choice parameter uncertainties. 
We first discuss our approach to construct uncertainty sets and different baseline approaches for the sake of comparison.
\subsection{Constructing Uncertainty Sets}
\label{sec:construct-uncertainty-set}
%We discuss our approach  to construct uncertainty sets to deal with the issue of choice parameter uncertainty, noting that a similar approach has been used in 
Inspired by  \cite{Rusmevichientong2012robust} in the context of robust assortment optimization, 
such an uncertainty set can be created for the situation that the market is heterogeneous, i.e., the market has several costumer types and the choice parameters would vary across them, but the proportion of each customer type is not known with certainty. 
To be more precise, let assume that the true parameters  for the underlying GEV choice model can be one of $K$ vectors $\{(\ba^{(1)}, \bb^{(1)}),\ldots, (\ba^{(K)}, \bb^{(K)})\}$, representing $K$ types of customers. For ease of notation, let $\bw^k = (\ba^{(k)}, \bb^{(k)})$ for all $k\in [K]$.  Let $\tau_1,\ldots,\tau_K \in [0,1]$ be the proportion of each customer type with $\sum_{k\in [K]}\tau_k = 1$. 
\iffalse
Once the choice parameters and the proportions $\tau_k, k\in [K]$, for each customer type are known exactly, a pricing optimization approach would be to maximize the average expected revenue $\max_{\bx \in \bbR^m_+}  \left\{\sum_{k\in[K]} \tau_k \Phi(\bx,\bw^k)\right\}$, which is generally not tractable for $m\geq 2$ \citep{Li2019productMMNL}. 
An alternative and tractable approach would be use the weighted average of the choice parameters $\widetilde{w} = \sum_{k\in [K]} \tau^k \bw^k$ and solve the deterministic problem \eqref{prob:pricing-deterministic}. 
\fi
We are interested in the situation that the proportions can be estimated somehow using historical data, but estimation may have errors and the  proportion estimates may not make good representation to the ``true'' ones.
%far away from the ``true'' proportions.
%if the proportions of the  $K$ customer types are not known with certainty but not too far from the prior proportions $\{\tau_k,k\in[K]\}$,
In this situation, an uncertainty set can be constructed around the proportion estimates as
\begin{equation}
\label{eq:uncertainty-set-1}
\cA = \left\{\bw = \sum_{k\in [K]} \lambda_k \bw^k\Big|\; \sum_{k\in [K]} \lambda_k = 1 \text{ and } \max_{k\in [K]}|\lambda_k-\tau_k  | \leq \epsilon \right\}.
\end{equation}
In the partially homogeneous case, as our results require partition-wise separable uncertainty sets, such an uncertainty set can be constructed in a similar way as follows. For each customer type $k$,  let $\bw^{k,n}$ be the vector of choice parameters of partition $n\in[N]$. The uncertainty set for each partition can be defined as
\begin{equation}
\label{eq:uncertainty-set-2}
\cA_n = \left\{\bw^{n} = \sum_{k\in [K]} \lambda^n_k \bw^{k,n}\Big|\; \sum_{k\in [K]} \lambda^n_k = 1 \text{ and } \max_{k\in [K]}|\lambda^n_k-\tau_k  | \leq \epsilon \right\},\; \forall n\in [N].
\end{equation}
Here $\epsilon\in[0,1]$ reflects an ``uncertainty level'' of the uncertainty set. Larger $\epsilon$ values provide larger uncertainty sets, corresponding to more conservative models that may help protect well against worst-case scenarios, but may lead to low average performance. On the other hand, smaller $\epsilon$ values provide smaller uncertainty sets and would lead to less conservative robust solutions, which may perform well in terms of average performance but would be worse in protecting bad scenarios of the choice parameters.
Adjusting $\epsilon$ would help the firm balance the worst-case protection and average performance.
Clearly, $\epsilon = 0$ corresponds to the deterministic case, i.e., the proportion of each customer type are given with certainty, and $\epsilon = 1$ reflects the situation that we are totally uncertain about how likely the proportion of each customer type is, and have to ignore the predefined  proportions $\{\tau_1,\ldots,\tau_K\}$. 

\subsection{Baseline Approaches}
We discuss tractable baseline approaches that would be used to solve the pricing problem when facing the issue of choice parameter uncertainty.  
A straightforward approach would be to  employ the mean values of the choice parameters  and solve the deterministic version. In this context, we know that the pricing problem is computationally tractable.
Alternatively, one may look at different possibilities of the choice parameters and define a mixed formulation where the market is divided into a finite number of market
segments and each segment is governed by a scenario of the choice parameters. However, one can show that the expected revenue in this context is no longer unimodal and the \textit{constant-markup} property identified for the GEV pricing problem no-longer holds, {even if there are only two market segments} \citep{Li2019productMMNL}.
As a result, this mixed version is not computationally tractable.

Another baseline approach  is to sample some choice parameters from the uncertainty set and use simulation to select a solution that provides best protection from worst-case scenarios. More precisely, let assume that the firm needs to make a pricing decision while being aware that the choice parameters may vary in an uncertainty set. In this context, the firm can sample some points from the uncertainty set 
and compute the corresponding optimal prices for each selection, using the deterministic approach. Then, for each price vector, the firm can sample a \textit{sufficiently large} number of vector of choice parameters from the uncertainty set, in order to evaluate how each price vector obtained performs when the choice parameters vary in the uncertainty set. This can be done by simply selecting the solution that gives the best worst-case profit among the samples. This approach may be computationally tractable with a reasonable number of samples, but would be much more computationally expensive than the robust and deterministic approaches. 
We refer to this as the sampling-based  approach.
One can show that solutions given by this approach will converge to those from the robust counterpart when the sample sizes grow to infinity.

In these experiments, we will compare our robust models (denoted as RO), which are computationally tractable, against the sampling-based approach (denoted as SA) and the deterministic one  with mean-value choice parameters (denoted as DET). We will employ two popular GEV models in the literature, i.e., the MNL and nested logit models. %, when the PSP are the same over all the products.
  {For the SA approach, we sample points uniformly from the uncertainty set since we do not make any assumption about the distribution of the choice parameters.
One can argue that the uniform distribution may not be the best choice in the case that the firm believes that it has some ideas (perhaps via estimation) about the distribution of the choice parameters. Nevertheless, estimating such a distribution is not  easy  in practice. A common approach in choice modeling is to assume that  the parameters follow some distributions (e.g. normal distribution) with unknown coefficients and try to estimate these  coefficients by maximum likelihood estimation \citep{Mcfadden2000mixed}. This approach, even though popular, does not guarantee that the distribution obtained is the \textit{true} distribution of the choice parameters, assuming that there exists a true distribution.
As such, the distribution of the choice parameters is typically only known ambiguously. Distributionally robust optimization is a \textit{robust} approach that is explicitly designed to handle this ambiguity \citep{Shapiro2018tutorialDRO}, which we keep for future research.}

%The goal here is to compare the robust (RO) approach against the standard deterministic (DET)  and the sampling-based (SA)  approaches under two popular GEV models in the literature, i.e., the MNL and nested logit models, when the PSP are the same over all the products. In this context, we know that the optimal prices for all the cases (RO, DET, and SA) have a constant-markup with respect to the product costs, and this constant markup can be computed by a closed-form formula for the DET and SA approaches and by binary search for the RO approach.  

\subsection{Experimental Settings}

We choose $m = 50$, 
%noting that in this experiment,  there is only one price sensitive parameter for all the products. 
 and  $K = 5$ (i.e., there are 5 customer types) and randomly choose the proportions  and the underlying choice parameter vectors $\{\bw^1,\ldots,\bw^5\}$. For each $\epsilon>0$ we define the uncertainty set as in \eqref{eq:uncertainty-set-1}.
%Given an uncertainty level $\epsilon>0$, we define the polyhedron uncertainty  $\cA = \{(\ba,b)|\ ||\ba - \ba_0||_1 + \beta||b- b_0||_1 \leq \epsilon\}$, where $||.||_1$ stands for the L1 norm and $\beta$ is the scale of the PSP $b$ with respect to $\ba$. 
The comparison is done as follows. For each $\epsilon$, we solve the corresponding robust problem and obtain a robust solution $\bx^{\RO}$. For the DET, we solve the deterministic model with the weighted average  parameters $\widetilde{\bw} = \sum_{k\in [K]} \tau_k \bw^k$ and obtain an optimal solution $\bx^\DET$.
For the SA approach, we sample randomly and uniformly $s_1$ points from the uncertainty set, and for each point compute the corresponding optimal prices, which have a constant markup over products. For each pricing solution, we again sample randomly and uniformly $1000$ choice parameters from $\cA$, and compute and pick a pricing solution with the largest worst-case expected revenue among the 1000 samples. We test this approach with $s_1 = 10$ and $s_1 = 50$ and denote the corresponding solutions as $\bx^{\textsc{SA10}}$, $\bx^{\textsc{SA50}}$, respectively. Larger $s_1$ can be chosen, but it would mean that the SA becomes way more expensive as compared to the RO and DET approaches. For example, if we choose $s_1 = 100$, the SA requires to solve 100 deterministic problems and compute $10^5$ expected revenues to obtain a pricing solution.

\subsection{Comparison Results}
We provide experimental results for the robust model under the nested logit model. The CPGF of the nested logit model is given as
$
G(\bY) = \sum_{n\in [N]} \left(\sum_{i\in C_n} Y_i^{\mu_n} \right)^{\mu/\mu_n},
$
where $[N]$ is the set of nests and for each $n\in [N]$, $C_n$ is the corresponding subset of the items, $\mu$ and $\mu_n$, $n\in\cN$ are the positive parameters of the nested logit model.
In this experiment, we separate the whole item set into 5 nests of the same size (10 items per each nest), i.e. $N = 5$ and $|
C_n| =10$ for all $n\in [N]$. To evaluate the performance of the three approaches when the choice parameters vary, given the uncertainty set defined above, we randomly and uniformly  sample 1000 parameters $(\ba, \bb)$ from the set $\cA$, and compute the expected revenues given by $\bx^\RO$, $\bx^\textsc{SA10}$, $\bx^\textsc{SA50}$, and $\bx^\DET$. So, for each solution, we get a distribution of expected revenues over 1000 samples. We then draw the histograms of  of the distributions to compare. We first provide experiments for the case of homogeneous PSP and then move to the case of partition-wise homogeneous PSP.

\subsubsection{Homogeneous PSP.}

\begin{figure}[htb]
 \centering
 \includegraphics[width=1.0\textwidth]{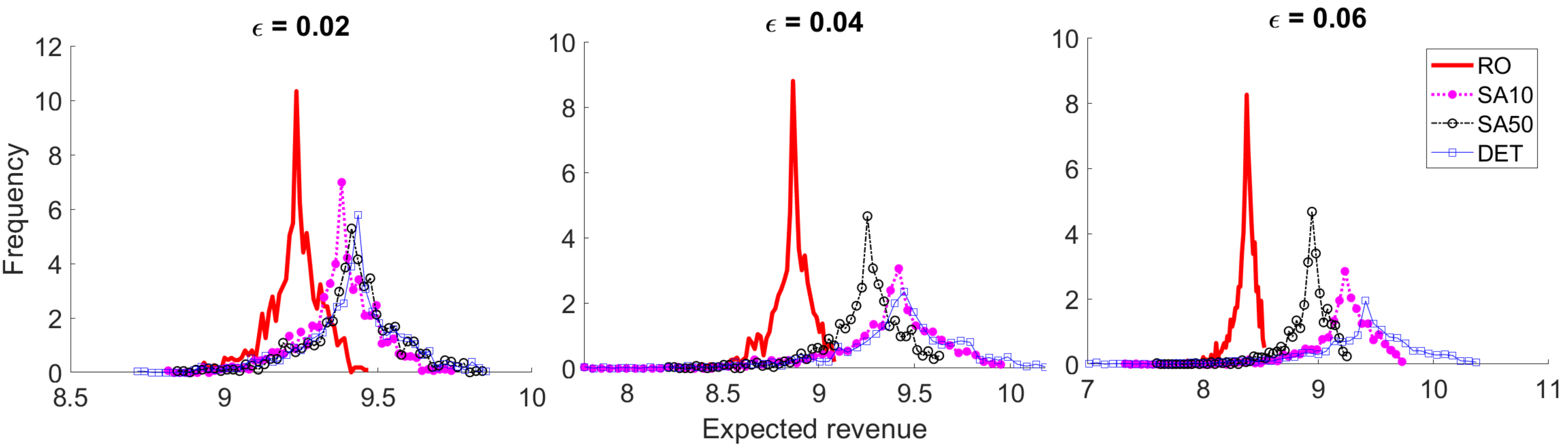}
 \caption{Comparison between revenue distributions given by optimal price vectors given by the robust (RO) and deterministic (DET) and sampling-based (SA10 and SA50) approaches, under the nested logit model and  different uncertainty levels $\epsilon$.}
 \label{fig:Nested-uncontrained}
\end{figure}

The histograms of the distributions obtained in 
Figure \ref{fig:Nested-uncontrained} for $\epsilon \in \{0.02,0.04, 0.06\}$.  
We see that the distributions  given by the RO approach always have higher peaks, lower variances and shorter tails, as compared to the other  approaches. The difference becomes clearer with larger $\epsilon$ 
This demonstrates the capability of the RO approach in giving  
not-too-low revenues.
In addition,  the sampling-based approach (SA10 and SA50) perform better then the DET in terms of protecting us against too low revenues. In this aspect, the  SA50 also performs better than the SA10, especially when $\epsilon$ increases.

In Table \ref{tb:RO-uncon-Nested}, we provide more details about the average and worst-case values of the distributions given by the three approach. In particular, we compute the ``\textit{percentile ranks}'' of the RO worst-case revenues, which indicates the percentages that the expected revenues given by the baseline approaches (DET, SA10 and SA50) are lower than the corresponding worst-case expected revenues given by the RO. For example, for $\epsilon = 0.1$, there are $23\%$ of the revenues given by the DET (over 1000 sampled revenues) are less than the corresponding RO worst-case revenue. Over $\epsilon \in\{0.02,\ldots,0.4\}$, the average percentile ranks of the RO worst-case revenues are 26.7\%, 9.5\% and 5.3\% for the DET, SA10 and SA50 approaches, respectively,  which clearly indicates  gains from the use of the RO approach. It can be seen that in terms of average revenue, the DET approach performs the best, followed by the the SA10, SA50 and RO approaches.  
In general, the baseline approaches (DET, SA10, SA50) always give higher average revenues, but lower worst-case revenues, which clearly indicates that the RO approach does a better job
in protecting us from worst-case situations, but also show the trade-off of being robust. Moreover, the results in Table  \ref{tb:RO-uncon-Nested} also tell us that if the firm cares more about the worst cases, a large $\epsilon$ can be chosen to have better protection against too low expected revenues. On the other hand, if average performance is of concern, then by choosing a small $\epsilon$, one can still get a protection from the robust solutions, but also get an average performance that is comparable to that of the solutions by the deterministic approach. 
This observation is also consistent with those from other robust work in the revenue management literature \citep{Li2019robustCVaR,Rusmevichientong2012robust}.

\begin{table}[htb]
\centering
\begin{tabular}{r|rrrr|rrrr|rrr}
\multicolumn{1}{c|}{$\epsilon$} & \multicolumn{4}{c|}{Average}   & \multicolumn{4}{c|}{Worst}  & \multicolumn{3}{c}{\begin{tabular}[c]{@{}l@{}}Percentile rank~\\of RO worst-case\end{tabular}}   \\ 
\cline{2-12}
\multicolumn{1}{l|}{} & \multicolumn{1}{l}{DET} & \multicolumn{1}{l}{SA10} & \multicolumn{1}{l}{SA50} & \multicolumn{1}{l|}{RO} & \multicolumn{1}{l}{DET} & \multicolumn{1}{l}{SA10} & \multicolumn{1}{l}{SA50} & \multicolumn{1}{l|}{RO} & \multicolumn{1}{l}{DET} & \multicolumn{1}{l}{SA10} & \multicolumn{1}{l}{SA50}  \\ 
\hline
0.02& 9.4& 9.4 & 9.3 & 9.2& 8.7& 8.7 & 8.7 & 8.7& 4  & 8   & 4\\
0.04& 9.3& 9.1 & 9.0 & 8.8& 7.5& 8.1 & 8.1 & 8.2& 13 & 8   & 4\\
0.06& 9.2& 9.0 & 8.7 & 8.3& 6.0& 7.1 & 7.4 & 7.5& 17 & 8   & 5\\
0.08& 9.1& 8.1 & 8.8 & 7.8& 4.2& 7.1 & 6.9 & 7.0& 18 & 7   & 7\\
0.10& 8.9& 8.7 & 7.8 & 7.3& 2.6& 4.7 & 5.6 & 6.4& 23 & 12  & 6\\
0.12& 8.8& 7.7 & 7.6 & 6.7& 1.4& 5.8 & 5.9 & 5.9& 26 & 7   & 6\\
0.14& 8.6& 8.0 & 7.2 & 6.2& 1.1& 4.6 & 5.0 & 5.3& 26 & 9   & 4\\
0.16& 8.3& 6.6 & 7.2 & 5.8& 0.9& 4.9 & 4.9 & 5.1& 34 & 7   & 4\\
0.18& 8.4& 7.3 & 6.6 & 5.3& 0.5& 3.7 & 4.3 & 4.7& 30 & 9   & 7\\
0.20& 8.3& 7.4 & 6.7 & 4.9& 0.2& 2.0 & 2.8 & 4.1& 31 & 10  & 6\\
0.22& 8.2& 7.1 & 6.4 & 4.5& 0.1& 1.4 & 2.0 & 3.9& 33 & 11  & 6\\
0.24& 8.0& 6.6 & 5.8 & 4.2& 0.1& 2.7 & 2.7 & 3.4& 32 & 8   & 7\\
0.26& 8.1& 7.9 & 6.0 & 3.9& 0.0& 0.1 & 0.8 & 3.1& 31 & 28  & 5\\
0.28& 8.0& 5.4 & 5.5 & 3.7& 0.0& 2.3 & 2.1 & 2.9& 32 & 8   & 6\\
0.30& 7.7& 5.1 & 5.1 & 3.5& 0.0& 2.1 & 2.2 & 2.3& 37 & 7   & 5\\
0.32& 7.9& 4.4 & 3.8 & 3.3& 0.0& 1.9 & 2.0 & 2.2& 31 & 7   & 4\\
0.34& 8.1& 5.4 & 3.9 & 3.1& 0.0& 1.2 & 1.5 & 1.8& 27 & 9   & 6\\
0.36& 7.8& 4.8 & 4.2 & 3.0& 0.0& 1.0 & 1.0 & 1.3& 29 & 7   & 5\\
0.38& 7.9& 5.3 & 4.2 & 2.8& 0.0& 0.7 & 0.9 & 1.3& 26 & 9   & 5\\
0.40& 7.7& 5.6 & 5.2 & 2.7& 0.0& 0.3 & 0.6 & 1.1& 33 & 11  & 4   
\end{tabular}
\caption{Comparison results for unconstrained robust (RO), deterministic (DET), and sampling-based (SA10 and SA50) pricing under the nested logit model with homogeneous PSP.}
\label{tb:RO-uncon-Nested}
\end{table}

\subsubsection{Partition-wise Homogeneous PSP. }
We provide comparison results for the case of partition-wise homogeneous PSP considered in Section \ref{sec:uncon-hete}. 
We use the same nested logit model with partition-wise decomposable CPGF 
specified above, i.e., $
G(\bY) = \sum_{n\in [N]} \left(\sum_{i\in C_n} Y_i^{\mu_n} \right)^{\mu/\mu_n}
$,  but the PSP are the same in each nest but different across nests. 
In this context, we know that the robust problem can be converted equivalently into  a convex optimization problem.
On the other hand, for the SA approach, if we select $s_1$ vectors of choice parameters from the uncertainty set, we need to solve $s_1$ convex optimization problems.  

We select $N = 5$ partitions of the same size. For each uncertainty level $\epsilon>0$ and for each partition (or nest) $n\in [N]$, we define a polyhedron uncertainty set  as in Section \ref{sec:construct-uncertainty-set} above.
%Vector $\ba_0 = \{\ba^{n0},\ n\in [N]\}$ is chosen similarly as in the previous section and  $b^{n0}$ is chosen in $[0.5,1.0]$. 
%The purpose here is to test the performance of the robust approach, so any reasonable values can be chosen. 
%$(0.535, 0.635, 0.335, 0.735, 0.454)$.
%Here, to simplify the experiments, we only consider uncertainty sets where the differences between the bounds and the mean values are the same over all the coordinates. 
Similarly to the previous section,  
we first solve the  deterministic problem by bisection with the weighted average  parameters $\widetilde{\bw} = \sum_{k\in [K]} \tau_k \bw^k$ to obtain a solution $\bx^{\DET}$.
Then, for each set $\cA^\epsilon$ we solve the RO problem  by convex optimization to obtain a robust solution $\bx^\RO$. 
We also sample $s_1=10$ and $s_1 = 50$ points from the uncertainty set for the  SA approach.

\begin{figure}[htb]
 \centering
 \includegraphics[width=0.9\textwidth]{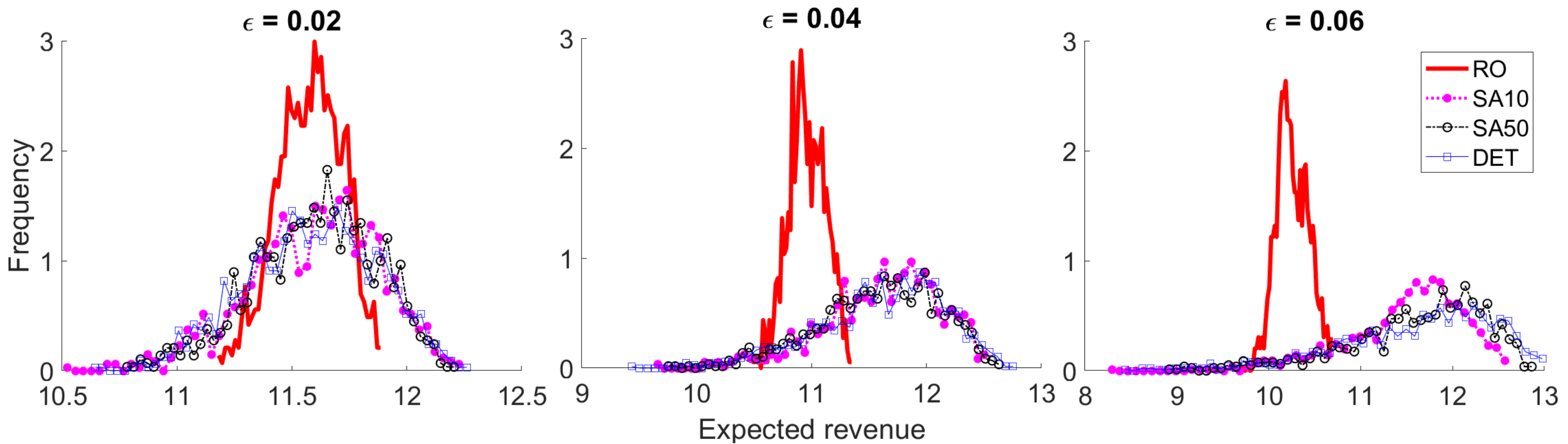}
 \caption{Distributions of the expected revenues under  partition-wise homogeneous PSP.}
 \label{fig:Nested-uncontrained-hete}
\end{figure}

To evaluate the performance of the solutions obtained,  we also sample 1000 points randomly and uniformly  from $\cA^n$, $n\in[N],$ and compute the  expected revenues given by  $\bx^\RO$, $\bx^\textsc{SA10}$, $\bx^\textsc{SA50}$, and $\bx^\DET$. The distributions of the expected revenue over 1000 samples with $\epsilon \in\{0.02,0.04,0.06\}$ are plotted in Figure \ref{fig:Nested-uncontrained-hete}. There is nothing surprising, as similarly to the previous experiments, distributions given by $\bx^\RO$ have small variances, higher peaks, shorter tails and higher worst-case revenues, as compared to those from $\bx^\textsc{SA10}$, $\bx^\textsc{SA50}$ and $\bx^\DET$. In Table \ref{tb:RO-uncon-MNL-nested-hete},we report in detail the average, maximum and worst-case revenues when $\epsilon$ increases from 0.02 to 0.4.
We also see that the RO approach always gives higher  worst-case revenues but lower average revenues, and the SA approaches also provide some protections against low revenues.
However, in this case, the percentile ranks for the DET  and SA approaches are significantly lower (3.45 on average). In particular, we see that there are some instances where the percentile ranks are only 3-\textit{th}, which means that only 3\% of the revenues are lower than the corresponding RO worst-case revenues. Nevertheless , the average revenues given by the SA50 are remarkably higher than those from the RO, especially when $\epsilon$ is large. From this view point, the RO seems too conservative. % in this case (partially homogeneous PSP with rectangular uncertainty sets), and large $\epsilon$ should not be chosen. 

\begin{table}
\centering
\begin{tabular}{r|rrrr|rrrr|rrr}
\multicolumn{1}{l|}{} & \multicolumn{4}{c|}{Average}             & \multicolumn{4}{c|}{Worst}            & \multicolumn{3}{c}{\begin{tabular}[c]{@{}l@{}}Percentile rank~\\of RO worst-case\end{tabular}}     \\ 
\cline{2-12}
\multicolumn{1}{l|}{} & \multicolumn{1}{l}{DET} & \multicolumn{1}{l}{SA10} & \multicolumn{1}{l}{SA50} & \multicolumn{1}{l|}{RO} & \multicolumn{1}{l}{DET} & \multicolumn{1}{l}{SA10} & \multicolumn{1}{l}{SA50} & \multicolumn{1}{l|}{RO} & \multicolumn{1}{l}{DET} & \multicolumn{1}{l}{SA10} & \multicolumn{1}{l}{SA50}  \\ 
\hline
0.02& 11.7  & 11.6   & 11.7   & 9.1   & 6.7   & 6.0& 6.4& 8.7   & 6 & 7  & 8   \\
0.04& 11.6  & 11.6   & 11.6   & 11.6  & 10.7  & 10.7   & 10.7   & 11.1  & 13& 15 & 11  \\
0.06& 11.6  & 11.6   & 11.6   & 10.9  & 9.5   & 9.8& 9.7& 10.4  & 8 & 7  & 7   \\
0.08& 11.7  & 11.5   & 11.7   & 10.3  & 8.0   & 8.9& 9.1& 9.7   & 7 & 4  & 5   \\
0.10 & 11.6  & 11.5   & 11.6   & 9.6   & 6.1   & 8.4& 7.9& 9.1   & 9 & 3  & 4   \\
0.12& 11.8  & 11.7   & 11.4   & 8.6   & 4.6   & 3.7& 7.1& 8.3   & 8 & 12 & 5   \\
0.14& 11.9  & 11.8   & 11.3   & 8.2   & 4.5   & 3.6& 6.9& 7.9   & 7 & 8  & 3   \\
0.16& 11.9  & 11.6   & 11.7   & 7.9   & 5.9   & 5.4& 4.6& 7.5   & 7 & 5  & 4   \\
0.18& 12.0  & 11.8   & 11.3   & 7.6   & 3.2   & 2.4& 5.0& 7.3   & 8 & 7  & 4   \\
0.20 & 12.0  & 11.7   & 11.7   & 7.3   & 1.3   & 1.5& 3.4& 7.0   & 7 & 4  & 5   \\
0.22& 12.0  & 11.7   & 11.3   & 7.1   & 2.7   & 2.8& 5.4& 6.7   & 6 & 5  & 3   \\
0.24& 11.9  & 11.7   & 10.9   & 6.8   & 1.5   & 4.4& 3.3& 6.5   & 8 & 4  & 3   \\
0.26& 12.0  & 11.1   & 11.5   & 6.7   & 0.5   & 5.3& 4.3& 6.4   & 8 & 3  & 3   \\
0.28& 12.1  & 11.1   & 11.9   & 6.6   & 1.3   & 4.7& 0.8& 6.3   & 7 & 3  & 6   \\
0.30 & 12.1  & 11.9   & 10.8   & 6.6   & 2.2   & 3.7& 6.2& 6.2   & 6 & 4  & 2   \\
0.32& 12.2  & 11.4   & 10.6   & 6.5   & 0.6   & 3.0& 5.0& 6.2   & 2 & 2  & 2   \\
0.34& 12.1  & 11.8   & 11.2   & 6.5   & 0.3   & 2.3& 4.6& 6.1   & 9 & 5  & 3   \\
0.36& 12.1  & 11.7   & 11.1   & 6.4   & 0.4   & 0.7& 3.2& 6.0   & 2 & 2  & 2   \\
0.38& 12.2  & 12.2   & 11.8   & 6.4   & 0.4   & 0.6& 1.2& 6.0   & 7 & 10 & 5   \\
0.40 & 12.2  & 11.7   & 10.2   & 6.4   & 0.2   & 4.8& 2.5& 6.0   & 8 & 3  & 5  
\end{tabular}
\caption{Comparison results for unconstrained robust (RO), deterministic (DET), and sampling-based (SA10 and SA50) pricing under a nested logit model with partition-wise homogeneous PSP.}
\label{tb:RO-uncon-MNL-nested-hete}
\end{table}

In summary, our experiments  show {gains}
from our robust models in protecting us from revenues that would be too low. The histograms  given by the robust models have higher peaks, smaller variances, higher worst-case revenues, but lower averages, as compared to their deterministic and sampling-based counterparts. This observation also shows the trade-off in being robust in making pricing decisions when the choice parameters are uncertain, and also consistent with  observations from other relevant studies in the revenue management literature \citep{Rusmevichientong2012robust,Li2019robustCVaR}.

\section{Conclusion}
\label{sec:conclude}
In this paper, we have considered robust versions of the pricing problem under GEV choice models, in  which the choice parameters are not given in advance but lie in an uncertainty set. These robust models are motivated by the fact that uncertainties may occur in the estimation procedure of the choice parameters. We have shown that when the problem is unconstrained  and the PSP are the same over all the products, the robust optimal prices have a constant markup with respect to the product costs and we have shown how to efficiently compute this constant markup by bisection. When the PSP are partition-wise homogeneous and the CPGF and the uncertainty set are also partition-wise separable, we have shown that the robust problem can be converted equivalently into a reduced optimization program, and the reduce problem can be solved conveniently by convex optimization. 

We have also considered the  pricing problem with over-expected-revenue-penalties as an alternative  to the constrained pricing problem.
%We have shown that, when the penalty parameters goes to infinity, the penalty term converges to zero and the optimal value converges to the expected revenue given by the constrained pricing problem. {Since there may be no fixed prices under which the purchase probabilities always satisfy the expected sale constraints when the choice parameters vary,
%the robust version of the pricing problem under over-expected-sale penalties is more appropriate to use in this context, as opposed to the robust constrained version.}
We have shown that under the same assumptions as in the case of partition-wise homogeneous PSP , the robust problem can be converted equivalently into a reduced one, which can be further solved by convex optimization. Experimental results based on the nested logit model have shown the {advantages}
of our robust model in providing protection against bad-case revenues. 
In future research, it would be interesting to look at distributionally robust versions of the pricing problem, which may help provide less conservative robust solutions as compared to the standard robust optimization approaches. We are also interested in robust approaches for the joint assortment and pricing problem under  GEV choice models.%, and a distributionally robust approach to mitigate the conservativeness of the current robust models.  

% Acknowledgments here
\ACKNOWLEDGMENT{This research is supported by the National Research Foundation, Prime Minister’s Office, Singapore under its Campus for Research Excellence and Technological Enterprise
(CREATE) program, Singapore-MIT Alliance for Research and Technology (SMART) Future Urban
Mobility (FM) IRG.}
%We thank an Associate Editor and two anonymous referees whose comments substantially helped us improve the previous version of the paper.}

\bibliographystyle{plainnat}
\bibliography{RefDRO,RefMathProg,RefRM,refsPricing}

\begin{thebibliography}{37}
\providecommand{\natexlab}[1]{#1}
\providecommand{\url}[1]{\texttt{#1}}
\expandafter\ifx\csname urlstyle\endcsname\relax
  \providecommand{\doi}[1]{doi: #1}\else
  \providecommand{\doi}{doi: \begingroup \urlstyle{rm}\Url}\fi

\bibitem[Ben-Akiva(1973)]{Ben1973structure}
Moshe~E Ben-Akiva.
\newblock \emph{Structure of passenger travel demand models.}
\newblock PhD thesis, Massachusetts Institute of Technology, 1973.

\bibitem[Ben-Akiva et~al.(1985)Ben-Akiva, Lerman, and Lerman]{Ben1985discrete}
Moshe~E Ben-Akiva, Steven~R Lerman, and Steven~R Lerman.
\newblock \emph{Discrete choice analysis: theory and application to travel
  demand}, volume~9.
\newblock MIT press, 1985.

\bibitem[Ben-Tal and Nemirovski(1998)]{Ben1998robustconvex}
Aharon Ben-Tal and Arkadi Nemirovski.
\newblock Robust convex optimization.
\newblock \emph{Mathematics of operations research}, 23\penalty0 (4):\penalty0
  769--805, 1998.

\bibitem[Ben-Tal and Nemirovski(2000)]{Ben2000robust}
Aharon Ben-Tal and Arkadi Nemirovski.
\newblock Robust solutions of linear programming problems contaminated with
  uncertain data.
\newblock \emph{Mathematical programming}, 88\penalty0 (3):\penalty0 411--424,
  2000.

\bibitem[Ben-Tal et~al.(2006)Ben-Tal, Boyd, and Nemirovski]{Ben2006extending}
Aharon Ben-Tal, Stephen Boyd, and Arkadi Nemirovski.
\newblock Extending scope of robust optimization: Comprehensive robust
  counterparts of uncertain problems.
\newblock \emph{Mathematical Programming}, 107\penalty0 (1-2):\penalty0 63--89,
  2006.

\bibitem[Bhat(1998)]{Bhat1998accommodating}
Chandra~R Bhat.
\newblock Accommodating variations in responsiveness to level-of-service
  measures in travel mode choice modeling.
\newblock \emph{Transportation Research Part A: Policy and Practice},
  32\penalty0 (7):\penalty0 495--507, 1998.

\bibitem[Daly and Bierlaire(2006)]{Daly2006general}
Andrew Daly and Michel Bierlaire.
\newblock A general and operational representation of generalised extreme value
  models.
\newblock \emph{Transportation Research Part B: Methodological}, 40\penalty0
  (4):\penalty0 285--305, 2006.

\bibitem[De~Klerk(2006)]{DeKlerk2006aspects}
Etienne De~Klerk.
\newblock \emph{Aspects of semidefinite programming: interior point algorithms
  and selected applications}, volume~65.
\newblock Springer Science \& Business Media, 2006.

\bibitem[Dong et~al.(2009)Dong, Kouvelis, and Tian]{Dong2009dynamic}
Lingxiu Dong, Panos Kouvelis, and Zhongjun Tian.
\newblock Dynamic pricing and inventory control of substitute products.
\newblock \emph{Manufacturing \& Service Operations Management}, 11\penalty0
  (2):\penalty0 317--339, 2009.

\bibitem[Fosgerau et~al.(2013)Fosgerau, McFadden, and
  Bierlaire]{Fosgerau2013GPGF}
Mogens Fosgerau, Daniel McFadden, and Michel Bierlaire.
\newblock Choice probability generating functions.
\newblock \emph{Journal of Choice Modelling}, 8:\penalty0 1--18, 2013.

\bibitem[Gallego and Hu(2014)]{Gallego2014dynamicPricing}
Guillermo Gallego and Ming Hu.
\newblock Dynamic pricing of perishable assets under competition.
\newblock \emph{Management Science}, 60\penalty0 (5):\penalty0 1241--1259,
  2014.

\bibitem[Gallego and Van~Ryzin(1997)]{Gallego1997multiproductPricing}
Guillermo Gallego and Garrett Van~Ryzin.
\newblock A multiproduct dynamic pricing problem and its applications to
  network yield management.
\newblock \emph{Operations research}, 45\penalty0 (1):\penalty0 24--41, 1997.

\bibitem[Gallego and Wang(2014)]{Gallego2014multiproduct}
Guillermo Gallego and Ruxian Wang.
\newblock Multiproduct price optimization and competition under the nested
  logit model with product-differentiated price sensitivities.
\newblock \emph{Operations Research}, 62\penalty0 (2):\penalty0 450--461, 2014.

\bibitem[Hogan(1973)]{Hogan1973}
William~W Hogan.
\newblock Point-to-set maps in mathematical programming.
\newblock \emph{SIAM review}, 15\penalty0 (3):\penalty0 591--603, 1973.

\bibitem[Hopp and Xu(2005)]{Hopp2005}
Wallace~J Hopp and Xiaowei Xu.
\newblock Product line selection and pricing with modularity in design.
\newblock \emph{Manufacturing \& Service Operations Management}, 7\penalty0
  (3):\penalty0 172--187, 2005.

\bibitem[Huh and Li(2015)]{Huh2015pricing}
Woonghee~Tim Huh and Hongmin Li.
\newblock Pricing under the nested attraction model with a multistage choice
  structure.
\newblock \emph{Operations Research}, 63\penalty0 (4):\penalty0 840--850, 2015.

\bibitem[Keller(2013)]{Keller2013}
Philipp~Wilhelm Keller.
\newblock \emph{Tractable multi-product pricing under discrete choice models}.
\newblock PhD thesis, Massachusetts Institute of Technology, 2013.

\bibitem[Koppelman and Wen(2000)]{koppelman2000paired}
Frank~S Koppelman and Chieh-Hua Wen.
\newblock The paired combinatorial logit model: properties, estimation and
  application.
\newblock \emph{Transportation Research Part B: Methodological}, 34\penalty0
  (2):\penalty0 75--89, 2000.

\bibitem[Li et~al.(2015)Li, Rusmevichientong, and Topaloglu]{li2015d_nested}
Guang Li, Paat Rusmevichientong, and Huseyin Topaloglu.
\newblock The d-level nested logit model: Assortment and price optimization
  problems.
\newblock \emph{Operations Research}, 63\penalty0 (2):\penalty0 325--342, 2015.

\bibitem[Li and Huh(2011)]{Li2011pricing}
Hongmin Li and Woonghee~Tim Huh.
\newblock Pricing multiple products with the multinomial logit and nested logit
  models: Concavity and implications.
\newblock \emph{Manufacturing \& Service Operations Management}, 13\penalty0
  (4):\penalty0 549--563, 2011.

\bibitem[Li et~al.(2018)Li, Webster, Mason, and Kempf]{Li2018mixed_pricing}
Hongmin Li, Scott Webster, Nicholas Mason, and Karl Kempf.
\newblock Product-line pricing under discrete mixed multinomial logit demand.
\newblock \emph{Manufacturing \& Service Operations Management}, 21\penalty0
  (1):\penalty0 14--28, 2018.

\bibitem[Li et~al.(2019)Li, Webster, Mason, and Kempf]{Li2019productMMNL}
Hongmin Li, Scott Webster, Nicholas Mason, and Karl Kempf.
\newblock Product-line pricing under discrete mixed multinomial logit demand:
  Winner—2017 m\&som practice-based research competition.
\newblock \emph{Manufacturing \& Service Operations Management}, 21\penalty0
  (1):\penalty0 14--28, 2019.

\bibitem[Li and Ke(2019)]{Li2019robustCVaR}
Xiaolong Li and Jiannan Ke.
\newblock Robust assortment optimization using worst-case cvar under the
  multinomial logit model.
\newblock \emph{Operations Research Letters}, 47\penalty0 (5):\penalty0
  452--457, 2019.

\bibitem[Mai et~al.(2017)Mai, Frejinger, Fosgerau, and Bastin]{Mai2017dynamic}
Tien Mai, Emma Frejinger, Mogens Fosgerau, and Fabian Bastin.
\newblock A dynamic programming approach for quickly estimating large
  network-based mev models.
\newblock \emph{Transportation Research Part B: Methodological}, 98:\penalty0
  179--197, 2017.

\bibitem[McFadden(1978)]{Mcfadden1978modeling}
Daniel McFadden.
\newblock Modeling the choice of residential location.
\newblock \emph{Transportation Research Record}, \penalty0 (673), 1978.

\bibitem[McFadden(1980)]{Mcfadden1980econometric}
Daniel McFadden.
\newblock Econometric models for probabilistic choice among products.
\newblock \emph{Journal of Business}, pages S13--S29, 1980.

\bibitem[McFadden and Train(2000)]{Mcfadden2000mixed}
Daniel McFadden and Kenneth Train.
\newblock Mixed mnl models for discrete response.
\newblock \emph{Journal of applied Econometrics}, 15\penalty0 (5):\penalty0
  447--470, 2000.

\bibitem[Rusmevichientong and Topaloglu(2012)]{Rusmevichientong2012robust}
Paat Rusmevichientong and Huseyin Topaloglu.
\newblock Robust assortment optimization in revenue management under the
  multinomial logit choice model.
\newblock \emph{Operations Research}, 60\penalty0 (4):\penalty0 865--882, 2012.

\bibitem[Shapiro(2018)]{Shapiro2018tutorialDRO}
Alexander Shapiro.
\newblock Tutorial on risk neutral, distributionally robust and risk averse
  multistage stochastic programming.
\newblock \emph{Optimization Online http://www. optimization-online.
  org/DB\_HTML/2018/02/6455. html}, 2018.

\bibitem[Small(1987)]{Small1987discrete}
Kenneth~A Small.
\newblock A discrete choice model for ordered alternatives.
\newblock \emph{Econometrica: Journal of the Econometric Society}, pages
  409--424, 1987.

\bibitem[Song and Xue(2007)]{Song2007demand}
Jing-Sheng Song and Zhengliang Xue.
\newblock Demand management and inventory control for substitutable products.
\newblock \emph{Working paper}, 2007.

\bibitem[Talluri and Van~Ryzin(2004)]{Talluri2004RM}
Kalyan Talluri and Garrett Van~Ryzin.
\newblock Revenue management under a general discrete choice model of consumer
  behavior.
\newblock \emph{Management Science}, 50\penalty0 (1):\penalty0 15--33, 2004.

\bibitem[Vovsha and Bekhor(1998)]{Vovsha1998link}
Peter Vovsha and Shlomo Bekhor.
\newblock Link-nested logit model of route choice: overcoming route overlapping
  problem.
\newblock \emph{Transportation research record}, 1645\penalty0 (1):\penalty0
  133--142, 1998.

\bibitem[Wen and Koppelman(2001)]{Wen2001generalized}
Chieh-Hua Wen and Frank~S Koppelman.
\newblock The generalized nested logit model.
\newblock \emph{Transportation Research Part B: Methodological}, 35\penalty0
  (7):\penalty0 627--641, 2001.

\bibitem[Whelan et~al.(2002)Whelan, Batley, Fowkes, and
  Daly]{whelan2002flexible}
GRTA Whelan, R~Batley, T~Fowkes, and A~Daly.
\newblock Flexible models for analyzing route and departure time choice.
\newblock \emph{Publication of: Association for European Transport}, 2002.

\bibitem[Zhang and Lu(2013)]{Zhang2013assessing}
Dan Zhang and Zhaosong Lu.
\newblock Assessing the value of dynamic pricing in network revenue management.
\newblock \emph{INFORMS Journal on Computing}, 25\penalty0 (1):\penalty0
  102--115, 2013.

\bibitem[Zhang et~al.(2018)Zhang, Rusmevichientong, and
  Topaloglu]{zhang2018multiproduct}
Heng Zhang, Paat Rusmevichientong, and Huseyin Topaloglu.
\newblock Multiproduct pricing under the generalized extreme value models with
  homogeneous price sensitivity parameters.
\newblock \emph{Operations Research}, 66\penalty0 (6):\penalty0 1559--1570,
  2018.

\end{thebibliography}

% Appendix here
% Options are (1) APPENDIX (with or without general title) or
% (2) APPENDICES (if it has more than one unrelated sections)
% Outcomment the appropriate case if necessary
%
% \begin{APPENDIX}{<Title of the Appendix>}
% \end{APPENDIX}
%
%or
%
\clearpage
 \begin{APPENDICES}
 %\section{<Title of Section A>}
 %\section{<Title of Section B>}
% etc

%%%%%%%%%%%%%%%%%%%%%%%%%%%%%
%
%  PROOF .....
%
%%%%%%%%%%%%%%%%%%%%%%%%%%%%

\section{Proofs}
\label{apd:sec-proofs}
This section provides some detailed proofs of the claims presented in the main part of the paper.

%%%%%%%%%%%%%%%%%%
%%%%%%%%%%%%%%%%%%

\subsection{Proof of Proposition \ref{prop:G-strictly-convex}}\label{proof:G-strictly-convex}

First, we consider function $f^G(\bs): \bbR^m \rightarrow \bbR_+$
\[
f^G(\bs) = G(Y_1,\ldots,Y_m), \text{ where } Y_i = e^{s_i},\ \forall i\in\cV
\]
We will prove that $f^G(\bs)$ is convex. Taking the first and second derivatives of $f^G(\bs)$  we obtain
\[
\frac{\partial f^G(\bs)}{\partial s_i} = \partial G_i(\bY)Y_i, 
\]
and
\[
\begin{aligned}
 \frac{\partial^2 f^G(\bs)}{\partial s_i\partial s_i} &=  \partial G_{ii}(\bY)Y^2_i +  \partial G_{i}(\bY)Y_i, \\
  \frac{\partial^2 f^G(\bs)}{\partial s_i\partial s_j} &=  \partial G_{ij}(\bY)Y_iY_j.  \\
\end{aligned}
\]
So we have
\[
\nabla^2 f^G(\bs) = \dg(\bY)\nabla^2 G(\bY)\dg(\bY) + \dg (\nabla G(\bY) \circ \bY), 
\]
where $\dg(\bY)$ is the square diagonal matrix with the elements of vector $\bY$ on the main diagonal.
The second term $\dg (\nabla G(\bY) \circ \bY) $ is always positive definite, \mtien{where $\circ$ is the element-by-element operator}. Moreover, $\dg(\bY)\nabla^2 G(\bY)\dg(\bY)$  is symmetric
and its $(i, j)$-th component is given by $Y_i\partial G_{ij}(\bY) Y_j$. For $i\neq j$, we have $\partial G_{ij} (\bY) \leq 0$ by the property of the GEV-CPGF $G$, so all off-diagonal entries of the matrix are non-positive. In addition,  $\sum_{j\in\cV} Y_j\partial G_{ij} (\bY) = 0$, so that each
row of the matrix sums to zero. Thus,  $\dg(\bY)\nabla^2 G(\bY)\dg(\bY)$ is positive semi-definite \citep[see Theorem A.6 in][]{DeKlerk2006aspects}. %\pj{not sure I understand why this sentence is here, also note that it isnt a complete sentence}\tm{I corrected the sentence}. 
So, $\nabla^2 f^G(\bs)$ is positive definite, or equivalently, $f^G(\bs)$ is strictly convex in $\bs$. This lead to the following inequality, for all $\bs^1,\bs^2\in\bbR^m$ and $\lambda \in (0,1)$
\[ \lambda f^G(\bs^1) + \lambda f^G(\bs^2) > f^G(\lambda \bs^1 + (1-\lambda )\bs^2).
\]
For all $(\ba^1,\bb^1), (\ba^2,\bb^2)\in \cA$, replace $s^1_i$ by $a^1_i - b^1_i(z+c_i)$ and $s^2_i$ by $a^2_i - b^2_i(z+c_i)$ we have
\[
  \lambda G(\bY|z,\ba^1,\bb^1) + \lambda G(\bY|z,\ba^2,\bb^2) > G(\bY|z,\lambda\ba^1 + (1-\lambda) \ba^2,\lambda \bb^1+ (1-\lambda)\bb^2),\ \forall \lambda \in (0,1)
\]
which means that  $G(\bY|z,\ba,\bb)$ is strictly convex in $\ba,\bb$. 

The continuity of $(\ba^*(z),\bb^*(z)$ is a direct result from the convexity of $G(\bY|z,\ba,\bb)$ and  the  Corollary 8.2 of  \cite{Hogan1973}.
This completes the proof.

%%%%%%%%%%%%%%%%%%%
%%%%%%%%%%%%%%%%%%
\subsection{Proof of Lemma \ref{lemma:G-bounded}}
\label{apdx:proof-Lemma-G-bounded}
Consider $f^G(\bs) = G(Y_1,\ldots,Y_m)$, where  $Y_i = e^{s_i},\ \forall i=1,\ldots,m$. Taking the derivative of $f^G(\bs)$ w.r.t. $s_i$ we have
\[
\frac{\partial f^G(\bs)}{\partial s_i} = \partial G_i(\bY)Y_i  \geq 0
\]
So, $f^G(\bs)$ is monotonic in every coordinate, meaning that given  any $\bs,\bs_0 \in \bbR^m$, $\bs \succeq s_0$, we have $f^G(\bs) \geq f^G(\bs_0)$. Moreover, it is clear  that  
\[
 \underline{\ba} - \overline{\bb}\circ (\bc+z\bbe) \preceq \ba - \bb \circ (\bc+z\bbe) \preceq \overline{\ba} - \underline{\bb}\circ (\bc+z\bbe), \ \forall z\in\bbR_+,(\ba,\bb) \in \cA.
\]
So, we obtain the following inequality
\[
G(\bY|z,\underline{\ba},\overline{\bb}) \leq  G(\bY|z,\ba,\bb) \leq G(\bY|z,\overline{\ba},\underline{\bb}), \ \forall (\ba,\bb) \in\cA,
\]
which completes the proof.

\subsection{Proof of Theorem \ref{theor;RO-hete-convexness-simplified-model}}
\label{proof:RO-hete-convexness-simplified-model}
The first lemma shows that $\underline{\cG}^n$ is invertible.
\begin{lemma}
\label{lemma:reduce-concave-lm1}
Given $n\in[N]$, for any $\alpha > 0$, there is a unique $z_n\in \bbR$ such that $\underline{\cG}^n(z_n) = \alpha$. 
\end{lemma}
\proof{Proof:}
From the properties of CPGF  (Remark \ref{propert:GEV-CPGF}), we see that function $G^n(\bY^n|z_n,\ba^n,\bb^n)$ is strictly monotonic-decreasing, so $\underline{\cG}^n(z_n)$ is also strictly monotonic-decreasing. Moreover, 
we have $\lim_{z_n \rightarrow +\infty }G^n(\bY^n|z_n,\ba^n,\bb^n) = 0$ and $\lim_{z_n \rightarrow -\infty }G^n(\bY^n|z_n,\ba^n,\bb^n) = \infty$. Thus,  $G^n(\bY^n|z_n,\ba^n,\bb^n)$ spans all over the set $\bbR_+$ when $z_n$ varies. On the other hand, $\cA^n$ is bounded, we will also have $\lim_{z_n \rightarrow +\infty }\underline{\cG}^n(z_n) = 0$ and $\lim_{z_n \rightarrow -\infty }\underline{\cG}^n(z_n) = \infty$. Since $\underline{\cG}^n(z_n) = 0$ is continuous and strictly  monotonic-decreasing, we easily obtain the desired result.  
\endproof
The above lemma allows us to define the inverse function of $\underline{\cG}^n(\cdot)$ as $(\underline{\cG}^n)^{-1}(\alpha):\bbR_+ \rightarrow \bbR$ such that $\underline{\cG}^n((\underline{\cG}^n)^{-1}(\alpha)) = \alpha$. As shown above, this function can be computed by bisection. Lemma \ref{lemma:zp-unique-deter} below shows how to to identify $\bz(\bp^G)$. 
\begin{lemma}
\label{lemma:zp-unique-deter}
Given any $\bp^G\in \cP^G$, $\bz(\bp^G)$ can be uniquely computed as
\[
\bz(\bp^G)_n = (\underline{\cG}^n)^{-1}\left(\frac{p^G_n}{1-\sum_{l\in[N]}  p^G_l}\right)
\]
\end{lemma}
\proof{Proof:}
From \eqref{eq:pG-z} we see that
\[
\sum_{l\in [N]} p^G_l =  \frac{\sum_{l\in[N]}\underline{\cG}^l(z_l)}{1 + \sum_{l\in[N]} \underline{\cG}^l(z_l)},
\]
So we have
\[
1 + \sum_{l\in[N]} \underline{\cG}^l(z_l) = \frac{1}{1-\sum_{l\in[N]}  p^G_l}.
\]
Thus, for any $n\in[N]$ 
\[
\underline{\cG}^l(z_n) = \frac{p^G_n}{1 + \sum_{l\in[N]} \underline{\cG}^l(z_l)} = \frac{p^G_n}{1-\sum_{l\in[N]}  p^G_l},
\]
which directly leads to the desired result. 
\endproof
We now move to the second claim of Theorem \ref{theor;RO-hete-convexness-simplified-model}, i.e., the convexity of $\cW(\bz(\bp^G))$. To support the proof, let use consider a deterministic version of \eqref{prob:RO-hete-simlified} in which all the choice parameter are given
$
\widetilde{\cW}(\widetilde{\bz}(\bp^G|\ba,\bb)) = \sum_{n\in[N]} \widetilde{\bz}(\bp^G|\ba,\bb)_n p^G_n,
$
where $\widetilde{\bz}(\bp^G|\ba,\bb) \in \bbR^N$ are a vector of constant-markups that archive vector $\bp^G$  as 
\begin{equation}
\label{eq:pG-z-fixed-ab}
p^{G}_n =  \frac{G^n(\bY^n| \widetilde{z}_n,\ba^n,\bb^n)}{1 + \sum_{l\in[N]} G^l(\bY^l| \widetilde{z}_l,\ba^l,\bb^l)},\;\forall n\in [N].
\end{equation}
\begin{lemma}
\label{lemma:Ro-hete-deter-simplified}
Given any $\bp^G \in \cP^G$, $\widetilde{\bz}(\bp^G|\ba,\bb)$ can be uniquely determined by solving a strictly convex optimization problem. Moreover, $\widetilde{\cW}(\widetilde{\bz}(\bp^G|\ba,\bb))$ is strictly concave in $\bp^G$.
\end{lemma}
\proof{Proof:}
Let $\Theta(\bz):\bbR^N\rightarrow\bbR^N$ such that $\Theta(\bz)_n = G^n(\bY^n| z_n)\Big/\left(1+\sum_{j\in[N]} G^j(\bY^j|z_j)\right)$, where $G^n(\bY^n|z_n) = G^n(\bY^n|z_n,{\ba},{\bb})$
but we omit the choice parameters $({\ba},{\bb})$ for notational simplicity.
We also denote by $\widetilde{\bb}$  a vector of size $N$ with entries $\widetilde{b}_n = \me{{\bb}^n}$. 
Consider the problem
\begin{equation}
\label{prob:ro-pen-gev-sub-eq1}
\min_{\bz \in \bbR^N}\left\{\ln\left(1+ \sum_{n\in[N]}G^n(\bY^n|z_n) \right) + \sum_{n\in[N]} p^G_n \widetilde{b}_n z_n \right\} 
\end{equation}
and we now show that \eqref{prob:ro-pen-gev-sub-eq1} is a strictly convex optimization problem and solving it will yield a solution $\bz^*$ such that $\Theta(\bz^*) = \bp^G$.
Note that
%The lemma below
% show that given any $\bp^G$, $\bz(\bp^G)$ can be uniquely determined by solving a strictly convex optimization problem.  
the structure of the  problem presented in this lemma is slightly different with those considered in Theorem 4.1 in \cite{zhang2018multiproduct}, so even though the proof of the lemma is quite similar, we provide its own proof for the sake of self-contained. 
% \begin{lemma}
%Given $\bp^G\in\Delta^N$, there is a unique vector $\bz(\bp^G)\in \bbR^N$ such that $\Theta(\bz(\bp^G)) = \bp^G$, and this $\bz(\bp^G)$ is the unique solution to the convex optimization problem.
%\end{lemma}
%\proof{Proof:}
To prove that \eqref{prob:ro-pen-gev-sub-eq1} is a strictly convex optimization problem, we will show that $\nabla^2\cQ(\bz)$ is  a positive definite matrix, where $\cQ(\bz) =  \ln\left(1+ \sum_{n\in[N]}G^n(\bY^n|z_n) \right)$. To simplify the proof and make use of previous results, let us denote $\bu(\bz):\bbR^N \rightarrow \bbR^N$ such that $u(\bz)_n = -\ln G^n(\bY^n|z_n)/{\widetilde{b}_n}$.
with this definition we have
\[
\frac{\partial u(\bz)_n}{ \partial z_n} = \frac{\sum_{i\in\cV_n} \partial G^n_i(\bY^n|z_n)Y_i\widetilde{b}_n}{ G^n(\bY^n|z_n) \widetilde{b}_n} = 1.
\]
The objective function now can be written as
\[
\cQ(\bu(\bz)) = \ln\left(1+ \sum_{n\in[N]} \exp({-u(\bz)_n}). \right) 
\]
Taking the derivative of $\cQ$ with respect to $z_n$ we obtain
\[
\begin{aligned}
\frac{\partial \cQ(\bu(\bz))}{\partial z_n} &= \left.\frac{\partial \cQ(\bu)}{\partial u_n} \right\rvert_{{\bu = \bu(\bz)}} \frac{\partial u(\bz)_n}{\partial z_n} = \left.\frac{\partial \cQ(\bu)}{\partial u_n} \right\rvert_{{\bu = \bu(\bz)}}.
\end{aligned}
\]
And if we take the second derivative with respective to $z_n,z_k$, $n,k\in[N]$ we get
\[
\frac{\partial^2 \cQ(\bu(\bz))}{\partial z_n \partial z_k} =  \left.\frac{\partial^2 \cQ(\bu)}{\partial u_n \partial u_k} \right\rvert_{{\bu = \bu(\bz)}},
\]
or equivalently $\nabla^2\cQ(\bz) =  \nabla^2_\bu \cQ(\bu)$, where $\bu = \bu(\bz)$. 
Moreover, $\cQ(\bu)$ is just a special objective function under the MNL model with $N$ products and all the PSP are equal to 1. As a result, $\nabla^2_\bu \cQ(\bu)$ is positive definite \cite[see Theorem 4.1][]{zhang2018multiproduct}, so $\nabla^2\cQ(\bz)$ is also positive definite, as desired. 

Now we know that \eqref{prob:ro-pen-gev-sub-eq1} is strictly convex, so it yields a unique solution. Moreover, one can  show that \eqref{prob:ro-pen-gev-sub-eq1} have finite optimal solutions. For any $n\in [N]$, taking the derivative of $\cQ(\bz)$ with respect to $z_n$ and set it to zero we obtain
\[
\frac{\sum_{n\in[N]} \sum_{i\in\cV_n} -\partial G^n_i(\bY^n|z_n) Y_i \widetilde{b}_n }{1+G(\bY)} = p_n^G\widetilde{b}_n,
\]
or equivalently, $\bp^G = \Theta(\bz)$. 
So, if $\bz(\bp^G)$ is the unique solution to \eqref{prob:ro-pen-gev-sub-eq1},  we always have $\bp^G = \Theta(\bz(\bp^G))$ as desired.

Next, we will show that $\widetilde{\cW}(\widetilde{\bz}(\bp^G|\ba,\bb))$ is a strictly concave function of $\bp^G$.
We also omit the choice parameters for notational convenience and denote 
$\widetilde{\cW}(\widetilde{\bz}) =  \sum_{n\in[N]} \widetilde{z}_n p^G_n$. 
We first see that $\bp^G$ is also a choice probability vector given by a MNL model with $N$ products with the utility vector $-\bu(\widetilde{\bz}(\bp^G))\circ \widetilde{\bb}$. So, if we denote $\bu'(\bp^G)$ be a mapping from $\bbR^N$ to $\bbR^N$ such that $p^G_n = \exp(- \widetilde{b}_n u'(\bp^G)_n )\Big/ \left(\sum_{n\in[N]} \exp(-\widetilde{b}_n u'(\bp^G)_n)\right)$, then we have $\bu(\widetilde{\bz}(\bp^G)) = \bu'(\bp^G)$
\[
\begin{aligned}
\frac{\partial  \widetilde{\cW}(\widetilde{\bz}(\bp^G))}{\partial p_n^G} &= \widetilde{z}(\bp^G)_n + \sum_{j\in[n]}\frac{p^G_j\partial \widetilde{z}(\bp^G)_j}{\partial p_n^G} \\
 &= \widetilde{z}(\bp^G)_n + \sum_{j\in[N]}p^G_j\frac{\partial \widetilde{z}(\bp^G)_j}{\partial u(\widetilde{\bz}(\bp^G))_j}\frac{\partial u(\widetilde{\bz}(\bp^G))_j}{\partial p_n^G}\\
 & =  \widetilde{z}(\bp^G)_n +  \sum_{j\in[N]}{p^G_j}\frac{\partial u'(\bp^G)_j}{\partial p_n^G}
\end{aligned}
\]
and
\[
\begin{aligned}
\frac{\partial^2  \widetilde{\cW}(\widetilde{\bz}(\bp^G))}{\partial p_n^G \partial p_k^G}
 & =  \frac{\partial \widetilde{z}(\bp^G)_n}{\partial p_k^G} +  \sum_{j\in[N]}{p^G_j}\frac{\partial^2 u'(\bp^G)_j}{\partial p_n^G \partial p_k^G} \\
 & =  \frac{\partial \widetilde{z}(\bp^G)_n}{\partial u(\widetilde{\bz}(\bp^G))_n} \frac{\partial u(\widetilde{\bz}(\bp^G))_n}{\partial p_k^G} +  \sum_{j\in[N]}{p^G_j}\frac{\partial^2 u'(\bp^G)_j}{\partial p_n^G \partial p_k^G} \\
 &= \frac{\partial u'(\bp^G)_n}{\partial p_k^G} +  \sum_{j\in[N]}{p^G_j}\frac{\partial^2 u'(\bp^G)_j}{\partial p_n^G \partial p_k^G}
\end{aligned}
\]
Moreover, if we denote $  \widetilde{\cW}'(\bp^G) = \sum_{n\in [N]} u'(\bp^G)_np^G_n$, we also have 
\[
\frac{\partial^2  \widetilde{\cW}'(\bp^G)}{\partial p_n^G \partial p_k^G} = \frac{\partial u'(\bp^G)_n}{\partial p_k^G} +  \sum_{j\in[N]}{p^G_j}\frac{\partial^2 u'(\bp^G)_j}{\partial p_n^G \partial p_k^G},\ \forall n,k\in[N].
\]
So, $\nabla^2  \widetilde{\cW}(\widetilde{\bz}(\bp^G)) = \nabla^2  \widetilde{\cW}'(\bp^G)$. We also see that $  \widetilde{\cW}'(\bp^G)$ is the expected revenue function (as a function of the purchase probabilities $\bp^G$) where there are $N$ products, the choice model is MNL, the PSP are $\widetilde{\bb}$ and the utility vector is $-u'(\bp^G)\circ \widetilde{\bb}$ ($\circ$ is the \textit{dot product}). So we know that $\nabla^2  \widetilde{\cW}'(\bp^G)$ is negative definite \citep{zhang2018multiproduct}, so $\widetilde{\cW}(\widetilde{\bz}(\bp^G))$ is  strictly concave in $\bp^G$.
This completes the proof.
\endproof

We now make a connection between $\bz(\bp^G)$ defined in \eqref{eq:pG-z} and \eqref{eq:pG-z-fixed-ab}. This is crucial to show the concavity of ${\cW}(\bz(\bp^G))$. We have the following lemma.
\begin{lemma}
\label{lemma:Ro-hete-dev-z}
Given any $\bp^G\in\cP^G$, we have the following equalities
\begin{itemize}
\item [(i)] $\bz(\bp^G) = \widetilde{\bz}(\bp^G|\ba^*(\bz), \bb^*(\bz))$
\item [(ii)] The first and second-order derivatives of $z(\bp^G)_n$, $n\in[N]$
\begin{align}
\frac{\partial z(\bp^G)_n }{\partial p^G_l} &=  \left.\frac{\partial \widetilde{z}(\bp^G|\ba,\bb)_n}{\partial p^G_l} \right\rvert_{{\substack{\ba = \ba^*(\bz)\\\bb= \bb^*(\bz)}}} \qquad \forall n\in[N] \nonumber \\
\frac{\partial^2 z(\bp^G)_n }{\partial p^G_l\partial p^G_k} &=  \left.\frac{\partial^2 \widetilde{z}(\bp^G|\ba,\bb)_n}{\partial p^G_l\partial p^G_k} \right\rvert_{{\substack{\ba = \ba^*(\bz)\\\bb= \bb^*(\bz)}}} \qquad \forall l,k\in[N]\nonumber
\end{align}
where  $(\ba^*(\bz), \bb^*(\bz))  = \{(\ba^{n*}(\bz), \bb^{n*}(\bz))|\; n\in [N]\}$ and $(\ba^{n*}(z_n), \bb^{n*}(z_n))=$ $ \text{argmin}_{(\ba^{n}, \bb^{n}) \in\cA^n}  G^n(\bY^n| z_n,\ba^n,\bb^n)$.
\end{itemize}
\end{lemma}
%\proof{Proof of Lemma \ref{lemma:Ro-hete-dev-z}:}
\proof{Proof:}
First, we see that since the uncertainty set  $\cA^n$ does not depend on $z_n$, the derivatives of $\underline{\cG}^n(z_n)$ can be computed as
\begin{align}
\frac{\partial \underline{\cG}^n(z_n)}{\partial z_n} = \left.\frac{\partial  G^n(\bY^n|z_n,\ba^n,\bb^n)}{\partial z_n} \right\rvert_{{\ba^{n} = \ba^{n*}(z_n);\:  \bb^{n} = \bb^{n*}(z_n)}} \label{eq:grad-cG} \\
\frac{\partial^2 \underline{\cG}^n(z_n)}{\partial^2 z_n} = \left.\frac{\partial^2  G^n(\bY^n|z_n,\ba^n,\bb^n)}{\partial^2 z_n} \right\rvert_{{\ba^{n} = \ba^{n*}(z_n);\:  \bb^{n} = \bb^{n*}(z_n),}}. \nonumber
\end{align}

For (i), we know that $\bz(\bp^G)$  is a unique solution to the following system 
\begin{equation}
\label{eq:lemma:Ro-hete-dev-z-eq1}
\underline{\cG}^n(z_n) = \frac{p^G_n}{1-\sum_{l\in[N]} p^G_l}, \; \forall n\in[N],
\end{equation}
and $ \widetilde{\bz}(\bp^G|\ba^*(\bz), \bb^*(\bz))$ is a unique solution to
\begin{equation}
\label{eq:lemma:Ro-hete-dev-z-eq2}
G^n(\bY^n| \widetilde{z}_n, \ba^*(\bz), \bb^*(\bz)) = \frac{p^G_n}{1-\sum_{l\in[N]} p^G_l}, \; \forall n\in[N],
\end{equation}
and note that $G^n(\bY^n| {z}_n, \ba^*(\bz), \bb^*(\bz)) = \underline{\cG}^n(z_n)$. This leads to the desired equality.

For (ii), we take the derivative of \eqref{eq:lemma:Ro-hete-dev-z-eq1} with respect to $p^G_j$, $j\in[N]$, and obtain
\begin{align}
\frac{\partial h^n(\bp^G)}{\partial p^G_l} &= \frac{\partial \underline{\cG}^n(z(\bp^G)_n)}{\partial p^G_l} =  \frac{\partial \underline{\cG}^n(z(\bp^G)_n)}{\partial z(\bp^G)_n} \frac{\partial z(\bp^G)_n}{\partial p^G_l}\nonumber\\
&= \left.\frac{G^n(\bY^n| {z}(\bp^G)_n, \ba, \bb)}{\partial {z}(\bp^G)_n} \right\rvert_{{\substack{\ba = \ba^*(\bz)\\\bb= \bb^*(\bz)}}}\frac{\partial z(\bp^G)_n}{\partial p^G_l}\label{eq:lemma:Ro-hete-dev-z-eq3}
\end{align}
where $h^n(\bp^G) =({p^G_n})/({1-\sum_{l\in[N]} p^G_l})$. We also take the first derivatives of  \eqref{eq:lemma:Ro-hete-dev-z-eq2} and obtain
\begin{equation}
 \label{eq:lemma:Ro-hete-dev-z-eq4}
  \frac{\partial h^n(\bp^G)}{\partial p^G_l}  = \frac{G^n(\bY^n| \widetilde{z}(\bp^G)_n, \ba^*(\bz), \bb^*(\bz))}{\partial \widetilde{z}(\bp^G)_n}\frac{\partial \widetilde{z}(\bp^G)_n}{\partial p^G_l}.
\end{equation}
Now we just combine \eqref{eq:lemma:Ro-hete-dev-z-eq3} and \eqref{eq:lemma:Ro-hete-dev-z-eq4} and the result  that $\bz(\bp^G) = \widetilde{\bz}(\bp^G)$ to have  the first equation of (ii). The second equation of   (ii) can be verified  similarly, as we just need to take the second-order derivatives of \eqref{eq:lemma:Ro-hete-dev-z-eq1} and \eqref{eq:lemma:Ro-hete-dev-z-eq2} and use the results from (i) and \eqref{eq:lemma:Ro-hete-dev-z-eq4} to obtain the desired equality.  
\endproof
We are now to provide a complete proof for Theorem \ref{theor;RO-hete-convexness-simplified-model}.
\proof{Proof of Theorem \ref{theor;RO-hete-convexness-simplified-model}:}
The first claim is already validated in Lemma \ref{lemma:zp-unique-deter}. To prove that $\cW(\bz(\bp^G))$ is strictly concave, we will show that its second derivative $\cW(\bz(\bp^G))$ is negative definite. This can be easily seem as
\begin{align}
\frac{\partial \cW(\bz(\bp^G))}{\partial p^G_n} &= z(\bp^G)_n +\sum_{l\in[N]} p^G_l \frac{\partial z(\bp^G)_l }{\partial p^G_n} \nonumber \\
\frac{\partial^2 \cW(\bz(\bp^G))}{\partial p^G_n\partial p^G_k} &=\frac{\partial z(\bp^G)_n}{\partial p^G_k}  +  \frac{\partial z(\bp^G)_k }{\partial p^G_n} +
\sum_{l\in[N]} p^G_l \frac{\partial^2 z(\bp^G)_l }{\partial p^G_n\partial p^G_k} \nonumber
\end{align}
Then  using Lemma \ref{lemma:Ro-hete-dev-z} we have
\begin{equation}
\label{eq:eq234}
\frac{\partial^2 \cW(\bz(\bp^G))}{\partial p^G_n\partial p^G_k} =  \left.\left(\frac{\partial \widetilde{z}(\bp^G|\ba,\bb)_n}{\partial p^G_k}  +  \frac{\partial \widetilde{z}(\bp^G|\ba,\bb)_k }{\partial p^G_n} +
\sum_{l\in[N]} p^G_l \frac{\partial^2 \widetilde{z}(\bp^G|\ba,\bb)_l }{\partial p^G_n\partial p^G_k}\right) \right\rvert_{{\substack{\ba = \ba^*(\bz)\\\bb= \bb^*(\bz)}}}
\end{equation}
It not difficult to see that the left hand side of \eqref{eq:eq234} is equal to $\partial^2 \widetilde{\cW}(\widetilde{\bz}(\bp^G|\ba^*(\bz),\bb^*(\bz)) )/ (\partial p^G_n\partial p^G_k)$, leading to
\[
\nabla^2 \cW(\bz(\bp^G)) = \nabla^2 \widetilde{\cW}(\widetilde{\bz}(\bp^G|\ba^*(\bz),\bb^*(\bz)) ).
\]
Since $\nabla^2 \widetilde{\cW}(\widetilde{\bz}(\bp^G|\ba^*(\bz),\bb^*(\bz)) )$ is always negative definite (Lemma \ref{lemma:Ro-hete-deter-simplified}), so is $\nabla^2 \cW(\bz(\bp^G))$. This completes the proof. 
\endproof

\subsection{Proof of Theorem \ref{theor:convert-simplified-model}}
\label{apdx:proof-concert-to-reduced-model}

We will make use of Lemmas \ref{lemma:RO-GEV-hete-bounds}-\ref{lemma:RO-GEV-hete-solutions} below to prove the claim.  Lemma \ref{lemma:RO-GEV-hete-bounds} shows that under any prices $\bz\in\bbR^N$, the adversary will force $G^n(\bY^n|z,\ba^n,\bb^n)$ to either its maximum or minimum value, for any $\n\in[N]$, with a note that both cases can occur. We further, in Lemmas \ref{lemma:psi-monotone} and \ref{lemma:RO-GEV-hete-solutions}, show that, under an optimal price solution, the adversary will always force $G^n(\bY^n|z_n,\ba^n,\bb^n)$ to their minimums. This is an essential claim to show the equivalence between the robust problem and the reduced one.

First, let
\[
\overline{\cG}^n(z_n) = \max_{(\ba^n,\bb^n) \in\cA^n}\; \Big\{ G^n(\bY^n| z_n,\ba^n,\bb^n)\Big\}.
\]
\begin{lemma}
\label{lemma:RO-GEV-hete-bounds}
Give a markup vector $\bz\in\bbR^N$, let  $(\ba^*,\bb^*) = \{(\ba^{n*}, \bb^{n*}),\ n\in[N]\} $ be a solution to the adversary's problem of \eqref{prob:RO-GEV-hete}, then for any $n\in [N]$, we have
\[
G^n(\bY^n|z,\ba^{n*},\bb^{n*}) = 
\begin{cases}
\underline{\cG}^n(z_n)  &\text{ if }\rho(\bz,\ba^*, \bb^*)< z_n \\
\overline{\cG}^n(z_n)  &\text{ if }\rho(\bz,\ba^*, \bb^*)> z_n \\
\end{cases}\qquad \forall n \in [N].
\]
Moreover, if there are a set of indexes $\cN\subset [N]$ such that $\rho(\bz,\ba^*, \bb^*) = z_n$, $\forall n\in\cN$, then all the solutions in the following set are optimal to the adversary's problem
\[
S^* = \{(\ba,\bb)\in \cA|\ G^n(\bY^n|z,\ba^l,\bb^{l}) = G^n(\bY^n|z,\ba^{l*},\bb^{l*}),\forall l\in[N], n\notin\cN\}.
\]
\end{lemma}
\proof{Proof:}
Given $n\in[N]$, let us denote 
\begin{align}
 A &=  \sum_{l\in[N],l\neq n} {z}_l G^l(\bY^l|{z}_l,{\ba^{l*}},{\bb^{l*}})\nonumber \\
 B &= 1+ \sum_{l\in[N],l\neq n}  G^l(\bY^l|{z}_l,{\ba^{l*}},{\bb^{l*}}).\nonumber
\end{align}
We can write 
\[
\rho(\bz,\ba^*, \bb^*) = \frac{A + {z}_n G^n(\bY^n|{z}_n,{\ba^{n*}},{\bb^{n*}}) }{B+  G^n(\bY^n|{z}_n,{\ba^{n*}},{\bb^{n*}})}
\]
We prove the lemma by considering  the following three cases:
\begin{itemize}
 \item[(i)] If $\rho(\bz,\ba^*, \bb^*)< z_n$, then for any $\gamma < G^n(\bY^n|{z}_n,{\ba^{n*}},{\bb^{n*}})$
 one can easily show the following inequality 
 \[
 \rho(\bz,\ba^*, \bb^*) =  \frac{A + {z}_n G^n(\bY^n|{z}_n,{\ba^{n*}},{\bb^{n*}}) }{B+  G^n(\bY^n|{z}_n,{\ba^{n*}},{\bb^{n*}})} > \frac{A + {z}_n \gamma }{B+  \gamma} 
 \]
 Since $(\ba^*, \bb^*)$ is a solution to the adversary problem \eqref{prob:RO-GEV-hete}, $G^n(\bY^n|{z}_n,{\ba^{n*}},{\bb^{n*}})$ must be equal to its minimum,  i.e., $G^n(\bY^n|{z}_n,{\ba^{n*}},{\bb^{n*}}) = \underline{\cG}^n(z_n)$.
  \item[(ii)]  If $\rho(\bz,\ba^*, \bb^*)> z_n$, then similarly to the previous case, we can show that, for any $\gamma > G^n(\bY^n|{z}_n,{\ba^{n*}},{\bb^{n*}})$
 one can easily show the following inequality 
 \[
 \rho(\bz,\ba^*, \bb^*) > \frac{A + {z}_n \gamma }{B+  \gamma}. 
 \]
Thus,$G^n(\bY^n|{z}_n,{\ba^{n*}},{\bb^{n*}})$ must be equal to its maximum,  i.e., $G^n(\bY^n|{z}_n,{\ba^{n*}},{\bb^{n*}}) = \overline{\cG}^n(z_n)$.
 \item[(iii)]  If $\rho(\bz,\ba^*, \bb^*) = z_n$, then for any $\gamma\in \bbR$ we have
 \[
 \rho(\bz,\ba^*, \bb^*) = \frac{A + {z}_n \gamma }{B+  \gamma},
 \]
 meaning that any solution $(\ba^n, \bb^n)\in \cA^n$ would be chosen minimize the adversary's objective function.
\end{itemize}
Combining the above three cases, we obtain the desired result. 
\endproof
Here we remark that the behavior of the adversary showed in Lemma \ref{lemma:RO-GEV-hete-bounds} depends on the values of  $\bz$ and both cases can occur (Remark \ref{remark:1}).
\begin{remark}
\label{remark:1}
\textit{
Given a constant-markup vector $\bz\in\bbR^N$, let $n_1 = \argmax_{n\in[N]} \{z_n\}$,
 the adversary would need to force  $G^{n_1}(\bY^{n_1}|z_{n_1},\ba^{{n_1}},\bb^{n_1})$ to its minimum over $\cA^{n_1}$. Moreover, if there is $n_2 \in[N]$ such that  the markup $z_{n_2} = 0$ and $z_n>0$ for all $n\neq n_2$, then the adversary would need to force  $G^{n_2}(\bY^{n_2}|z_{n_2},\ba^{n_2},\bb^{n_2})$ to its maximum over $\cA^{n_2}$. So, in general, the adversary would force $G^{n}(\bY^{n}|z_{n},\ba^{{n}},\bb^n)$ to either its minimum or its maximum value, and both cases can occur. } 
\end{remark}
The remark is easy to verify, as we see that if $n_1 = \argmax_{n\in[N]} \{z_n\}$, then 
\[
 \rho(\bz,\ba^*, \bb^*) \leq z_{n_1}  \frac{\sum_{n\in[N]} G^n(\bY^n| z_n,\ba^{n*},\bb^{n*})}{1 + \sum_{n\in[N]} G^n(\bY^n|z_n,\ba^{n*},\bb^{n*})} <z_{n_1} 
\]
Thus,  according to Lemma \ref{lemma:RO-GEV-hete-bounds}, the adversary would need to force  $G^{n_1}(\bY^{n_1}|z_{n_1},\ba^{{n_1}},\bb^{n_1})$ to its minimum over $\cA^{n_1}$. On other hand, if $z_{n_2} = 0$ and $z_n>0$ for all $n\neq n_2$, then $z_{n_2}< \rho(\bz,\ba^*, \bb^*)$, meaning that adversary would need to force  $G^{n_2}(\bY^{n_2}|z_{n_2},\ba^{n_2},\bb^{n_2})$ to its maximum over $\cA^{n_2}$.

We now further characterize the adversary problem under an optimal prices $\bz^*$ by showing that  $\rho(\bz^*,{\ba}^*,{\bb}^*) \leq z^*_n$ for all $n\in[N]$, meaning that the adversary always  force $G^{n}(\bY^{n}|z^*_{n},\ba^{{n}},\bb^{n})$ to its minimum. 
Before doing this,  let us define a function $\psi(\bz|\bu): \bbR^N\rightarrow \bbR$ parameterized by a binary vector $\bu \in\{0,1\}^N$, in such a way that $\psi(\bz|\bu)$ has the following form
\begin{equation}
\label{eq:define-psi}
\psi(\bz|\bu)  =\frac{\sum_{n\in[N]} z_n \theta^n(z_n)}{1+\sum_{n\in[N]} \theta^n(z_n)}, 
\end{equation}
where $\theta^n(z_n) = \underline{\cG}^n(z_n)$ if $u_n=0$ or $\theta^n(z_n) = \overline{\cG}^n(z_n)$ if $u_n=1$. A binary vector $\bu$ can be referred to as a configuration of the adversary's objective function and Lemma \ref{lemma:RO-GEV-hete-bounds} tells us that, for any $\bz\in \bbR^N$, there is $\bu\in\{0,1\}^N$ such that the adversary's objective function can be  written as $\min_{(\ba,\bb)\in\cA} \rho(\bz,\ba, \bb) = \psi(\bz|\bu)$.
We also denote  $\bbe^k$ as a vector of size $N$ with zero elements except the $k$-element that is equal to 1, for any $k\in [N]$.

We need Lemma \ref{lemma:psi-monotone} below to support the main claim.  
\begin{lemma}
\label{lemma:psi-monotone}
Given $\bz\in \bbR^N$ and $\bu\in\{0,1\}^N$, if there is $k\in [N]$ such that $\psi(\bz|\bu) \geq  z_k$, then given any $\epsilon>0$ we have $\psi(\bz|\bu) < \psi(\bz + \epsilon \bbe^k|\bu)$.
\end{lemma}
\proof{Proof:}
For notational brevity, let $A = \sum_{n\in[N],  n\neq k} z_n \theta^n(z_n)$ and $B = 1+\sum_{n\in[N],  n\neq k} \theta^n(z_n)$. We write 
\begin{align}
\psi(\bz|\bu) &= \frac{A+ z_k \theta^k(z_k)}{B+ \theta^k(z_k)} \nonumber \\
\psi(\bz + \epsilon \bbe^k|\bu) &= \frac{A+ (z_k+ \epsilon) \theta^k(z_k+\epsilon)}{B+ \theta^k(z_k+\epsilon)} \nonumber
\end{align}
Thus we have
\begin{align}
\psi(\bz + \epsilon \bbe^k|\bu) - \psi(\bz|\bu)&=  \frac{A\theta^k(z_k) + B(z_k+\epsilon)\theta^k(z_k+\epsilon) +\epsilon \theta^k(z_k)\theta^k(z_k+\epsilon) - A\theta^k(z_k+\epsilon) - B z_k \theta^k(z_k)}{(B+ \theta^k(z_k))(B+ \theta^k(z_k+\epsilon)) }\nonumber\\
&=  \frac{(A-Bz_k) (\theta^k(z_k)- \theta^k(z_k+\epsilon))  + B\epsilon\theta^k(z_k+\epsilon) +\epsilon \theta^k(z_k)\theta^k(z_k+\epsilon) }{(B+ \theta^k(z_k))(B+ \theta^k(z_k+\epsilon)) }\nonumber\\
&> \frac{(A-Bz_k) (\theta^k(z_k)- \theta^k(z_k+\epsilon)) }{(B+ \theta^k(z_k))(B+ \theta^k(z_k+\epsilon)) }.\label{eq:psi-mono-eq1}
\end{align}
Moreover, from the assumption $\psi(\bz|\bu) \geq z_k$, we can easily see that $A\geq Bz_k$. On the other hand,  we know that $\underline{\cG}^k(z_k)$ and $\overline{\cG}^k(z_k)$ are monotonic decreasing in $z_k$, so $\theta^k(z_k)$ is also monotonic-decreasing in $z_k$. Thus $(A-Bz_k) (\theta^k(z_k)- \theta^k(z_k+\epsilon)) \geq 0$. Combine this with \eqref{eq:psi-mono-eq1} we have $\psi(\bz + \epsilon \bbe^k|\bu) > \psi(\bz|\bu)$ as desired. 
\endproof

We are now to show that under an optimal price vector $\bz^*$, the adversary needs to force each component $G^n(\cdot)$ to its minimum. The proof idea is to show that if it is not the case, then we can always find another price solution $\bz'$ that yields a better worst-case profit.   
%%%%%%%%%%%%%%%%%%%%%%%%%%
%%%%%%%%%%%%%%%%%%%%%%%%%
\iffalse
\ctien{The following needs to be revised ... }
\begin{lemma}
\label{lemma:RO-GEV-hete-rho-monotonic}
Given $\overline{\bz} \in \bbR^N$, $({\ba}, {\bb})\in \cA$ and $k\in [N]$, if $\rho(\overline{\bz},{\ba},{\bb})\geq \overline{z}_k$, then 
\[
\left.\frac{\partial \rho(\bz,\ba,\bb)}{\partial z_k} \right\rvert_{{\bz = \overline{\bz}}} >0.
\]
\end{lemma}
\proof{Proof:}
For notational brevity, let $A = \sum_{n\in[N]} \overline{z}_n G^n(\bY^n|\overline{z}_n,{\ba^n},{\bb^n})$ and $B = 1+ \sum_{n\in[N]}  G^n(\bY^n|\overline{z}_n,{\ba^n},{\bb^n})$. Taking the first derivative of $\rho(\bz,\ba,\bb)$ w.r.t. $z_k$ we obtain
\[
\begin{aligned}
\left.\frac{\partial \rho(\bz,\ba,\bb)}{\partial z_k} \right\rvert_{{\bz = \overline{\bz}}} &= \frac{\left(G^k(\overline{z}_k,{\ba^k},{\bb^k}) - \overline{z}_k\sum_{i\in\cV_k}b_iY_i\partial G^k_i(\bY^k)\right) B + A\left(\sum_{i\in\cV_k}b_iY_i\partial G^k_i(\bY^k)\right) }{B^2}\\
& = \frac{1}{B^2}\left( G^k(\overline{z}_k,{\ba^k},{\bb^k}) +  \left(\sum_{i\in\cV_k}b_iY_i\partial G^k_i(\bY^k)\right) \left(A - \overline{z}_kB \right)\right) >0.
\end{aligned}
\]
The last inequality is due to $\rho(\overline{\bz},{\ba},{\bb}) = A/B \geq \overline{z}_k$.
\endproof
\fi

\begin{lemma}\label{lemma:RO-GEV-hete-solutions}
Under robust optimal prices $\bz^*\in\bbR^{N}$, let  $(\ba^*,\bb^*) = \{(\ba^{n*}, \bb^{n*}),\ n\in[N]\} $ be a solution to the adversary's problem of \eqref{prob:RO-GEV-hete}, then for any $n\in [N]$, we have
\[
G^n(\bY^n|z^*_n,\ba^{n*},\bb^{n*}) = 
\underline{\cG}^n(z_n).  \]
\end{lemma}
\proof{Proof:}
%Let $\bz^*$ be an  optimal solution to \eqref{prob:RO-GEV-hete}, we first prove that $(\underline{\ba},\overline{\bb})$ is optimal to the adversary's problem under $\bz^*$.
%Let 
Let $f(\bz) = \text{argmin}_{(\ba,\bb) \in \cA}\ \rho(\bz,\ba,\bb)$
and  $\cU^*$ be the set of parameter $\bu$ such that 
$f(\bz^*) = \psi(\bz|\bu)$. 
 To prove the equality, we just need to show that  $f(\bz^*)< z^*_n$ for all $n \in[N]$. By contradiction, assume that there exists $n \in [N]$ such that 
$\rho(\bx^{*},\ba^*,\bb^*) \geq z^*_n$. Let
\[
\begin{aligned}
 k &= \text{argmax}\{z^*_n|\ n\in [N], f(\bz^*) > z^*_n \}\\
 h &=\text{argmin}\{z^*_n|\ n\in [N], f(\bz^*)<  z^*_n \}.
\end{aligned}
\]
Indeed, $h$ always exists because $\rho(\bz^{*},\ba^*,\bb^*)<\max_{n\in[N]}z^*_n$.
We consider two following cases
\begin{itemize}
 \item[(i)] If such $k$ exists. We have $
z^*_h >\rho(\bz^{*},\ba^*,\bb^*) > z^*_k
$ and for any $l\neq h$ and $l\neq k$ we have either $\rho(\bz^{*},\ba^*,\bb^*) = z^*_l$ or $z^*_k\geq z^*_l $ or $z^*_h\leq z^*_l$. Moreover, the  function $f(\bz)$ is continuous in $\bx$ \citep[Theorem 7, ][]{Hogan1973}. So, there is $\delta>0$ such that 
\[
z^*_h > f(\bz^{*}+ t\bbe^k)> z^*_k + t,\ \forall t\in [0,\delta], 
\]
As a result, for any $l\in[N]$ such that $z^*_l \geq z^*_h$ we have  $z^*_l > f(\bz^{*}+ t\bbe^k)$ and if $z^*_l \leq z^*_k$ we have  $z^*_l < f(\bz^{*}+ t\bbe^k)$. From Lemma \eqref{lemma:RO-GEV-hete-bounds}, this means that there is a parameter $\Bar{\bu} \in\bU^*$ and $t\in(0,\delta)$ such that
\begin{align}
f(\bz^*+ t\bbe^k) &= \psi(\bz^*+ t\bbe^k|\Bar{\bu})\nonumber\\
f(\bz^*) &= \psi(\bz^*|\Bar{\bu})\nonumber
\end{align}
Moreover, from Lemma \ref{lemma:psi-monotone}, we see that $\psi(\bz^*+ t\bbe^k|\Bar{\bu})> \psi(\bz^*|\Bar{\bu})$, or equivalently, $f(\bz^*+ t\bbe^k)>f(\bz^*)$,
which is contradictory to the assumption that $\bz^*$ is a robust solution to \eqref{prob:RO-GEV-hete}. 
\item[(ii)] If such $k$ does not exist, then $f(\bz^*) = z^*_n$ and  for any $l\in[N]$, either $f(\bz^*)<z^*_h\leq z^*_l$ or $f(z^*) = z^*_l$. Similar to the previous case, we also have the result that there exists $\delta>0$ such that $z^*_l >f(\bz^*+t\bbe^n)$ for any $l\in[N]$ and $l\neq n$  and for all $t\in(0,\delta)$, which also leads to the result that there is $\Bar{\bu}\in\bU^*$ such that $f(\bz^*+t\bbe^n) = \psi(\bz^*+t\bbe^n|\Bar{\bu})$. Using Lemma \ref{lemma:psi-monotone} and the fact that $f(\bz^*) = z^*_n$, we have $\psi(\bz^*+t\bbe^n|\Bar{\bu}) > \psi(\bz^*|\Bar{\bu})$ for a $t\in (0,\delta)$. Thus, $f(\bz^*+t\bbe^n)>f(\bz^*)$, which is also contradictory to the assumption that $\bz^*$ is a robust solution to \eqref{prob:RO-GEV-hete}. 
\end{itemize}
So, in closing, we can claim that 
$f(\bz^*) < z^*_n$, for all $n\in[N]$. Thus,  from Lemma \ref{lemma:RO-GEV-hete-bounds} we obtain the desired result.
\endproof

Now we know that under the optimal prices, the adversary will force each function $G^n(\bY^n|z^*_n,\ba^{n},\bb^{n})$ to its minimum over $\cA^n$, suggesting that we may be able to convert the robust problem into  the maximization problem in \eqref{prob:RO-hete-simlified}. With all the lemmas above, we are ready to prove Theorem \ref{theor:convert-simplified-model}.

\proof{Proof of Theorem \ref{theor:convert-simplified-model}:}
 We need to prove that if $\bz^*$ is a robust optimal solution to \eqref{prob:RO-GEV-hete}, then it is also optimal to \eqref{prob:RO-hete-simlified}, and vice-versa.
 \iffalse
To facilitate our exposition, let define the set
\[
\Xi = \left\{\theta \Big|\; \theta = \frac{\sum_{n\in[N]} z^*_n y_n}{1 + \sum_{n\in[N]} y_n},\; \forall \by\in \bbR^N, y_n \in \{\underline{\cG}^n(z_n), \overline{\cG}^n(z_n)\}  \right\}
\]

\here{}

In other words, $\Xi$ contains all adversary's objective values where functions $G^n(\bY^n|z,\ba^{n},\bb^{n})$ are equal to their maximums or minimums. We also define
\begin{equation}\label{eq:definition-tau}
\tau = \min_{\theta\in \Xi}\left\{ \theta - w(\bz^*)\Big|\ \theta > w(\bz^*) \right\}
\end{equation}
\fi
By contradiction, assume that $\bz^*$ is optimal to \eqref{prob:RO-GEV-hete} but $\bz^* \notin \text{argmax}_{\bz \in \bbR^{N}} \cW(\bz)$. Since the problem $\max_{\bz \in \bbR^{N}} \cW(\bz)$ has a unique local optimum (Proposition \ref{prop:RO-hete-fixed-point}), $\nabla_\bz \cW(\bz) \neq 0$.  Thus, there always exits a vector $\pmb{\epsilon}\in\bbR^N\neq 0$  and a constant $\delta >0$ such that
\begin{equation}\label{eq:ineq-optimal}
\cW(\bz^*) < \cW(\bz^*+ t\pmb{\epsilon}) ,\forall t\in(0,\delta). 
\end{equation}
Moreover, according to Lemma \ref{lemma:RO-GEV-hete-solutions}, we know that $\cW(\bz^*) <z_n^*$ for all $n\in [N]$. Since $\cW(\bz)$ and $f(\bz)$ are continuous in $\bz$ (recall that $f(\bz)  = \min_{\ba,\bb} \rho(\bz,\ba,\bb)$) and $f(\bz^*) = \cW(\bz^*)$, we  can always choose $\delta_1 \in (0,\delta)$ such that
\begin{equation}
f(\bz^* + t\pmb{\epsilon}) < z^*_n < z^*_n + t\epsilon_1,\;\forall n\in[N], t\in(0,\delta_1).
\end{equation}
Thus, using Lemma \ref{lemma:RO-GEV-hete-bounds}, we see that the adversary under prices $\bz^* + t_1\pmb{\epsilon}$ will also force $G^n(\bY^n|\bz^* + t_1\pmb{\epsilon},\ba,\bb)$ to be equal to their minimum values. Hence, we have  
\[
f(\bz^* + t\pmb{\epsilon}) = \cW(\bz^* + t\pmb{\epsilon}) \stackrel{(i)}{>} \cW(\bz^*) = f(z^*),
\]
where (i) is due to \eqref{eq:ineq-optimal}. This is contradictory against the assumption that $\bz^*$ is optimal to \eqref{prob:RO-GEV-hete}. So, $\bz^*$ needs to be optimal to \eqref{prob:RO-hete-simlified} as well. For the opposite side, Proposition \ref{prop:RO-hete-fixed-point} already tells us that \eqref{prob:RO-hete-simlified} always yields a unique optimal solution $\bz^*$. Thus, this solution is also optimal to \eqref{prob:RO-GEV-hete}. 
\endproof

%%%%%%%%%%%%%%%%%%%%
%%%%%%%%%%%%%%%%%%%%

\subsection{Proof of Proposition \ref{prop:RO-hete-fixed-point}}
\label{proof:RO-hete-fixed-point}
We first prove the equality \eqref{eq:temp3452}. Taking the first-order derivatives of the objective function  $\cW(\bz)$, we have
\begin{align}
\frac{\partial \cW(\bz)}{\partial z_n} &=  \frac{(1+\sum_{l\in[N]}\underline{\cG}^l(z_l) ) ( \underline{\cG}^n(z_n) + z_n \partial \underline{\cG}^n(z_n)/\partial z_n  )  + (\sum_{l\in[N]} z_l\underline{\cG}^l(z_l)) (\partial \underline{\cG}^n(z_n)/\partial z_n) }{(1+\sum_{l\in[N]}\underline{\cG}^l(z_l) )^2} 
\nonumber \\  
&=\left.\frac{(1+\sum_{l\in[N]}{G}^l(\bY^l|z_l,\ba,\bb) ) ( {G}^n(\bY^n|z_n,\ba,\bb) + z_n \partial {G}^n(\bY^n|z_n,\ba,\bb)/\partial z_n  ) }{(1+\sum_{l\in[N]}{G}^l(\bY^l|z_l,\ba,\bb) )^2}  \right\rvert_{{\substack{\ba = \ba^*(\bz)\\\bb= \bb^*(\bz)}}}\nonumber \\ 
&\qquad +\left.\frac{(\sum_{l\in[N]}z_l{G}^l(\bY^l|z_l,\ba,\bb) ) ( \partial {G}^n(\bY^n|z_n,\ba,\bb)/\partial z_n  ) }{(1+\sum_{l\in[N]}{G}^l(\bY^l|z_l,\ba,\bb) )^2}  \right\rvert_{{\substack{\ba = \ba^*(\bz)\\\bb= \bb^*(\bz)}}}\nonumber \\
&= \left. \frac{\partial \widetilde{\cW}(\widetilde{\bz}(\bp^G|\ba,\bb))}{\partial z_n} \right\rvert_{{\substack{\ba = \ba^*(\bz)\\\bb= \bb^*(\bz)}}}\nonumber
\end{align}
Now, since \eqref{prob:RO-hete-simlified} is an unconstrained problem, if $\bz^*$ is an optimal solution to \eqref{prob:RO-hete-simlified}, we have
\[
\frac{\partial \cW(\bz^*)}{\partial z_n} = 0, \;\forall n\in[N], 
\]
Moreover, From Lemma \eqref{lemma:Ro-hete-dev-z} - (i), we see that $\widetilde{\bz}(\bp^G|\ba^*(\bz^*),\bb^*(\bz^*)) = \bz(\bp^G)$, so $\bz^*$ is also a solution to the system
\begin{equation}
\label{eq:temp1234} 
 \left. \frac{\partial \widetilde{\cW}(\widetilde{\bz}(\bp^G|\ba,\bb))}{\partial z_n} \right\rvert_{{\substack{\ba = \ba^*(\bz^*)\\\bb= \bb^*(\bz^*)}}} = 0,\;\forall n\in[N].\nonumber
\end{equation}
Note that $\widetilde{\cW}(\widetilde{\bz}(\bp^G|\ba^*(\bz^*),\bb^*(\bz^*)))$ is the objective function of the deterministic pricing problem with partition-wise homogeneous PSP considered in \cite{zhang2018multiproduct}. Thus, using Theorem C1 of \cite{zhang2018multiproduct} we see that $\bz^*$ has to satisfy the system of equalities
\[
\begin{cases}
R(z^*)  =\sum_{n\in[N]} \frac{1}{\me{\bb^{n*}(\bz^*)}} {G}^n(\bY^n|z_n,\ba^*(\bz^*),\bb^*(\bz^*)) \\
z_n = \frac{1}{\me{\bb^{n*}(\bz^*)}} + R(\bz^*),\; \forall n\in[N].
\end{cases}
\]
Thus, we have
\[
\begin{aligned}
z_n^* &= \frac{1}{\me{\bb^{n*}(\bz^*)}} + \sum_{l\in[N]} \frac{1}{\me{\bb^{l*}(\bz^*)}} {G}^l(\bY^l|z^*_l,\ba^*(\bz^*),\bb^*(\bz^*)) \\
&= \frac{1}{\me{\bb^{n*}(\bz^*)}} + \sum_{n\in[N]} \frac{1}{\me{\bb^{l*}(\bz^*)}} \underline{\cG}^l(z^*_l), 
\end{aligned}
\]
which is also the desired equality \eqref{eq:temp3452}. 

We now prove that \eqref{eq:temp3452} always yields a unique local optimum. By contradiction, assume that there are to points $\bz_1,\bz_2 \in \bbR^N$ such that $\bz_1\neq \bz_2$ and $\frac{\partial \cW(\bz_1)}{\partial z_n} =  \frac{\partial \cW(\bz_2)}{\partial z_n} = 0$ for all $n\in[N]$. Let $\bp^G_1$  and $\bp^G_2$ be two vectors of purchase probabilities defined by \eqref{eq:pG-z} under $\bz_1,\bz_2$, respectively. From Lemma  \ref{lemma:zp-unique-deter} we have $\bp^G_1 \neq \bp^G_2$. Moreover, taking the derivatives of $\cW(\bz(\bp^G))$ with respect to $p^G_n$, $n\in[N]$ we get  
\[
\frac{\partial \cW(\bz(\bp^G_1))}{\partial p^G_n} = \sum_{l\in[N]}\left. \frac{\partial \cW(\bz)}{\partial z_l}\right\rvert_{\bz=\bz_1}  \times \left. \frac{\partial z(\bp^G)_l}{\partial p^G_n}\right\rvert_{\bp^G=\bp^G_1} = 0. 
\]
Similarly, we also have ${\partial \cW(\bz(\bp^G_2))}/{\partial p^G_n} = 0$, implying that both $\bp^G_1$ and $\bp^G_2$ are local optimal solutions to the problem $\max_{\bp^G}\cW(\bz(\bp^G))$, which is contrary to the claim that $\cW(\bz(\bp^G))$ is \textit{strictly concave} in $\bp^G$ (Theorem \ref{theor;RO-hete-convexness-simplified-model}). This completes the proof.  

%%%%%%%%%%%%%%%%%%%%%%%%%%%%%
%%%%%%%%%%%%%%%%%%%%%%%%%%%%%%
\subsection{Proof of Proposition \ref{prop:gradient-cW}}
\label{apdx:proof-gradient-cW}
Taking the gradient of $\cW(\bz(\bp^G))$ with respect to $p^G_n$, $n\in[N]$, we get 
\begin{equation}
\label{eq:proof6-eq1}
\frac{\partial \cW(\bz(\bp^G))}{\partial p^G_n}  = z(\bp^G)_n + \sum_{k\in[N]} \frac{\partial z(\bp^G)_k}{\partial p^G_n} p^G_k.
\end{equation}
Now, from Lemma \ref{lemma:zp-unique-deter}, for any $k\in [N]$, we have $\underline{\cG}^k(z(\bp^G)_k) = p^G_k/(1-\bbe^\T \bp^G)$. Taking the first-derivatives of the both sides with respect to $p^G_n$ we have
\begin{align}
\frac{p^G_k}{(1-\bbe^\T \bp^G)^2} + \frac{\bbI(k=n)}{(1-\bbe^\T \bp^G)} &= \frac{\partial \underline{\cG}^k(z(\bp^G)_k) }{\partial z_k}  \frac{ z(\bp^G)_k}{\partial p^G_n} \label{eq:proof6-eq2}. 
\end{align}
Moreover, we can write 
\begin{align}
\frac{\partial \underline{\cG}^k(z_k) }{\partial z_k} &= \left.\frac{\partial  G^k(\bY^k|z_k,\ba^k,\bb^k)}{\partial z_k} \right\rvert_{\substack{\ba^{k} = \ba^{k*}(z_k) \\ \bb^{k} = \bb^{k*}(z_k)}}
\nonumber \\
&=  \left.  \frac{ - \sum_{i\in\cV_k}  \partial G^k_i(\bY^k|z_k,\ba^k,\bb^k) Y^i \me{\bb^k}}{\partial z_k} \right\rvert_{\substack{\ba^{k} = \ba^{k*}(z_k) \\ \bb^{k} = \bb^{k*}(z_k)}}\nonumber \\
&=-\me{\bb^{k*}(z_k)} \underline{\cG}^k(z_k) \nonumber \\
&= \frac{-\me{\bb^{k*}(z_k)} p^G_k}{1-\bbe^\T\bp^G }.  \label{eq:proof6-eq3}
\end{align}
Combine \eqref{eq:proof6-eq2} and \eqref{eq:proof6-eq3} we have
\begin{equation}
\label{eq:proof6-eq4}
\frac{z(\bp^G)_k}{\partial p^G_n} = -\frac{\bbI[k=n]}{\me{\bb^{k*}(z_k)} p^G_k} - \frac{1}{\me{\bb^{k*}(z_k)} (1-\bbe^\T \bp^G)}.
\end{equation}
We substitute \eqref{eq:proof6-eq4} into \eqref{eq:proof6-eq1} and get
\[
\frac{\partial \cW(\bz(\bp^G))}{\partial p^G_n}  =  z(\bp^G)_n - \frac{1}{\me{\bb^{k*}(z_n)}} -\frac{1}{ (1-\bbe^\T \bp^G)} \sum_{k\in [N]} \frac{p^G_k}{\me{\bb^{k*}(z_k)}},
\]
as desired.

\section{\mtien{Robust Pricing with Over-expected-sale Penalties}}
\label{sec:robust-penalties}
Motivated by applications in inventory considerations \citep{Gallego2014dynamicPricing}, we study a robust model for the  pricing problem with expected sale requirements under uncertain choice parameters $(\ba,\bb)$. We first show in Section \ref{sec:robust-con} below
that there may be no fixed prices such that the corresponding  expected sale constraints are always satisfied when the choice parameters vary in the uncertainty set. It motivates us to consider a new  robust model with over-expected-sales penalties in Section \ref{sec:over-expected-sale Penalties}. We provide the proofs of the results in this section in  Section
 \ref{sec:proof-over-expected-sale Penalties}.

\subsection{Robust Pricing with Expected Sale Constraints}
\label{sec:robust-con}

Motivated by the fact that the expected profit is concave in the purchase probabilities, previous studies \citep{zhang2018multiproduct,Keller2013} show that it is convenient to consider the pricing problem with expected sale constraints.
Technically speaking, 
given a GEV-CPGF $G(\bY)$, price vector $\bx\in \bbR^m$ andparameters $(\ba, \bb)\in\bbR^{2m}$, let us define  the vector of purchase probabilities of products $\bp$ with entries $p_i = P_i(\bx,\ba,\bb|G)$.  We also let $\bx(\bp|\ba,\bb,G)$ be the  denote the prices that achieve the purchase probabilities $\bp$. The deterministic version of the constrained pricing problem can be formulated  as 
\begin{equation}\label{prob:det-nopenalties-inner-convexOpt}
\underset{\bp\in \cP}{\text{max}}\qquad \left(\sum_{i\in\cV}\bx(\bp|\ba,\bb)_i-c_i\right)p_i . 
\end{equation}
where $\cP\in\bbR^m$ is a convex set such that for all $\bp\in\cP$, $\sum_{i\in\cV}p_i \leq 1$. One the optimal purchase probabilities $\bp$ is specified, we can obtain the \textit{optimal prices}  $\bx(\bp|\ba,\bb,G)$ by solving a convex optimization problem.
A natural robust version of the constrained pricing problem can be formulated as 
\begin{equation}\label{eq:robust-constrained}
 \max_{\bp \in \cP} \left\{\phi(\bp) =  \min_{(\ba,\bb)\in \cA}  \sum_{i\in\cV}\left( \bx(\bp|\ba,\bb,G)_i-c_i\right)p_i  \right\}.
\end{equation}
Even though it is not difficult to show \eqref{eq:robust-constrained} is computationally tractable under rectangular or some polyhedrons uncertainty sets, the issue here is that the final decision is a price vector, not  purchase probabilities. So even if we get an optimal purchase probabilities $\bp$ from the robust model, it is not clear how to compute the corresponding optimal  prices under $(\ba,\bb)$ uncertainty. On the other hand, one can show that given any prices $\bx$, there may be $(\ba,\bb)\in \cA$ such the resulting purchase probability vector $\bp = P(\bx,\ba,\bb|G)$ that does not belong to the feasible set (i.e. the expected sale constraints are not satisfied).  
All these make the robust version in \eqref{eq:robust-constrained} inappropriate to use. This is the reasonwe  propose an alternative robust model in Section \ref{sec:robust-penalties}, in which instead of requiring the purchasing probabilities to satisfy some constraints, we add a penalty cost to the objective function.  

Alternatively, in some situations the firm may face uncertainties occurring in the inventory, leading to uncertain expected sale constraints. A robust model may require  the expected sales constraints to be satisfied for all the scenarios that may occur, i.e., $\bp \in  \cP(\xi)$, for all $\xi\in\Xi$. Such a robust model can be formulated as
\begin{align}
	\underset{\bp}{\text{max}}\qquad &  \left(\sum_{i\in\cV}\bx(\bp|\ba,\bb)_i-c_i\right)p_i & \label{prob:rb-ES-constraints} \\
	 \text{subject to} \qquad & (\bap^t(\xi))^\transpose \bp \leq r_t(\xi) & \forall \xi\in\Xi \nonumber\\
	  &  \sum_{i\in \cV}p_i\leq 1,\ \bp\geq 0 &  \nonumber
\end{align}
where $(\bap^t(\xi),r_t(\xi))$, $\forall t$, are the parameters of the expected sale constraints, which are not certain in the context and depend on a random vector $\xi\in\Xi$.
Since the objective function is concave and all the constraints are linear in $\bp$, the above problem is generally tractable \citep{Ben1998robustconvex}. A simple but useful setting is that the parameter of the expected sale constraints vary in a rectangular uncertainty set, i.e., $\underline{\bap}^t \preceq \bap^t \preceq \underline{\bap}^t$ and $\underline{r}^t \leq r^t \leq \overline{r}^t$ for all $t\in [T]$. In this context, one can show that \eqref{prob:rb-ES-constraints} is equivalent to the following convex optimization problem
\begin{align}
	\underset{\bp}{\text{max}}\qquad &  \left(\sum_{i\in\cV}\bx(\bp|\ba,\bb)_i-c_i\right)p_i & \label{prob:rb-ES-constraints-equiv} \\
	 \text{subject to} \qquad & (\overline{\bap}^t)^\transpose \bp \leq \underline{r}_t & \nonumber\\
	  &  \sum_{i\in \cV}p_i\leq 1,\ \bp\geq 0 &  \nonumber
\end{align}
Other uncertainty sets may be considered, i.e., polyhedron or ellipsoidal ones, and we refer the reader to \cite{Ben1998robustconvex} for details.

%\pj{we should indicate here in words what we mean by not being appropriate}\tm{I clarified the issue}

\subsection{Robust Pricing with Over-expected-sale Penalties}
\label{sec:over-expected-sale Penalties}
We propose a version with over-expected-sale penalties, which allows us to handle both the expected sale requirements and the uncertainty issue.  
Our idea  is to put the expected sale constraints to the objective function, i.e., we do not force the purchase probabilities to be in a feasible set, but instead we add penalties for purchase probabilities violating the constraints.  
More precisely, we consider the objective function 
$\Phi(\bx,\ba,\bb) - \sum_{t=1}^T \lambda_t \max\{0, (\bap^t)^\transpose \bp - r_t\}$, 
where  $\lambda_t\geq 0$, $t=1,\ldots,T$, are penalty parameters \mtien{and  $\bp \in\bbR^m$ is a vector of purchase probabilities with entries $p_i =  P_i(\bY(\bx,\ba,\bb)$, $\forall i\in\cV$.}
In this objective function, if a constraint is violated, i.e.,  $(\bap^t)^\transpose \bp > r_t$, then a cost $- \lambda_t \max\{0, (\bap^t)^\transpose \bp - r_t\}$ is added to the expected revenue. In general, if we choose $\lambda_t$ large enough, we will need a vector of purchase probabilities satisfying all the expected sale constraints to obtain  high objective values. The deterministic pricing problem under the above objective function is
\begin{equation}\label{prob:DET-PENalties}
\max_{\bx \in \bbR^m}\left\{\Phi(\bx,\ba,\bb) - \sum_{t=1}^T \lambda_t \max\{0, (\bap^t)^\transpose \bp - r_t\}\right\}.
\end{equation}
In general a solution to \eqref{prob:DET-PENalties} does not have a constant markup over products, even when  the PSP are all homogeneous
\begin{proposition} The problem with penalties 
\label{prop:no-constant-markup}
\eqref{prob:DET-PENalties} does not have a constant markup over products, even when the PSP are homogeneous.
\end{proposition}
 For this reason, the results presented in this section are not a generalized version of  those shown in Sections \ref{sec:robust-homo} and \ref{sec:uncon-hete} when the penalty parameters $\bld$ equals zero. Since we consider the pricing problem with over-expected-sale penalties, we do not face the issue of violating the expected sale constraints when the choice parameters vary in the uncertainty set.We consider the robust version of \eqref{prob:DET-PENalties} under $(\ba,\bb)$ uncertainty 
\begin{equation}\label{prob:robust-constrained-penalties}
 \max_{\bx\in\bbR^m} \min_{(\ba,\bb) \in \cA}\left\{ \cK(\bx,\ba,\bb)  = \Phi(\bx,\ba,\bb) - \sum_{t=1}^T \lambda_t \max\{0, (\bap^t)^\transpose \bp - r_t\} \right\}
\end{equation}
In this version, we consider the settings that the PSP are partition-wise homogeneous, and the CPGF and uncertainty set are partition-wise separable, as in Section \ref{sec:uncon-hete}.  
The adversarial problem of \eqref{prob:robust-constrained-penalties} is way more  difficult to solve as compared to the robust versions considered in the previous sections, as the objective function now is not differentiable. Moreover, a solution to the deterministic problem \eqref{prob:DET-PENalties} would not have a constant-markup style, thus  a robust solution to \eqref{prob:robust-constrained-penalties} would not have either. The adversary's problem under such a non-constant-markup solution seems not possible to handle tractably. 
For that reason, we consider a robust version in which we only seek constant-markup solutions. Moreover, to have a tractable structure, we also need to further assume that the expected-sale parameters $\bap^t$, $t\in [T]$ are partition-wise homogeneous, i.e., for any partition $n\in[N]$, $\alpha^t_i = \alpha^t_j$, $\forall i,j\in\cV_n, t\in[T]$. Let denote by $\bd^t$ a vector in $\bbR^N$ such that $\bd^t_n = \alpha^t_j$, for any $j\in\cV_n, n\in[N],t\in[T]$.
The  robust  problem with all the above settings becomes 
\begin{equation}
\label{prob:R-penalty-z}
\max_{\bz \in \bbR^N} \min_{(\ba,\bb) \in\cA} \left\{\cL(\bz,\ba,\bb) =  \frac{\sum_{n\in[N]} z_n G^n(\bY^n|z_n,\ba^n,\bb^n)}{1 + \sum_{n\in[N]}  G^n(\bY^n|z_n,\ba^n,\bb^n} - \sum_{t=1}^T \lambda_t \max\{0, (\bd^t)^\transpose \widetilde{\bp}^G(\bz,\ba,\bb) - r_t\}\right\},
\end{equation}
where 
$\widetilde{\bp}^G(\bz,\ba,\bb)\in\bbR^N$ with entries
\[
\widetilde{\bp}^G(\bz,\ba,\bb)_n =  \frac{ G^n(\bY^n|z_n,\ba^n,\bb^n)}{1 + \sum_{l\in[N]}  G^l(\bY^l|z_l,\ba^l,\bb^l}.
\]
Theorem \ref{theor:ro-pen-gev} below states that  \eqref{prob:R-penalty-z} can be solved by convex optimization. 
\begin{theorem}{\bf(Robust solutions for the robust pricing problem  with over-expected-sale penalties).}
\label{theor:ro-pen-gev}
 If $\bz^*$ is optimal to the problem
\begin{equation}\label{prob:RO-pen-GEV} 
\max_{\bz \in \bbR^N}\left\{\cH(\bz) =  \frac{\sum_{n\in[N]} z_n \underline{\cG}^n(z_n)}{1 + \sum_{n\in[N]} \underline{\cG}^n(z_n)} - \sum_{t=1}^T \lambda_t \max\{0, (\bd^t)^\transpose \widetilde{\bp}^G - r_t\}\right\},
\end{equation}
where  $\bp^G(\bz)$ is of size $N$ with entries $\widetilde{p}^G(\bz)_n = \underline{\cG}^n(z_n)\Big/\left(1+\sum_{j\in[N]} \underline{\cG}^l( z_l)\right)$, 
then the prices $\bx^*\in\bbR^m$ such that $x^*_i = c_i+z^*_n$, $\forall n\in[n],i\in\cV_n$ is optimal to the robust problem  $\max_{\bx\in X}\min_{(\ba,\bb)\in\cA} \cK(\bx,\ba,\bb)$.
 Moreover, the objective function of \eqref{prob:RO-pen-GEV} is concave in $\widetilde{\bp}^G$. 
\end{theorem}
We provide the proof in  Appendix \ref{sec:proof-theor-ro-pen-gev}. In general, the robust problem \eqref{prob:R-penalty-z}  is challenging to handle because of the  term $\sum_{t=1}^T \lambda_t \max\{0, (\bd^t)^\transpose \bp^G - r_t\}$, which makes the objective function no-longer differentiable in $\bz$. However, if we look at  the  subset $\cT\in[T]$  such that the constraints are violated only in $\cT$, we can write the objective function as 
\begin{align}
 \cH(\bz) &=\bz^\T \bp^G -\sum_{t\in\cT} \lambda_t(\bd^t)^\transpose \bp + \sum_{t\in\cT}\lambda_t r_t\nonumber \\
 &=\sum_{n\in [N]}(z_n-\sum_{t\in\cT} \lambda_t\bd^t_n) p^G_n+\sum_{t\in\cT}\lambda_t r_t\label{eq:proof4-eq1}
\end{align}
and  note that $\sum_{n\in [N]}(z_n-\sum_{t\in\cT} \lambda_t\bd^t_n) p^G_n$  is also an expected revenue with shifted item costs $c'_i = \sum_{t\in\cT} \lambda_t\bd^t_n +c_i$, $n\in [N], i\in \cV_n$.
We leverage this observation and follow the spirit of the proof of Theorem  \ref{theor:main-hete} to prove the results. The main idea is to show that, under an optimal prices of \eqref{prob:R-penalty-z}, the adversary will also force the each component $G^n(\cdot)$ to its minimums. From this, we can show an equivalence between the reduced problem \eqref{prob:RO-pen-GEV} and the robust one \eqref{prob:R-penalty-z}. The proof is indeed more complicated as we have to handle the term 
$\max\{0, (\bd^t)^\transpose \bp^G - r_t\}$.

The limitation of Theorem \ref{theor:ro-pen-gev} is that it only returns best solutions among those that have a constant markup in each partition, and all the expected sale parameters in each partition need to be homogeneous. Relaxing these assumption would make the robust problem challenging to handle  (see the discussion before \eqref{prob:R-penalty-z}). Moreover, we believe that the theorem is still useful in  contexts where the firm only wants to make pricing decisions for each group of products and only impose expected sale requirements for the whole groups instead of each single product in the groups.

An interesting and important question here is how the robust optimal value and optimal solutions change when  the penalty parameters $\bld$ increase.
To answer this, let us consider the following constrained problem
\begin{align}
	\underset{\bp^G}{\text{max}}\qquad &\cW(\bz(\bp^G))  & \label{prob:constrained-pG} \\
	 \text{subject to} \qquad & (\bd^t)^\transpose \bp^G  \leq r_t &  \nonumber\\
	  &  \sum_{n\in [N]}p^G_n \leq 1&  \nonumber\\
	  &  \bp^G\geq 0. &  \nonumber
\end{align}
We also define 
$\varphi^{\RO,\bld}$ as the optimal value of the robust problem in \eqref{prob:robust-constrained-penalties} under penalty parameters $\lambda$,  $\overline{\varphi}$ as the optimal value of the constrained problem  \eqref{prob:constrained-pG} , $\bx^{\RO,\bld}$ is an robust solution to \eqref{prob:robust-constrained-penalties} and $\bp^{G,\bld}$ is the purchase probabilities given by the robust solution $\bx^\bld$ in the worst-case.
Proposition \ref{thr:convergence-ro-pen} below tells us how the optimal value of the robust problem with over-expected-sale penalties  \eqref{prob:ro-over-ES-penalties-prob1} when the parameters $\bld$ increase.  
\begin{theorem}{\bf(Convergence of the robust optimal value when the penalty parameters $\bld$ increase).} 
\label{thr:convergence-ro-pen}
Given any $\epsilon>0$, we have 
\begin{itemize}
\item[(i)] If we select $\bld$ such that
$\min_t \lambda_t \geq (\Delta^* - \overline{\varphi})/\epsilon$ then $\sum_{t}\max\{0,(\bd^t)^\transpose \bp^{G,\bld} - r_t\}\leq \epsilon$, where $\Delta^* = \max_{\bz\in \bbR^N}\cW(\bz)$.
\item[(ii)]Assume that there are positive constant $L_i,l_i$, $i\in \cV$ such that $Y_i\partial G_i(\bY)$ is bounded from above  by $L_iY_i^{li}$ for all prices $\bx \geq 0$, then
if we select $\epsilon$ such that
\[
  \epsilon \leq {\min_{t\in[T],n\in[N]} \{d^t_n|\ d^t_i> 0\}}\min_t\left\{\frac{r_t}{(\bd^t)^\transpose \textbf{1}}\right\},
\]
then $\varphi^{\RO,\bld} - \overline{\varphi}$ can be bounded  as
\[
0\leq \varphi^{\RO,\bld} - \overline{\varphi} \leq \max\left\{ \max_{\substack{(\ba,\bb) \in\cA \\ n\in[N] \\  i\in\cV_n}} \left\{\frac{a_i - b_ic_i}{b_i} -  \frac{1}{b_il_i |\cV_n|}\log \frac{\delta(\epsilon)}{L_i} \right\},0\right\}\frac{N\epsilon}{\min_{t,n\in[N]} \{d^t_n|\ d^t_n> 0\}}
\]
where $\delta(\epsilon) = \min_t\left\{\frac{r_t}{(\bd^t)^\transpose \textbf{1}}\right\} - \frac{\epsilon}{\min_{t\in[T],n\in[N]} \{d^t_i|\ d^t_i> 0\}}$. This upper bound converges to zero linearly when $\epsilon$ tends to zero (i.e.,  $\min_t \{ \lambda_t\}$ goes to infinity).
\end{itemize}
\end{theorem}
The proof can be found in Appendix \ref{proof:convergence-pen}. Here, it is not difficult to validate that the assumption in Theorem \ref{thr:convergence-ro-pen}--(ii) holds for all the well-known GEV models in the literatures. For examples, for the MNL, $Y_i\partial G_i(\bY) = Y_i$. For a nested logit mode specified by $
G(\bY) = \sum_{n\in \cN} \left(\sum_{i\in C_n} \sigma_{in}Y_i^{\mu_n} \right)^{1/\mu_n}
$, where $\cN$ is the set of nests, $C_n$ is the corresponding nest and $\mu, \mu_n$ are some parameters, we have $Y_i\partial G_i(\bY) = Y_i^{\mu_n} \left(\sum_{j\in C_n} \sigma_{jn}Y_j^{\mu_n} \right)^{1/\mu_n-1}$. If $\mu_n>1$ then  $Y_i\partial G_i(\bY)\leq \sigma_{in}^{1/\mu_n-1} Y_i$ and if $\mu_n<1$ then $Y_i\partial G_i(Y)\leq L_nY_i^{\mu_n}$, where $L_n$ is an upper bound of $\left(\sum_{j\in C_n}\sigma_{jn} Y_j^{\mu_n} \right)^{1/\mu_n-1}$ for all $\bx\in\bbR^m_i$, which always exists.
For a more general GEV model, we note that $\partial G_{ij}(\bY)\leq 0$ and $Y_j\geq 0$ for all $i,j\in\cV, i\neq j$. As a result, we have $\partial G_{i}(\bY) \leq \partial G_{i}(\widetilde{\bY}^i)$, where $\widetilde{\bY}^i$ is a vector of size $m$ with entries $\widetilde{Y}^i_i = Y_i$ and $\widetilde{Y}^i_j = 0$ for all $j\neq i$.   Thus, $\partial G_{i}(\widetilde{\bY}^i)$ is a function of only $Y_i$. For a more complicated GEV model such as the network GEV model \citep{Daly2006general,Mai2017dynamic}, we can easily upper-bound $Y_i\partial G_{i}(\widetilde{\bY}^i)$ by a function of form $L_iY_i^{l_i}$, where $L_i,l_i>0$.

%In general, to bound the gap between two optimal values $\varphi^{\RO,\bld}$ and $\overline{\varphi}$,  we make use of  the concavity of $\cW(\bz(\bp^G))$ and the $\epsilon$-feasible model $\max_{\substack{\bp^{G}\in\cP^G}}\left\{\cW(\bz(\bp^G))\Big| \ (\bd^t)^\transpose \bp^G \leq r_t+ \epsilon,\ \forall t \right\}$ to make a connection between the constrained problem and the one with over-expected-sale penalties.

In Theorem \ref{thr:convergence-ro-pen},  (ii) tells us explicitly that the penalty term will converge to zero when the parameters $\bld$ are large enough. It also provides an estimate for $\min_t\{\lambda_t\}$ to get arbitrarily small  penalty costs. The second bound (ii) provides an upper-bound for the gap between the optimal expected revenues given by the constrained pricing problem and the  pricing problem with over-expected-sale penalties, and this upper bound converges to zero linearly when $\epsilon$ goes to zero. So in general, Problem \ref{thr:convergence-ro-pen} can be viewed as a generalized version of the constrained pricing problem, in the sense that if we select the penalty parameters $\bld$ large enough, then we will get a solution that is similar to the one from the constrained problem, and if we set $\bld = 0$ then we come back to the unconstrained problem. Thus, the formulation in \eqref{prob:R-penalty-z} provides a more flexible way to handle expected sale requirements.

\subsection{Proofs of the Results in Section \ref{sec:over-expected-sale Penalties} }
\label{sec:proof-over-expected-sale Penalties}

\subsubsection{Proof of Proposition \ref{prop:no-constant-markup}. }
We will give a counter example to illustrate the claim. For the sake of illustration, we only consider a pricing problem under the MNL model with 2 products and  homogeneous PSP. We also consider only one expected sale constraint as $\alpha_1p_1\leq r_t$, where $r_t/\alpha_1$ is very small.
Let $\bp^\lambda = (p^\lambda_1,p^\lambda_2)$ be a solution  to \eqref{prob:det-nopenalties-inner-convexOpt} under penalty parameter $\lambda$.  When $\lambda$  goes to infinity,  Theorem \ref{theor:det-over-ES-penalties} tells us that a solution to the pricing problem with penalties converge to a solution to the constrained pricing problem. Thus, for any $\epsilon>0$  arbitrarily small, we can chose $\lambda$ large enough such that $\alpha_1p^\lambda_1\leq r_t+\epsilon$. So, if we choose $\epsilon$  and $r_t/\alpha_1$ to be very small, then $p^\lambda_1$ would be very close to zero. Since $p_1^\lambda= \exp(a_1-bx^\lambda_1)\Big/\left(1+ \exp(a_1-bx^\lambda_1)+ \exp(a_2-bx^\lambda_2)\right)$ (where $b$ is the PSP of the two products, $\bx^\lambda$ is an optimal price solution to the pricing problem under penalty parameter $\lambda$), we have $\lim_{p_1^\lambda\rightarrow 0} x^\lambda_1(\bp) = +\infty$, meaning that to have an arbitrarily small probability $p_1^\lambda$, we need to increase the price of Product 1 to infinity. On the other hand, $p^\lambda_2$ does not affect the penalty term and we have $\lim_{x_2\rightarrow+\infty}(x_2-c_2)\exp(a_2-bx_2)\Big/\left(1+ \exp(a_1-bx^\lambda_1)+ \exp(a_2-bx_2)\right) = 0$. Thus,  to maximize the objective function, the solution $x^\lambda_2$ needs to be finite. So, in summary, we can create an example yielding a solution $(x^\lambda_1, x^\lambda_2)$ such that $x^\lambda_1$ can be arbitrarily large and $x^\lambda_2$ is bounded from above. Thus, $(x^\lambda_1, x^\lambda_2)$ would  not have the constant-markup style.

%%%%%%%%%%%%%%%%%%%%%%%%%%%%
%%%%%%%%%%%%%%%%%%%%%%%%%%%%
\subsubsection{Proof of Theorem \ref{theor:ro-pen-gev}.}
\label{sec:proof-theor-ro-pen-gev}

First, let $f(\bz)  = \min_{(\ba,\bb)\in\cA)} \cL(\bz,\ba,\bb).$
In Lemma \ref{lemma:rb-pen-lm1} below, we show that  given any prices $\bx\in\bbR^m$, The adversary's problem will force each component $G^n(\bY^n|z_n,\ba^{n},\bb^{n})$ to either its minimums or maximums. This result is similar to the case of partition-wise PSP without penalties considered in Section \ref{sec:uncon-hete}. The proof is however  more complicated as it involves the term  $\sum_{t=1}^T \lambda_t \max\{0, (\bd^t)^\transpose \bp^G - r_t\}$.
\begin{lemma}
\label{lemma:rb-pen-lm1}
Given any $\bz\in\bbR^N$, there is a solution  $(\ba^*,\bb^*) = \{(\ba^{n*},\bb^{n*})|\; n\in [N]\}$  to the corresponding adversary's problem \eqref{prob:R-penalty-z} such that
\[
G^n(\bY^n|z_n,\ba^{n*},\bb^{n*}) \in \Big\{\underline{\cG}^n(z_n), \overline{\cG}^n(z_n)\Big\}.
\]
\end{lemma}
\proof{Proof:}
 We  denote by $\cT$ a subset of $[T]$ such that $(\bd^t)^\transpose  \widetilde{\bp}^G(\bz,\ba^*, \bb^*) \geq r_t$ for all $t\in\cT$ and $(\bd^t)^\transpose \widetilde{\bp}^G(\bz,\ba^{*}, \bb^{*}) < r_t$ if $t\notin \cT$. The adversary's optimal value at $\bz$ becomes
\begin{align}
\cL(\bz,\ba^*,\bb^*) &= \sum_{n\in[N]} z_n \widetilde{p}^G(\bz,{\ba^*},{\bb^*})_n - \sum_{t\in\cT}\lambda_t(\bd^t)^\transpose  \widetilde{\bp}^G(\bz,{\ba^*},{\bb^*}) + \sum_{t\in\cT}\lambda_t r_t\nonumber\\
 &= \frac{\sum_{n\in[N]} \left(z_n - \sum_{t\in \cT}\lambda_t d^t_n\right)G^n(\bY^n|z_n,\ba^{n*},\bb^{n*}) }{1+\sum_n G^n(\bY^n|z_n,\ba^{n*},\bb^{n*})  } + \sum_{t\in\cT}\lambda_t r_t,\nonumber
\end{align}
 For notational brevity, let 
\[
\begin{aligned}
\rho^* &=  \frac{\sum_{n\in[N]} \left(z_n - \sum_{t\in \cT}\lambda_t d^t_n\right)G^n(\bY^n|z_n,\ba^{n*},\bb^{n*}) }{1+\sum_n G^n(\bY^n|z_n,\ba^{n*},\bb^{n*})  }\\ \cI_1 &= \{n\in [N]|\ \rho^* < z_n - \sum_{t\in \cT}\lambda_t d^t_n\}\\ 
\cI_2 &= \{n\in [N]|\ \rho^* >z_n - \sum_{t\in \cT}\lambda_t d^t_n \}\\
\cA^\bz &= \{(\ba,\bb)\in\cA|\ G^n(\bY^n|z_n,\ba^{n},\bb^{n}) = \underline{\cG}^n(z_n) \text{ if } n\in \cI_1,\ G^n(\bY^n|z_n,\ba^{n},\bb^{n}) = \underline{\cG}^n(z_n) \text{ if } n\in \cI_2 \}
\end{aligned}
\] 
From Lemma \ref{lemma:RO-GEV-hete-bounds}, if  $(\ba^*, \bb^*) \notin \cA^\bz$,  then for any $(\ba,\bb)\in\cA^\bz$ we have
\[
\begin{aligned}
\cL(\bz,\ba^*,\bb^*)&> \sum_n z_n \widetilde{p}_n(\bz,{\ba},{\bb}) - \sum_{t\in\cT}\lambda_t(\bd^t)^\transpose  \widetilde{\bp}^G(\bz,{\ba},{\bb}) + \sum_{t\in\cT}\lambda_t r_t \\
&\geq \cL(\bz,\ba,\bb),
\end{aligned}
\]
which is contradictory to the assumption that $(\ba^*,\bb^*)$ is optimal to the adversary's problem. So we have $(\ba^*, \bb^*) \in \cA^\bz$. On the other hand, Lemma \ref{lemma:RO-GEV-hete-bounds} tells us that if we take any point $(\ba,\bb)\in\cA^\bz$ such that $G^n(\bY^n|z_n,\ba,\bb) \in \Big\{\underline{\cG}^n(z_n), \overline{\cG}^n(z_n)\Big\}$ for all $n\notin \cI_1 \cup \cI_2$,  we also have
\[
\begin{aligned}
\cL(\bz,\ba^*,\bb^*)  &=  \sum_n \left(z_n-\sum_{t\in\cT}\lambda_t d^t_n\right)\widetilde{p}^G(\bz,{\ba},{\bb})_n  + \sum_{t\in\cT} \lambda_t r_t \\
&\geq  \sum_n \left(z_n\right)\widetilde{p}^G_n(\bz,{\ba},{\bb}) - \sum_{t=1}^T \lambda_t \max\{0, (\bd^t)^\transpose \widetilde{\bp}^G(\bz,{\ba},{\bb}) - r_t\} = 
\cL(\bz,\ba,\bb).  
\end{aligned}
\]
Since $\cL(\bz,\ba^*,\bb^*)= \min_{(\ba,\bb)\in\cA}\cL(\bz,\ba,\bb)$, we have $\cL(\bz,\ba^*,\bb^*) = \cL(\bz,\ba,\bb)$, meaning that $(\ba,\bb)$ is also optimal to the adversary's problem under prices $\bx$. This completes the proof. 
\endproof
The lemma above tells us that a optimal solution to the adversary's problem for which the adversary will force $G^n(\bY^n|z_n,\ba,\bb)$ to their minimum or maximum values. 
%This allows us to simplify the robust problem into those of the form
%\ctien{Write possible formulation here}. The number of possible problems to to consider is still huge ($2^{2m}$). 
The next lemma further characterizes an important property of the robust optimal prices, which states that, under a robustly optimal solution, the  adversary will always force $G^n(\bY^n|z_n,\ba,\bb)$ to their minimums. This claim is similar to the claim in Lemma \ref{lemma:RO-GEV-hete-solutions}. The proof is however more challenging due to the term $\sum_{t=1}^T \lambda_t \max\{0, (\bd^t)^\transpose \bp^G - r_t\}$.  

First, let $\bz^*$ be a robust optimal solution to the robust problem and $(\ba^*,\bb^*)$ be an optimal solution to the adversary problem such that $G^n(\bY^n|z^*_n,\ba^{n*},\bb^{n*}) \in \Big\{\underline{\cG}^n(z^*_n), \overline{\cG}^n(z^*_n)\Big\}$. We  also  denote by $\cT^*$ a subset of $[T]$ such that $(\bd^t)^\transpose  \widetilde{\bp}^G(\bz^*,\ba^{*}, \bb^{*}) \geq r_t$ for all $t\in\cT^*$ and $(\bd^t)^\transpose \widetilde{\bp}^G(\bz^*,\ba^{*}, \bb^{*}) < r_t$ if $t\notin \cT^*$, and let
\[
\rho^* =\frac{\sum_{n\in[N]} \left(z_n - \sum_{t\in \cT}\lambda_t d^t_n\right)G^n(\bY^n|z^*_n,\ba^{n*},\bb^{n*}) }{1+\sum_n G^n(\bY^n|z^*_n,\ba^{n*},\bb^{n*})  }.
\]
\begin{lemma}\label{lemma:rb-pen-lm2}
$\rho^* < z^*_n - \sum_{t\in \cT^*}\lambda_t d^t_i$, for all $n \in [N]$.
\end{lemma}
\proof{Proof:}
Let $k = \text{argmax}_{k\in[N]} z^*_k - \sum_{t\in \cT^*}\lambda_t d^t_k$.
By contradiction, assume that  $\rho^* \geq  z^*_k - \sum_{t\in \cT^*}\lambda_t d^t_k$.
According to Lemma \ref{lemma:rb-pen-lm1}  and to facilitate the exposition, we also parameterize the adversary's objective function by a vector $\bu\in\{0,1\}^N$ as
\[
\cL(\bz|\bu) = \psi(\bz|\bu) -  \sum_{t=1}^T \lambda_t \max\{0, (\bd^t)^\transpose \widetilde{\bp}^G(\bz|\bu) - r_t\} 
\]
where $\psi(\bz|\bu)$ is defined in \eqref{eq:define-psi}  and $\widetilde{\bp}^G(\bz|\bu)$ is defined as
\[
\widetilde{\bp}^G(\bz|\bu)_n = \frac{\theta^n(z_n)}{1+\sum_{n\in[N]} \theta^n (z_n)}, 
\]
where $\theta^n(z_n) = \underline{\cG}^n(z_n)$ if $u_n=0$ and $\theta^n(z_n) = \overline{\cG}^n(z_n)$ otherwise. 
We also let
\begin{align}
\delta &= \min_{\bu \in\{0,1\}^N} \left\{\cL(\bz^*|\bu) - f(\bz^*))\Big|\ \cL(\bz^*|\bu)>f(\bz^*)  \right\},\label{eq:rb-pen-def-delta}
\end{align}
with a note that we set $\delta = +\infty$  if the corresponding searching set is empty. Let $\bU^*$ be the set of parameter $\bu^*$ such that  $\cL(\bz^*|\bu^*) = f(\bz^*)$. 
For any $\bu^*\in\bU^*$ and any $\epsilon>0$ we easily  have $\widetilde{\bp}^G(\bz^*+\epsilon\bbe^k|\bu^*)_n > \widetilde{\bp}^G(\bz^*|\bu^*)_n$. Moreover, from Lemma \ref{lemma:psi-monotone} we have $\psi(\bz^*+ \epsilon \bbe^k|\bu^*) > \psi(\bz^*+ \epsilon \bbe^k|\bu^*)$, leading to the fact that, for any $\epsilon>0$ we have $\cL(\bz^*+\epsilon \bbe^k|\bu^*)> \cL(\bz^*|\bu^*)$.
Moreover, since  $\cL(\bz|\bu^*)$ and $f(\bz)$ are continuous in $\bz$, we always can select $\epsilon>0$ small enough such that
\begin{align}
\cL(\bz^*|\bu^*) <\cL(\bx+\epsilon \bbe^k|\bu^*)\label{eq:rb-penalties-cond1}\\
|f(\bz^*) - f(\bz^*+\epsilon \bbe^k)|<\delta/2\label{eq:rb-penalties-cond2}\\
\Big|\cL(\bz^*+\epsilon \bbe^k|\bu^\epsilon) - \cL(\bz^*|\bu^\epsilon) \Big| <\delta/2,\label{eq:rb-penalties-cond3}
\end{align}
where $\bu^\epsilon$ is a configuration of  the adversary's problem under prices $\bz+\epsilon \bbe^k$.
Applying the triangular inequality with  \eqref{eq:rb-penalties-cond2} and \eqref{eq:rb-penalties-cond3} we have
\begin{align}
 |f(\bz^*) -\cL(\bz^*|\bu^\epsilon)|\leq |f(\bz^*) - f(\bz^*+\epsilon\bbe^k)| + | f(\bz^*+\epsilon\bbe^k)- \cL(\bz^*|\bu^\epsilon)|<\delta.\nonumber  
\end{align}
So, according to the definition of $\delta$ in \eqref{eq:rb-pen-def-delta}, we have $f(\bz^*)  = \cL(\bz^*|\bu^\epsilon)$, meaning that $\bu^\epsilon \in  \bU^*$. So, we can always choose a configuration $\Bar{\bu} \in \bU^*$ such that $\Bar{\bu}$ is also a configuration for the adversary's problem under $\bz^*+\epsilon \bbe^k$, i.e., $\cL(\bz^*|\Bar{\bu}) = f(\bz^*)$ and $\cL(\bz^*+\epsilon\bbe^k|\Bar{\bu}) = f(\bz^*+\epsilon\bbe^k)$. 
Together with \eqref{eq:rb-penalties-cond1}, we have 
\[
f(\bz^*) = \cL(\bz^*,\Bar{\bu})<  \cL(\bz^*+\epsilon\bbe^k|\Bar{\bu}) = f(\bz^*+\epsilon\bbe^k),
\]
which is contradictory to our initial assumption that $\bz^*$ is a robust optimal solution. So our contradiction hypothesis is untrue and this completes the proof. 
\endproof
We are now ready to show the proof of Theorem \ref{theor:ro-pen-gev}. 
\proof{Proof of Theorem \ref{theor:ro-pen-gev}.}
From Lemma \ref{lemma:rb-pen-lm1} and \ref{lemma:rb-pen-lm2}, we have that if $\bz^*$ is a robust optimal solution, then the adversary will force all $G^n(\bY^n|z^*_n,\ba^{n},\bb^{n}) $ to their minimum values, i.e., $f(\bz^*) = \cL(\bz^*|\bu^0)$, where $\bu^0$ is a vector of size $N$ with all zero entries. Moreover, $f(\bz^*) < \cL(\bz^*|\bu)$ for any $\bu\in\{0,1\}^N$, $\bu\neq \bu^0$. We now need to prove that $\bz^*$ is also optimal to the maximization problem $\max_{\bx} \cL(\bz|\bu^0)$. In this case, the function $\cL(\bz|\bu^0)$ is not differentiable in $\bz$, so we cannot use the techniques in the proof of Theorem  \ref{theor:convert-simplified-model} above. Fortunately, if we consider the objective function $\cL(\bz|\bu^0)$ as a function of the purchase probabilities ${\bp}^G$, then we can show that this function is strictly concave in ${\bp}^G$. To facilitate this point, lets us define
\[
\cF(\bp^G|\bu^0) = \cL(\bz|\bu^0) = \bz(\bp^G|\bu^0)^\transpose \bp^G  - \sum_{t=1}^T \lambda_t \max\{0, (\bd^t)^\transpose \bp^G - r_t\}
\]
We know that the first term $\bz(\bp^G|\bu^0)^\transpose \bp^G$ is strictly concave in $\bp^G$  (Theorem \ref{theor;RO-hete-convexness-simplified-model}) and it is not difficult to show that $-\sum_{t=1}^T \lambda_t \max\{0, (\bd^t)^\transpose \bp^G - r_t\}$ is concave in $\bp^G$. As a result, $\cF(\bp^G|\bu^0)$ is strictly concave in $\bp^G$.

Now, let $\bp^{G*}$ be the purchase probabilities given by prices $\bz^*$. We will prove that $\bp^{G*} = \text{argmax}_{\bp^G} \cF(\bp^G|\bu^0)$. 
We omit $\bu^0$ for notational simplicity.
By contradiction, assume that  $\widetilde{\bp} = \text{argmax}_{\bp^G} \cF(\bp^G)$ and $\cF(\widetilde{\bp})> \cF({\bp^{G*}})$. Since $\cF(\bp^G)$ is strictly concave in $\bp^G$, we have, for any $t\in(0,1)$,
\[
t \cF(\widetilde{\bp})+ (1-t)\cF({\bp^{G*}}) <\cF(t\widetilde{\bp} + (1-t)\bp^{G*})
\]
Since $\cF(\widetilde{\bp}) \geq \cK(t\widetilde{\bp} + (1-t)\bp^{G*})$, we have $\cF({\bp^{G*}}) <\cF(t\widetilde{\bp} + (1-t)\bp^{G*})$ for all $t\in(0,1)$. This also mean that for any $\epsilon>0$, we always can find a point $\bp \in \cP^N$ such that $||\bp^{G*}-\bp||\leq \epsilon$ and $\cF({\bp^{G*}})<\cF({\bp})$.
Since $\bp^G(\bz|\bu^0)$ ($\bp^G$ as a function of $\bz$) is continuous  in $\bz$, this also means that given any $\epsilon>0$, there always exists $\bx\in \bbR^m$ such that $||\bx-\bx^*||\leq \epsilon$ and $\cL(\bx^*|\bu^0)<\cL(\bx|\bu^0)$. 

Now, similarly to the proof of Lemma \ref{lemma:rb-pen-lm2}, let
\begin{equation}
\label{eq:proof-theor-pen-eq1}
 \delta = \min_{\bu\in\{0,1\}^N}\left\{ \cL(\bz^*|\bu) - \cL(\bz^*|\bu^0)\Big|\ \cL(\bz^*|\bu) > \cL(\bz^*|\bu^0) \right\}
\end{equation}
Since $f(\bz)$ and $\cL(\bz|\bu)$ are continuous in $\bz$, there is an $\epsilon>0$ such that, for all $\bz\in\bbR^m$, $||\bz^*-\bz||\leq \epsilon$ 
\begin{equation}
 \begin{cases}
 |f(\bz^*) - f(\bz)| <\delta/2 \\
 |\cL(\bz^*|\bu^\bz)  - \cL(\bz|\bu^\bz)|<\delta/2,
 \end{cases}
\end{equation}
where $(\bu^\bz)$ is a configuration vector of the adversary's problem \eqref{prob:robust-constrained-penalties} under prices $\bz$.
As a result, for all $\bz$ such that $||\bz^*-\bz||\leq \epsilon$
\[
|\cL(\bz^*|\bu^0) - \cL(\bz^*|{\bu^\bz})|<\tau.
\]
Combine this with \eqref{eq:proof-theor-pen-eq1}, since
 $\bu^0$ is the unique configuration for the adversary's problem under prices $\bz^*$, we have that for all $\bz$ such that $||\bz^*-\bz||\leq \epsilon$, ${\bu^\bz} = \bu^0$. Moreover, we have shown that given any $\epsilon>0$, there exists $\overline{\bz}$ such that $||\overline{\bz}-\bz^*||\leq \epsilon$ and $\cL(\overline{\bz}|\bu^0)>\cL(\bz^*|\bu^0)$. So, if we choose $\epsilon>0$ and  small enough, we have 
 \[
 f(\Bar{\bz}) = \cL(\Bar{\bz}|\bu^{\Bar{\bz}}) = \cL(\Bar{\bz}|\bu^0) >\cL({\bz^*}|\bu^0) = f(\bz^*).
 \]
This is contradictory to the  assumption that $\bz^*$ is a robust optimal solution.
So,  our contradiction hypothesis that $\bp^{G*}$ is not optimal to $\max_{\bp^G} \cF(\bp^G|\bu^0)$ is untrue, meaning  $\bz^*$ is optimal to $\max_\bz \cL(\bz|\bu^0)$. Moreover, $\max_{\bp^G} \cF(\bp^G|\bu^0)$ always yields a unique solution as the objective function is strictly concave, 
$\max_\bz \cL(\bz|\bu^0)$ also yields a unique solution and this solution is also a unique robust  optimal solution to \eqref{prob:R-penalty-z}. 

We need a final step to complete the proof. We can easily see that  Problem \eqref{prob:RO-pen-GEV} can be converted equivalently as 
\begin{align}
	\underset{\bp^G,\by}{\text{max}}\qquad &\cW(\bz(\bp^G)) - \sum_{t=1}^T \lambda_t y_t & \nonumber \\
	 \text{subject to} \qquad & (\bd^t)^\transpose \bp^G - y_t \leq r_t &  \nonumber\\
	  &  \sum_{n\in [N]}p^G_n \leq 1&  \nonumber\\
	  &  \bp^G,\by\geq 0. &  \nonumber
\end{align}
which is a convex optimization problem, as $\cW(\bz(\bp^G))$ is strictly concave \eqref{theor;RO-hete-convexness-simplified-model}. 

%%%%%%%%%%%%%%%%%%%%%%%%%%%
%%%%%%%%%%%%%%%%%%%%%%%%%%%%

\subsubsection{Proof of Theorem \ref{thr:convergence-ro-pen}}
\label{proof:convergence-pen}

First, let us consider the deterministic version of the pricing problem with penalties, which can be formulated as the convex optimization problem
\begin{align}
	\underset{\bp,\by}{\text{max}}\qquad &  \sum_{i\in\cV}\left( \bx(\bp|\ba,\bb,G)_i-c_i\right)p_i - \sum_{t=1}^T \lambda_t y_t & \label{prob:DET-PENalties-inner-convexOpt} \\
	 \text{subject to} \qquad & (\bap^t)^\transpose \bp - y_t \leq r_t &  \nonumber\\
	  &  \sum_{i\in \cV}p_i \leq 1&  \nonumber\\
	  &  \bp,\by\geq 0. &  \nonumber
\end{align}
Before moving to a robust version, we  investigate some characteristics of the deterministic pricing problem with penalties \eqref{prob:DET-PENalties-inner-convexOpt}. First, let us denote by $v^*$ and $\bp^*$ the optimal value and optimal solution of the standard pricing problem under expected sale constraints  and $v^{\pmb{\lambda}}$ and $\bp^{\pmb{\lambda}}$ the optimal value and optimal solution to the pricing problem with over-expected-sale penalties \eqref{prob:DET-PENalties-inner-convexOpt}. Theorem \ref{theor:det-over-ES-penalties} below shows that the expected value given by \eqref{prob:DET-PENalties-inner-convexOpt} will converges to the optimal expected revenue given by the constrained pricing problem %\eqref{prob:det-nopenalties-inner-convexOpt}
when $\lambda_t$, $\forall t\in[T]$, increase to infinity. 
\begin{theorem}[Convergence of the optimal value when the penalty parameters increase]
\label{theor:det-over-ES-penalties}
For any $\epsilon>0$, we have
\begin{itemize}
 \item[(i)] For any  $\pmb{\lambda}^1, \pmb{\lambda}^2 \in\bbR^T_+$ such that $\pmb{\lambda}^1 - \pmb{\lambda}^2 = \epsilon \textbf{1}$, $$\sum_{t}\max\{0,(\bap^t)^\transpose \bp^{\pmb{\lambda}^1} - r_t\} \leq \sum_{t}\max\{0,(\bap^t)^\transpose \bp^{\pmb{\lambda}^2} - r_t\},$$
where $\textbf{1}$ is a unit vector of appropriate size.
 \item[(ii)] $v^{\pmb{\lambda}} \geq v^*$ for all $\pmb{\lambda}\in\bbR^T_+$ and if $\lambda_0 = \min_{t\in [T]}\lambda_t \geq (\Delta^* - v^*)/\epsilon$ then $\sum_{t}\max\{0,(\bap^t)^\transpose \bp^\bld - r_t\}\leq \epsilon$, where $\Delta^* = \max_{\bx}\Phi(\bx, \ba,\bb)$.
 \item[(iii)] Assume that there are positive constant $L_i,l_i$, $i\in \cV$ such that $Y_i\partial G_i(\bY)$ is bounded from above  by $L_iY_i^{li}$ for all prices $\bx \geq 0$, then for any $\epsilon$ such that 
 \[
 \epsilon \leq {\min_{t,i} \{\alpha^t_i|\ \alpha^t_i> 0\}}\min_t\left\{\frac{r_t}{(\bap^t)^\transpose \textbf{1}}\right\},
 \]
  then if we choose $\lambda_0\geq (\Delta^* - v^*)/\epsilon$, we can upper-bound $|v^\bld - v^*|$ as
  \[
  |v^\bld - v^*|\leq \max\left\{\max_i\left\{\frac{a_i}{b_i} - \frac{1}{b_il_i}\log \frac{\delta(\epsilon)}{L_i}\right\}, 0\right\} \frac{m\epsilon}{\min_{t,i} \{\alpha^t_i|\ \alpha^t_i> 0\}},
  \]
  where $\delta(\epsilon) = \min_t\left\{\frac{r_t}{(\bap^t)^\transpose \textbf{1}}\right\} - \frac{\epsilon}{\min_{t,i} \{\alpha^t_i|\ \alpha^t_i> 0\}}$, and this upper bound converges to zero linearly when $\epsilon$ tends to zero. 
\end{itemize}

\end{theorem}
\proof{Proof:}
First, for notational simplicity we denote $R(\bp) = \sum_{i\in\cV}\left( \bx(\bp|\ba,\bb,G)_i-c_i\right)p_i$.
For (i), we have the following inequalities 
\begin{align}
 R(\bp^{\bld^1})  &- \sum_{t} \lambda^1_t \max\{0,(\bap^t)^\transpose \bp^{\pmb{\lambda}^1} - r_t\} \geq  R(\bp^{\bld^2})  - \sum_{t}\lambda^1_t\max\{0,(\bap^t)^\transpose \bp^{\pmb{\lambda}^2} - r_t\}\nonumber \\
 & =  R(\bp^{\bld^2})  - \sum_{t}\lambda^2_t\max\{0,(\bap^t)^\transpose \bp^{\pmb{\lambda}^2} - r_t\} -\epsilon \sum_{t}\max\{0,(\bap^t)^\transpose \bp^{\pmb{\lambda}^2} - r_t\} \nonumber \\
 & \geq  R(\bp^{\bld^1})  - \sum_{t}\lambda^2_t\max\{0,(\bap^t)^\transpose \bp^{\pmb{\lambda}^1} - r_t\} -\epsilon \sum_{t}\max\{0,(\bap^t)^\transpose \bp^{\pmb{\lambda}^2} - r_t\} \nonumber
\end{align}
So, we have
\[
\epsilon \sum_{t}\max\{0,(\bap^t)^\transpose \bp^{\pmb{\lambda}^2} - r_t\} \geq \sum_t (\lambda^1_t-\lambda^2_t)\max\{0,(\bap^t)^\transpose \bp^{\pmb{\lambda}^1} - r_t\},
\]
which leads to the desired inequality.

For (ii), since $(\bap^t)^\transpose \bp^* \leq  r_t$ for all $t$,  given $\bld\in\bbR^T_+$, we have 
\[
v^{\bld} \geq R(\bp^*) - \sum_{t}\lambda_t \max\{0,(\bap^t)^\transpose \bp^* - r_t\} = R(\bp^*) = v^*.
\]
Moreover, since
$v^{\bld} - v^* = R(\bp^{\bld})- v^* - \sum_{t}\lambda_t \max\{0,(\bap^t)^\transpose \bp^{\bld} - r_t\}$, we have
\begin{equation}\label{eq:det-over-ES-penalties-eq1}
R(\bp^{\bld})- v^* \geq \sum_{t}\lambda_t \max\{0,(\bap^t)^\transpose \bp^{\bld} - r_t\}\geq \lambda_0 \sum_{t} \max\{0,(\bap^t)^\transpose \bp^{\bld} - r_t\} 
\end{equation}
The left hand side of \eqref{eq:det-over-ES-penalties-eq1} is less  than $\Delta^*-v^*$, so if we choose $\lambda_0\geq (\Delta^*-v^*)/\epsilon$ then $\sum_{t} \max\{0,(\bap^t)^\transpose \bp^{\bld} - r_t\} \leq \epsilon$ as desired.

We move to (iii). 
As shown previously, we  can choose $\lambda_0$ such that  $\max\{0,(\bap^t)^\transpose \bp^{\bld} - r_t\} \leq \epsilon$  or $(\bap^t)^\transpose \bp^{\bld} \leq r_t+ \epsilon$ for all $t\in[T]$. We now consider the following problem 
\begin{equation}\label{prob:det-over-ES-penalties-prob1}
\max_{\substack{\bp\geq 0 \\\sum_i p_i\leq 1}}\left\{R(\bp)\Big| \ (\bap^t)^\transpose \bp \leq r_t+ \epsilon,\ \forall t \right\}
\end{equation} 
and denote by ${\bp}^\epsilon$ as an optimal solution to \eqref{prob:det-over-ES-penalties-prob1}. Since $\bp^{\bld}$ is feasible to \eqref{prob:det-over-ES-penalties-prob1} we have $R({\bp}^\epsilon) \geq R(\bp^{\bld})$. Moreover, if we define $\cP := \{\bp\in\bbR^m|\ p_i\geq 0, \sum_i p_i\leq 1,\ (\bap^t)^\transpose \bp^{\bld} \leq r_t,\forall t\in [T]\}$, then $v^*\geq R(\bp)$ for all $\bp\in\cP$.  
Therefore, we have
\begin{equation}\label{eq:det-over-ES-penalties-eq2}
 |v^\bld - v^*| \leq R({\bp}^\epsilon) - R(\bp),\ \forall \bp \in \cP.
\end{equation}
We will show that there is $\bp\in \cP$ such that $||\bp^\epsilon - \bp||$ can be arbitrarily small when $\epsilon$ decreases, which allows us to use the \textit{Mean Value Theorem} to bound $|R({\bp}^\epsilon) - R(\bp)|$. 
If $\bp^\epsilon\in \cP$, then the result is obvious and we have $|v^\bld - v^*| = 0$. Now assume that $\bp^\epsilon\notin \cP$, let $\cT:= \{t\in[T]|\ (\bap^t)^\transpose \bp^\epsilon > r_t\}$ and for any $t\in\cT$ we select $i_t = \text{argmax}_{i\in\cV} \{ p^\epsilon_i|\ \alpha^t_i>0 \}$. Then we denote $\cI = \{i_t|\ t\in \cT\}$. We pick a $\widetilde{\bp}$ such that 
\begin{equation}\label{eq:det-over-ES-penalties-eq3}
\begin{cases}
\widetilde{p}_{i} = p^\epsilon_{i} - {\epsilon}/({\min_{t,j} \{\alpha^t_j|\ \alpha^t_j> 0\}}),\ \forall i\in\cI \\
\widetilde{p}_{j}  = p^\epsilon_{j},\ \forall j\notin \cI. 
\end{cases}
\end{equation}
With this selection, we see that, for any $t\in\cT$ 
\begin{equation}\label{eq:proof-th4-eq2}
\begin{aligned}
(\bap^t)^\transpose \widetilde{\bp} &\leq (\bap^t)^\transpose {\bp}^\epsilon - \alpha^t_{i_t} \epsilon/({\min_{t,i} \{\alpha^t_i|\ \alpha^t_i> 0\}}) \\
&\leq  (\bap^t)^\transpose {\bp}^\epsilon - \epsilon \leq r_t.
\end{aligned}
\end{equation}
And indeed for any $t\notin\cT$ we have $(\bap^t)^\transpose \widetilde{\bp} \leq (\bap^t)^\transpose {\bp^\epsilon}\leq r_t$.
Furthermore, for any $t\in\cT$, we have  $(\bap^{t})^\transpose \bp^\epsilon>r_{t}$. Combine this with the fact that  $i_t = \text{argmax}_{i\in\cV} \{ p^\epsilon_i|\ \alpha^t_i>0 \}$ we have 
\[
\left(\sum_{i}\alpha^t_i \right)p^{\epsilon}_{t_t} \geq (\bap^{t})^\transpose \bp^\epsilon > r_t.
\]
So, under the assumption on the selection of $\epsilon$, we have the chain of inequalities
\[
p^\epsilon_{i_t} >\frac{r_{t}}{(\bap^{t})^\transpose \textbf{1}}\geq \min_t\left\{\frac{r_t}{(\bap^t)^\transpose \textbf{1}}\right\} \geq \frac{\epsilon}{\min_{t,i} \{\alpha^t_i|\ \alpha^t_i> 0\}},
\]
meaning that $\widetilde{\bp}> 0$. So, combine with \eqref{eq:proof-th4-eq2} we have $\widetilde{\bp}\in\cP$. 

Moreover, for any point $\bp' \in [\widetilde{\bp},\bp^\epsilon]$ and any $i\in \cI$, we have 
\begin{align}
p'_{i} \geq \widetilde{p}_{i} &= p^\epsilon_i - \frac{\epsilon}{\min_{t,i} \{\alpha^t_i|\ \alpha^t_i> 0\}} \nonumber\\
&> \min_t\left\{\frac{r_t}{(\bap^t)^\transpose \textbf{1}}\right\} - \frac{\epsilon}{\min_{t,i} \{\alpha^t_i|\ \alpha^t_i> 0\}} := \delta(\epsilon). 
\end{align}
So, if we  denote $\bx' = \bx(\bp',G)$ (i.e., the prices that result in purchase probabilities $\bp'$). 
For any $i\in\cI$, under the assumption that $Y_i\partial G_i(\bY)\leq L_iY_i^{l_i}$ we have
\[
\begin{aligned}
L_{i}Y_{i}(\bx')^{l_{i}}& \geq  p'_i (1+G(\bY(\bx'))) \\
&\geq p'_i \geq \delta(\epsilon),
\end{aligned}
\]
where $\bY(\bx')$ is a vector of size $m$ with entries $Y_j(\bx') = \exp(a_j-b_jx'_j)$, $\forall j\in\cV$. So we have
\[
x'_i \leq \frac{a_i}{b_i} - \frac{1}{b_il_i}\log \frac{\delta(\epsilon)}{L_i}.
\]
Moreover, if we look at the gradient of $R(\bp)$ at $p'_{i}$.  According to Theorem 4.3 in \cite{zhang2018multiproduct} we have
\begin{equation}\label{eq:det-over-ES-penalties-eq4}
 \nabla_{\bp} R(\bp')_{i} \leq x'_{i}  \leq \max_j\left\{\frac{a_j}{b_j} - \frac{1}{b_jl_j}\log \frac{\delta(\epsilon)}{L_j}\right\}.
\end{equation}
Now, we look at $|R(\bp^\epsilon) - R(\widetilde{\bp})|$ and by combining \eqref{eq:det-over-ES-penalties-eq3}, \eqref{eq:det-over-ES-penalties-eq4}, the \textit{Mean Value Theorem} tells us that there is $\bp'\in[\widetilde{\bp},\bp^\epsilon]$
\begin{align}
|R(\bp^\epsilon) - R(\widetilde{\bp})| &= \sum_{i\in \cI} \nabla_\bp R(\bp')_{i} |\widetilde{p}_{i} - p^\epsilon_{i}|\nonumber\\
&\leq \max_i\left\{\frac{a_i}{b_i} - \frac{1}{b_il_i}\log \frac{\delta(\epsilon)}{L_i}\right\} \frac{m\epsilon}{\min_{t,i} \{\alpha^t_i|\ \alpha^t_i> 0\}}\label{eq:det-over-ES-penalties-eq5}
\end{align}
Combine \eqref{eq:det-over-ES-penalties-eq5} with \eqref{eq:det-over-ES-penalties-eq2} and recall that $\widetilde{\bp}\in\cP$, we have
\[
|v^\bld - v^*|\leq \max_i\left\{\frac{a_i}{b_i} - \frac{1}{b_il_i}\log \frac{\delta(\epsilon)}{L_i}\right\} \frac{m\epsilon}{\min_{t,i} \{\alpha^t_i|\ \alpha^t_i> 0\}}.
\]
Combine with the case $\bp^\epsilon \in\cP$, we obtain the desired bound, which definitely converge to zero when $\epsilon$ tends to zero, as desired.\hfill Q.E.D.
\endproof

Now we are ready for the main proof. 

\proof{Proof of Theorem \ref{thr:convergence-ro-pen}:}
Using a similar evaluation as in \eqref{eq:det-over-ES-penalties-eq1} we can have
\[
\cW(\bz(\bp^{G,\bld})) - \overline{\varphi}  \geq \sum_{t}\lambda_t \max\{0,(\bd^t)^\transpose \bp^{G,\bld} - r_t\}\geq \lambda_0 \sum_{t} \max\{0,(\bd^t)^\transpose \bp^{G,\bld} - r_t\} 
\]
and we also wee that the left hand side of the above is less than $\Delta^*-\overline{\varphi}$. Thus, if we choose $\lambda_0 \geq (\Delta^*-\overline{\varphi})/\epsilon$, the we have $\sum_{t}\lambda_t \max\{0,(\bd^t)^\transpose \bp^{G,\bld} - r_t\} \leq \epsilon$.

To prove the second claim of the corollary, we can choose $\lambda_0$ such that $\sum_{t}\lambda_t \max\{0,(\bd^t)^\transpose \bp^{G,\bld} - r_t\} \leq \epsilon$ and consider $\bp^{G,\epsilon}$ as an optimal solution to the following problem
\begin{equation}\label{prob:ro-over-ES-penalties-prob1}
\max_{\substack{\bp^{G}\in\cP^G}}\left\{\cW(\bz(\bp^G))\Big| \ (\bd^t)^\transpose \bp^G \leq r_t+ \epsilon,\ \forall t \right\}.
\end{equation} 
 Let us define $\widetilde{\cP}^G = \{\bp^G\in\cP^G|(\bd^t)^\transpose \bp^G \leq r_t,\ \forall t\}$. Then, we have 
 \begin{equation}
 \label{eq:ro-proof123}
|\varphi^{\RO,\bld} - \overline{\varphi}| \leq \cW(\bz(\bp^{G,\epsilon})) -  \cW(\bz(\bp^G)),\;\forall \bp^G\in \widetilde{\cP}^G.
 \end{equation} 
We now try to bound $\cW(\bz(\bp^{G,\epsilon})) -  \cW(\bz(\bp^G))$ using the \textit{Mean Value Theorem}.  We see that if $\bp^{G,\epsilon}\in \widetilde{\cP}^G$ then $\varphi^{\RO,\bld} = \overline{\varphi}$. Otherwise,  assume that $\bp^{G,\epsilon}\notin \widetilde{\cP}^G$, let $\cT:= \{t\in[T]|\ (\bd^t)^\transpose \bp^{G,\epsilon} > r_t\}$ and for any $t\in\cT$ we select $n_t = \text{argmax}_{n\in [N]} \{ p^{G,\epsilon}_n|\ d^t_n>0 \}$. We also denote $\cI = \{n_t|\ t\in \cT\}$. We pick a vector  $\widetilde{\bp}^G$ such that 
\begin{equation}\label{eq:ro-over-ES-penalties-eq3}
\begin{cases}
\widetilde{p}^G_{n} = p^{G,\epsilon}_{n} - {\epsilon}/({\min_{t,k\in[N]} \{d^t_k|\ \alpha^t_k> 0\}}),\ \forall k\in\cI \\
\widetilde{p}^G_{k}  = p^{G,\epsilon}_{k},\ \forall k\notin \cI. 
\end{cases}
\end{equation}
Then  for any $t\in\cT$ we have 
\begin{equation}\label{eq:proof-th4-eq2-ro}
\begin{aligned}
(\bd^t)^\transpose \widetilde{\bp}^G &\leq (\bd^t)^\transpose {\bp}^{G,\epsilon} - d^t_{n_t} \epsilon/({\min_{t,k} \{d^t_k|\ d^t_k> 0\}}) \\
&\leq  (\bd^t)^\transpose {\bp}^{G,\epsilon} - \epsilon \leq r_t.
\end{aligned}
\end{equation}
Now, for any $t\notin\cT$ we have $(\bd^t)^\transpose \widetilde{\bp}^G \leq (\bd^t)^\transpose {\bp^{G,\epsilon}}\leq r_t$.
Furthermore, for any $t\in\cT$, we have  $(\bd^{t})^\transpose \bp^{G,\epsilon}>r_{t}$. Combine this with the selection of $n_t$ as $n_t = \text{argmax}_{n\in[N]} \{ p^{G,\epsilon}_n|\ d^t_n>0 \}$ we have 
\[
\left(\sum_{n}d^t_n \right)p^{G,\epsilon}_{n_t} \geq (\bd^{t})^\transpose \bp^{G,\epsilon} > r_t.
\]
So, from the selection of $\epsilon$, we have 
\[
p^\epsilon_{n_t} >\frac{r_{t}}{(\bd^{t})^\transpose \textbf{1}}\geq \min_t\left\{\frac{r_t}{(\bd^t)^\transpose \textbf{1}}\right\} \geq \frac{\epsilon}{\min_{t,n} \{d^t_n|\ d^t_n> 0\}}.
\]
Thus, 
 $\widetilde{\bp}^G\in \widetilde{\cP}^G$. Moreover, for any 
for any point $\bp' \in [\widetilde{\bp}^G,\bp^{G,\epsilon}]$ and any $n\in \cI$, we have 
\begin{align}
p'_{n} \geq \widetilde{p}_{n} &= p^{G,\epsilon}_n - \frac{\epsilon}{\min_{t,n\in[N]} \{d^t_n|\ d^t_n> 0\}} \geq \min_t\left\{\frac{r_t}{(\bd^t)^\transpose \textbf{1}}\right\} - \frac{\epsilon}{\min_{t,n} \{d^t_n|\ d^t_n> 0\}} := \delta^G(\epsilon). \nonumber
\end{align}
Hence, if we denote $\bz' =\bz(\bp')$ (i.e., the prices that result in the purchase probabilities $\bp'$), then for any $n\in\cI$, under the assumption that $Y_i\partial G_i(\bY)\leq L_iY_i^{l_i}$ we have $$G^n(\bY^n)  = \sum_{i\in \cV_n}Y_i\partial G_i(\bY) \leq \sum_{i\in \cV_n} L_iY_i^{l_i} \leq |\cV_n|L_{i_n}Y_{i_n}^{l_{i_n}}, $$
where $i_n\in \cV_n$ is chosen such that $L_{i_n}Y_i^{l_{i_n}}  = \max_{i\in \cV_n} L_iY_i^{l_i}$. We  have
\[
\begin{aligned}
L_{i_n}Y_{i_n}(\bz')^{l_{i_n}}& \geq  p'_i /|\cV_n| \\
&\geq  \delta(\epsilon)/|\cV_n|,
\end{aligned}
\]
Moreover, we can write  $Y(\bz')_{i_n} = \exp(a_{i_n}-b_{i_n}z'_n - b_{i_n}c_{i_n})$, for a given $(\ba,\bb\in \cA)$. So we have
\[
z'_n \leq  \max_{(\ba,\bb) \in\cA} \max_{i\in \cV_n}\left\{\frac{a_i - b_ic_i}{b_i} -  \frac{1}{b_il_i |\cV_n|}\log \frac{\delta(\epsilon)}{L_i} \right\}.
\]
Moreover, looking at the gradient of $\cW(\bz(\bp))$ at $p'_{n}$, we also see that there is a vector of parameters $(\ba,\bb)\in\cA$ such that 
\begin{equation}\label{eq:ro-over-ES-penalties-eq4}
 \nabla_{\bp} \cW(\bz(\bp'))_{n} \leq z'_{n}  \leq \max_{(\ba,\bb) \in\cA} \max_{i\in \cV_n}\left\{\frac{a_i - b_ic_i}{b_i} -  \frac{1}{b_il_i |\cV_n|}\log \frac{\delta(\epsilon)}{L_i} \right\}.
\end{equation}
Now, using the \textit{Mean Value Theorem}, there is $\bp'\in[\widetilde{\bp}^G,\bp^{G,\epsilon}]$ such that 
\begin{align}
|\cW(\bz(\bp^{G,\epsilon})) - \cW(\bz(\widetilde{\bp}^G))| &= \sum_{n\in \cI} \nabla_{\bp} \cW(\bz(\bp'))_{n} |\widetilde{p}^G_{n} - p^{G,\epsilon}_{n}|\nonumber\\
&\leq \max_{(\ba,\bb) \in\cA} \max_{n\in [N], i\in \cV_n}\left\{\frac{a_i - b_ic_i}{b_i} -  \frac{1}{b_il_i |\cV_n|}\log \frac{\delta(\epsilon)}{L_i} \right\}\frac{N\epsilon}{\min_{t,n\in[N]} \{d^t_n|\ d^t_n> 0\}}\label{eq:ro-over-ES-penalties-eq5}
\end{align}
Combine \eqref{eq:ro-over-ES-penalties-eq5} with \eqref{eq:ro-proof123} and recall that $\widetilde{\bp}\in\widetilde{\cP}^G$, we have
\begin{equation}
\label{eq:ro-proof-2147}
 |\varphi^{\RO,\bld} - \overline{\varphi}| \leq \max_{(\ba,\bb) \in\cA} \max_{n\in[N],i\in \cV_n}\left\{\frac{a_i - b_ic_i}{b_i} -  \frac{1}{b_il_i |\cV_n|}\log \frac{\delta(\epsilon)}{L_i} \right\}\frac{N\epsilon}{\min_{t,n\in[N]} \{d^t_n|\ d^t_n> 0\}}
\end{equation}
Combine with the case $\bp^{G,\epsilon} \in \widetilde{\cP}^G$, we obtain the desired bound. Since $\cA$ is bounded, the left hand side of \eqref{eq:ro-proof-2147} will always converges to zero when $\epsilon$ tends to zero, as desired.
\endproof
%%%%%%%%%%%%%%%%%%%%%%%%%%%%%%%%%%%%%%%
%
% ROBUST CONSTRAINED
%
%%%%%%%%%%%%%%%%%%%%%%%%%%%%%%%%%%%%%%

%%%%%%%%%%%%%%%%%%%%%%%%%%%%
%
% ROBUST CONSTRAINED PRICING
%%%%%%%%%%%%%%%%%%%%%%%%%%%%

\subsection{Experiments}
\label{apd:results}

%We first provide  experiments to show how the constrained pricing model performs when the choice parameters are uncertain, and how the deterministic pricing model with over-expected-sale penalties works, as compared to the constrained counterpart.
Our goal here is to illustrate  how the robust model with over-expected-sale penalties performs, as compared to other baseline approaches, i.e., deterministic and sampling-based counterparts.
We employ the same nested logit model with partition-wise homogeneous PSP considered above. 
We create one expected sale constraint (i.e., $T=1$) in such a way that the optimal prices from the unconstrained problem do not satisfy the expected sale constraint.  
We solve the deterministic problem with  the weighted average parameters   $\widetilde{\bw} = \sum_{k\in [K]} \tau_k \bw^k$ to obtain a solution $\bx^{\DET}$. 
Then, for each uncertainty level $\epsilon>0$ we solve the RO problem  by convex optimization to obtain a robust solution $\bx^\RO$. For the sampling-based approach, we also sample $10$ and $50$ points from the uncertainty set to get solutions $\bx^{\textsc{SA10}}$ and $\bx^{\textsc{SA50}}$, respectively.
We do not select a large sample size for the sampling-based approach due to the fact that the number of points $s_1$ is also the number of convex optimization problems to be solved, and these optimization problems, even-though computationally tractable, are still expensive to be done.

To evaluate the performance of the solutions obtained, similarly to the other cases, we sample randomly and uniformly 1000 points from $\cA$ and compute the corresponding expected revenues given by the four solutions  $\bx^{\textsc{DET}}$, $\bx^{\textsc{SA10}}, \bx^{\textsc{SA50}}$ and $\bx^{\RO}$.  
The distributions of the profit values (the expected revenue minus the penalty cost) for different $\lambda$ and $\epsilon$ are plotted in Figure  \ref{fig:Pen-hist}, where similar observations apply.
The histograms given by $\bx^\RO$ always have higher peaks, smaller variances, shorter tails and get tighter as $\epsilon$ increases,  as compared to the other solutions. 
 The histograms given by $\bx^{\textsc{SA50}}$ are quite similar to those from  $\bx^{\textsc{SA10}}$ and  also have higher peaks and smaller variances, as compared to those from $\bx^{\textsc{DET}}$
%This seems however not the case for $\lambda = 0.6$, especially when $\lambda=0.6$ and $\epsilon = 1.5$, the SA10 gives remarkably low profit values as compared to the other approach.
%More detailed results are provided in Table \ref{tb:RO-PEN-comparison}.
In general, we also see that the RO approach always gives higher worst-case  but lower average profits. 
The SA10 and SA50 also provide some protections against worst-case scenarios. The protection becomes better when the same size increases, which is rational given the fact that the minimax equality holds, thus a solution to the SA will approach a robust solution when the same size grows. %The percentile ranks of the RO worst-case profits are  small for the SA50 approaches, which indicates the conservativeness of the RO approach. However, as being said, the computational cost required  to perform the SA is proportional to the same size and could be very expensive when the sample size is large.   

\begin{figure}[htb]
 \centering
 \includegraphics[width=1\textwidth]{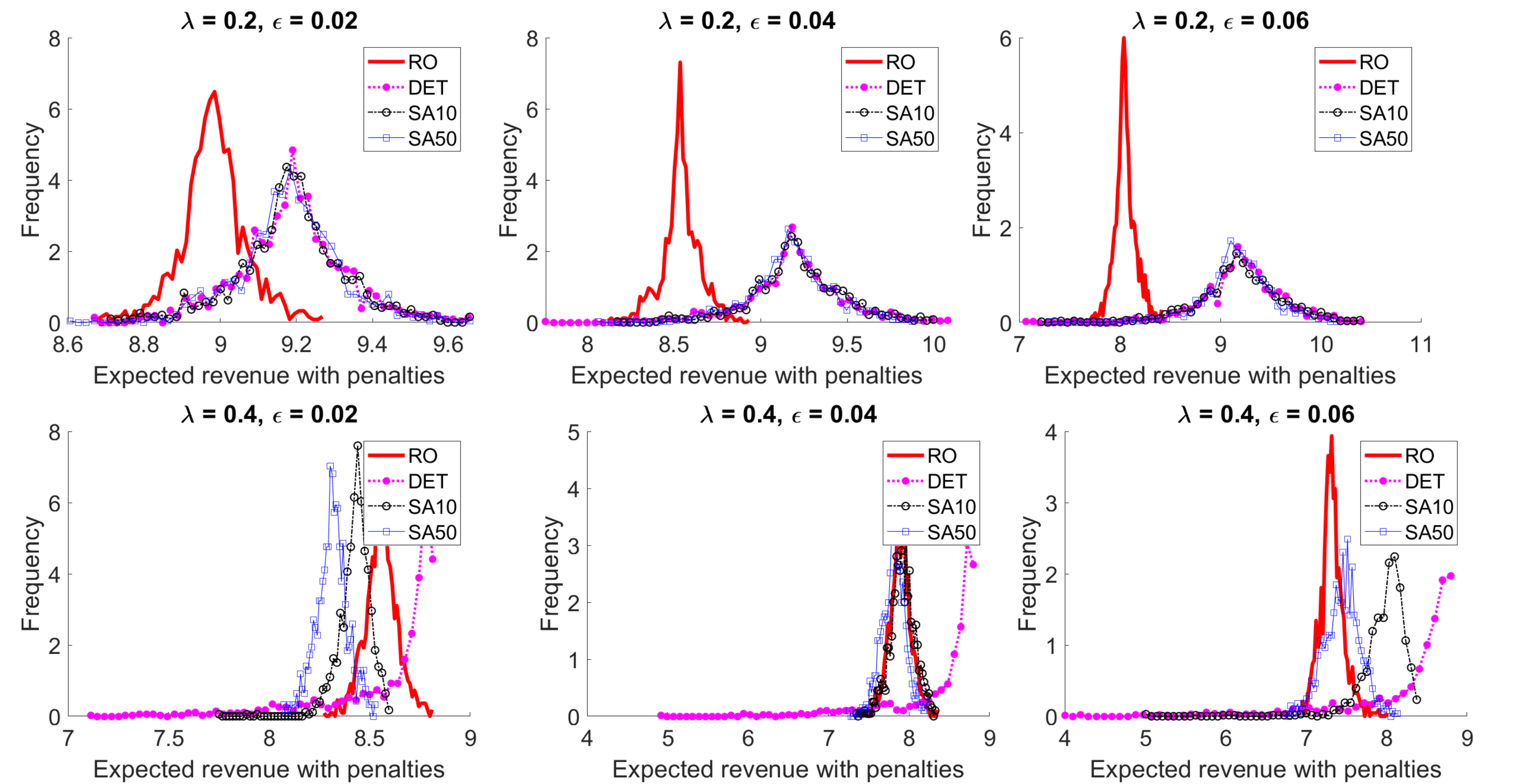}
 \caption{Distributions of the profit values under over-expected-sale penalties given by $\bx^\RO$, $\bx^{\textsc{DET}}$, $\bx^\textsc{SA10}$,  $\bx^{\textsc{SA50}}$.}
 \label{fig:Pen-hist}
\end{figure}

 \end{APPENDICES}

%% Here starts the e-companion (EC)
%%%%%%%%%%%%%%%%%%%%%%%%%%%%%%%%%%%%%%%%%%%%%%%%%%%%%%%%%%
%\ECSwitch
%\ECHead{Proofs of Statements}
%\section{...}
%\ECDisclaimer
%%%%%%%%%%%%%%%%%%%%%%%%%%%%%%%%%%%%%%%%%%%%%%%%%%%%%%%%%%

%%% Main head for the e-companion

%%%%%%%%%%%%%%%%%
\end{document}